\documentstyle[12pt]{article}

\catcode`\@=11
\def\marginnote#1{}

\newcount\hour
\newcount\minute
\newtoks\amorpm
\hour=\time\divide\hour by60
\minute=\time{\multiply\hour by60 \global\advance\minute by-\hour}
\edef\standardtime{{\ifnum\hour<12 \global\amorpm={am}%
        \else\global\amorpm={pm}\advance\hour by-12 \fi
        \ifnum\hour=0 \hour=12 \fi
        \number\hour:\ifnum\minute<10 0\fi\number\minute\the\amorpm}}
\edef\militarytime{\number\hour:\ifnum\minute<10 0\fi\number\minute}

\def\draftlabel#1{{\@bsphack\if@filesw {\let\thepage\relax
   \xdef\@gtempa{\write\@auxout{\string
      \newlabel{#1}{{\@currentlabel}{\thepage}}}}}\@gtempa
   \if@nobreak \ifvmode\nobreak\fi\fi\fi\@esphack}
        \gdef\@eqnlabel{#1}}
\def\@eqnlabel{}
\def\@vacuum{}
\def\draftmarginnote#1{\marginpar{\raggedright\scriptsize\tt#1}}

\def\draft{\oddsidemargin -.5truein
        \def\@oddfoot{\sl preliminary draft \hfil
        \rm\thepage\hfil\sl\today\quad\militarytime}
        \let\@evenfoot\@oddfoot \overfullrule 3pt
        \let\label=\draftlabel
        \let\marginnote=\draftmarginnote
   \def\@eqnnum{{\rm (\theequation)}\rlap{\kern\marginparsep\tt\@eqnlabel}%
\global\let\@eqnlabel\@vacuum}  }

\def\l@refs#1#2{\ifnum \c@tocdepth >\z@
    \addpenalty\@secpenalty
    \addvspace{1.0em \@plus\p@}%
    \setlength\@tempdima{1.5em}%
    \begingroup
      \parindent \z@ \rightskip \@pnumwidth
      \parfillskip -\@pnumwidth
      \leavevmode \bfseries
      \advance\leftskip\@tempdima
      \hskip -\leftskip
      #1\nobreak\hfil \nobreak\hb@xt@\@pnumwidth{\hss #2}\par
    \endgroup
  \fi}

\newcounter{app}

\def\app{\setcounter{equation}{0}
\def\theequation{A\arabic{app}.\arabic{equation}}\par
   \addvspace{4ex}
   \@afterindentfalse
  \secdef\@app\@dapp}
\newcommand\@app{\@startsection {app}{1}{0ex}%
                                   {-3.5ex \@plus -1ex \@minus -.2ex}%
                                   {2.3ex \@plus.2ex}%
                                   {\normalfont\Large\bf}}

\def\@dapp#1{%
{\parindent \z@ \raggedright  \bf #1}\par\nobreak}
\def\l@app#1#2{\ifnum \c@tocdepth >\z@
    \addpenalty\@secpenalty
    \addvspace{1.0em \@plus\p@}%
    \setlength\@tempdima{8.5em}%
    \begingroup
      \parindent \z@ \rightskip \@pnumwidth
      \parfillskip -\@pnumwidth
      \leavevmode \bfseries
      \advance\leftskip\@tempdima
      \hskip -\leftskip
      #1\nobreak\hfil \nobreak\hb@xt@\@pnumwidth{\hss #2}\par
    \endgroup\fi}
\newcounter{sapp}[app]

\def\sapp{\def\theequation{A\arabic{app}.\arabic{equation}}\par
   \@afterindentfalse
  \secdef\@sapp\@dsapp}
\newcommand\@sapp{\@startsection{sapp}{2}{\z@}%
                                     {-3.25ex\@plus -1ex \@minus -.2ex}%
                                     {1.5ex \@plus .2ex}%
                                     {\normalfont\large\bfseries}}

\def\@dsapp#1{%
{\parindent \z@ \raggedright  \bf #1}\par\nobreak}
\newcommand{\l@sapp}{\@dottedtocline{2}{1.5em}{3em}}

\def \b {\beta}
\def \d {\delta}

\def \tr{{\rm tr}}

\def \det{{\rm det}}

\def \log {{\rm log}}

\def \be {begin {equation}}
\def \ee {end{equation}}
\def\iw{intertwining operator }
\def\iws{intertwining operators }
\def\tilde{\widetilde}
\def\bar{\overline}
\def\hat{\widehat}
\def\*{\star}
\def\({\left(}		
\def\){\right)}		
\def\[{\left[}		
\def\]{\right]}

\def\frac#1#2{{#1 \over #2}}

\def\d{\partial}

\def\2pi{\hbox{$2\pi i$}}

\def\dsl{\raise.15ex\hbox{/}\kern-.57em\partial}
\def\Dsl{\,\raise.15ex\hbox{/}\mkern-.13.5mu D}
\def\be{\beta}
\def\al{\alpha}

\def\de{\delta}		\def\De{\Delta}
\def\om{\omega}

\font\numbers=cmss12
\font\upright=cmu10 scaled\magstep1
\def\stroke{\vrule height8pt width0.4pt depth-0.1pt}
\def\topfleck{\vrule height8pt width0.5pt depth-5.9pt}
\def\botfleck{\vrule height2pt width0.5pt depth0.1pt}
\def\Zmath{\vcenter{\hbox{\numbers\rlap{\rlap{Z}\kern 0.8pt\topfleck}\kern
2.2pt
                   \rlap Z\kern 6pt\botfleck\kern 1pt}}}
\def\Qmath{\vcenter{\hbox{\upright\rlap{\rlap{Q}\kern
                   3.8pt\stroke}\phantom{Q}}}}
\def\Nmath{\vcenter{\hbox{\upright\rlap{I}\kern 1.7pt N}}}
\def\Cmath{\vcenter{\hbox{\upright\rlap{\rlap{C}\kern
                   3.8pt\stroke}\phantom{C}}}}
\def\Rmath{\vcenter{\hbox{\upright\rlap{I}\kern 1.7pt R}}}
\def\Z{\ifmmode\Zmath\else$\Zmath$\fi}
\def\Q{\ifmmode\Qmath\else$\Qmath$\fi}
\def\N{\ifmmode\Nmath\else$\Nmath$\fi}
\def\C{\ifmmode\Cmath\else$\Cmath$\fi}
\def\R{\ifmmode\Rmath\else$\Rmath$\fi}

\def\be{ \begin{eqnarray} }
\def\ee{ \end{eqnarray} }

\def\d{\partial}
\def\p{\partial\over\partial}

\def\bea{\begin{eqnarray}}
\def\eea{\end{eqnarray}}
\def\nn{\nonumber}
\def\beq{\begin{equation}}
\def\eeq{\end{equation}}
\def\ba{\beq\new\begin{array}{c}}
\def\ea{\end{array}\eeq}
\def\be{\ba}
\def\ee{\ea}
\def\stackreb#1#2{\mathrel{\mathop{#2}\limits_{#1}}}
\def\Tr{{\rm Tr}}

\def\f{1\over}

\makeatletter
\newdimen\normalarrayskip              
\newdimen\minarrayskip                 
\normalarrayskip\baselineskip
\minarrayskip\jot
\newif\ifold             \oldtrue            \def\new{\oldfalse}
\def\arraymode{\ifold\relax\else\displaystyle\fi} 
\def\eqnumphantom{\phantom{(\theequation)}}     
\def\@arrayskip{\ifold\baselineskip\z@\lineskip\z@
     \else
     \baselineskip\minarrayskip\lineskip2\minarrayskip\fi}
\def\@arrayclassz{\ifcase \@lastchclass \@acolampacol \or
\@ampacol \or \or \or \@addamp \or
   \@acolampacol \or \@firstampfalse \@acol \fi
\edef\@preamble{\@preamble
  \ifcase \@chnum
     \hfil$\relax\arraymode\@sharp$\hfil
     \or $\relax\arraymode\@sharp$\hfil
     \or \hfil$\relax\arraymode\@sharp$\fi}}
\def\@array[#1]#2{\setbox\@arstrutbox=\hbox{\vrule
     height\arraystretch \ht\strutbox
     depth\arraystretch \dp\strutbox
     width\z@}\@mkpream{#2}\edef\@preamble{\halign
\noexpand\@halignto
\bgroup \tabskip\z@ \@arstrut \@preamble \tabskip\z@ \cr}%
\let\@startpbox\@@startpbox \let\@endpbox\@@endpbox
  \if #1t\vtop \else \if#1b\vbox \else \vcenter \fi\fi
  \bgroup \let\par\relax
  \let\@sharp##\let\protect\relax
  \@arrayskip\@preamble}
\def\eqnarray{\stepcounter{equation}%
              \let\@currentlabel=\theequation
              \global\@eqnswtrue
              \global\@eqcnt\z@
              \tabskip\@centering
              \let\\=\@eqncr
              $$%
 \halign to \displaywidth\bgroup
    \eqnumphantom\@eqnsel\hskip\@centering
    $\displaystyle \tabskip\z@ {##}$%
    \global\@eqcnt\@ne \hskip 2\arraycolsep
         $\displaystyle\arraymode{##}$\hfil
    \global\@eqcnt\tw@ \hskip 2\arraycolsep
         $\displaystyle\tabskip\z@{##}$\hfil
         \tabskip\@centering
    &{##}\tabskip\z@\cr}
\def\theequation{\thesection.\arabic{equation}}
\newfont{\hr}{msbm10}
\newfont{\ams}{msam10}
\def\b{$\tau$}
\font\teneufm=cmmib10
\font\seveneufm=cmmib7
\font\fiveeufm=cmmib5
\def\bfit#1{{\textfont1=\teneufm\scriptfont1=\seveneufm
\scriptscriptfont1=\fiveeufm
\mathchoice{\hbox{$\displaystyle#1$}}{\hbox{$\textstyle#1$}}
{\hbox{$\scriptstyle#1$}}{\hbox{$\scriptscriptstyle#1$}}}}

\def\balpha{{\bfit\alpha}}
\def\bbeta{{\bfit\beta}}

\def\bmu{{\bfit\mu}}

\def\bphi{{\bfit\phi}}
\def\blambda{{\bfit\lambda}}
\def\brho{{\bfit\rho}}
\def\bj{{\bfit j}}
\def\bfe{{\bfit e}}

\def\tr#1{{\rm tr}\kern-3pt\left[#1\right]}

\def\bea{\begin{eqnarray}}
\def\eea{\end{eqnarray}}
\def\bqa{\begin{eqnarray}}
\def\eqa{\end{eqnarray}}

\def\nn{\nonumber}
\def\beq{\begin{equation}}
\def\eeq{\end{equation}}
\def\ba{\beq\new\begin{array}{c}}
\def\ea{\end{array}\eeq}
\def\be{\ba}
\def\ee{\ea}
\def\2{{1\over 2}}
\def\stackreb#1#2{\mathrel{\mathop{#2}\limits_{#1}}}

\def\f{1\over}
\def\Tr{{\tr}}

\def\tr{{\hbox{tr}}}

\def\d{\partial}

\def\b{$\tau$}

\def\R{\bfrm R}
\def\Z{\bfrm Z}
\def\C{\bfrm C}
\def\G{\bfrm G}

\begingroup\ifx\undefined\newsymbol \else\def\input#1 {\endgroup}\fi
\input amssym.def \relax
\input amssym
\font\teneufm=cmmib10
\font\seveneufm=cmmib7
\font\fiveeufm=cmmib5
\def\bfit#1{{\textfont1=\teneufm\scriptfont1=\seveneufm
\scriptscriptfont1=\fiveeufm
\mathchoice{\hbox{$\displaystyle#1$}}{\hbox{$\textstyle#1$}}
{\hbox{$\scriptstyle#1$}}{\hbox{$\scriptscriptstyle#1$}}}}
\font\tenrmbf=cmbx12
\font\sevenrmbf=cmbx7
\font\fivermbf=cmbx5
\def\bfrm#1{{\textfont1=\tenrmbf\scriptfont1=\sevenrmbf
\scriptscriptfont1=\fivermbf
\mathchoice{\hbox{$\displaystyle#1$}}{\hbox{$\textstyle#1$}}
{\hbox{$\scriptstyle#1$}}{\hbox{$\scriptscriptstyle#1$}}}}

\newcommand{\ti}[1]{\tilde{#1}}

\newcommand{\ms}{\mapsto}

\newcommand{\beql}[1]{\be\label{#1}}
\newcommand{\eq}{\ee}

\begin{document}
\setcounter{footnote}0
\begin{center}
\hfill ITEP/TH-58/97\\
\hfill FIAN/TD-13/97\\
\hfill hep-th/9711006\\
\vspace{0.3in}
{\LARGE\bf $\tau$-function within group theory approach}\\
\vspace{0.1in}
{\LARGE\bf and its quantization}
\\
\bigskip
\bigskip
{\Large
A.Mironov
\footnote{E-mail address:
mironov@lpi.ac.ru, mironov@itep.ru}}
\\
\bigskip
{\it Theory Department, P.N.Lebedev Physics
Institute, Moscow, Russia}\hfill
\centerline{and}
{\it ITEP, Moscow, Russia}\hfill
\end{center}
\bigskip 

\begin{abstract}
This is a review of the results related to generalizations of the notion of
$\tau$-function and integrable hierarchies and to their interpretation within
the group theory framework that admits an immediate quantization procedure.
Different group theory structures underlying integrable system are discussed
in detail as well as their quantum deformations.
\end{abstract}

\tableofcontents

\setcounter{footnote}0
\section{Introduction}
\setcounter{equation}{0}
During last years, integrability gets more and more important in particle
physics and string theory \cite{Mardis,Mirdis}. One of the main lessons that
can be learned from the numerous concrete examples is that the {\it
classical} integrability in {\it quantum} field theories turns out to be
surprisingly typical.

One of the first examples of using essentially integrable systems in this
field is attempts to perturbatively deform conformal field theories
along the "integrable" directions \cite{Zam1}. This is ultimately
necessary for constructing the non-perturbative solutions in string
theory. The perturbative deformations certainly describe only the vicinity of
the critical points, while not exposing the global structure of the
configuration space of string theory. Still, it is the papers
\cite{Zam1} (see also \cite{Sasaka}) that have
already pointed out the {\it integrable} deformations of conformal field
theories. Moreover, within the string theory approach it has been also
expected that the sum of the whole perturbative series for the amplitudes
(and the partition function) can be described by a quantum integrable system
\cite{Mor2}. At this point, it is of great importance that the effective
action for the {\it quantum} integrable system is a $\tau$-fucntion of some
{\it classical} integrable system\footnote{Surprisingly enough, integrating
some objects like $\tau$-function \cite{Saito} over the
universal moduli space in string theory \cite{Mor4} leads to a
$\tau$-function given on the same universal moduli space.}. This amazing
phenomenon has been already realized in concrete examples
\cite{AR}-\cite{Kor}, while never been explained.

It is interesting that all the examples can be parted into two
large classes of the systems naturally depending on times (coupling
constants) and on the so called Miwa variables, these latter ones being
just eigenvalues of an (infinite) external matrix. An instance of the first
class system is a functional of the correlators in the quantum non-linear
Schr\"odinger model, which is a $\tau$-function of some classical integrable
non-linear equation of the same kind \cite{BIKS}. A typical second class
example  is given by the partition function in the six-vertex model with
non-trivial boundary conditions. This partition function turns out to be the
$\tau$-function of the $2d$ Toda system expressed in Miwa variables (that
are, in essence, the boundary conditions) \cite{Kor}.

However, the most important case of the integrability phenomenon under
consideration is given by the matrix models \cite{Mor3,Marrev,Mirdis}.
Among other things, in the matrix models the both types of time variables
are realized. Indeed, the partition function of the discrete matrix
models is a $\tau$-function, the times being coupling constants, while that
of the continuous matrix models is a $\tau$-function depending on an external
matrix that gives rise to Miwa variables.

In fact, the classical integrability in effective field theory is quite
typical and non-specific, say, for two dimensions. A good example is given by
the solution of the $4d$ $N=2$ supersymmetric Yang-Mills theory, which has
been recently constructed by N.Seiberg and E.Witten \cite{WS}. Classical
integrability in this case has been revealed in the paper \cite{WSI1}
(for recent progress see \cite{Mardis,Mir97,GGM} and references therein).
Since integrability looks so typical, one could expect for the effective
action in typical field theory to be a $\tau$-function provided the variables
are choosen properly and all possible interactions are taken into account.
In practice, these latter are usually restricted by some additional
requirements (renormalizability etc).

It was already noted above that, although the integrability phenomenon has
been investigated in a series of examples, its origin still remains
mysterious. In the present paper, we deal with issues related to the general
properties of quantum and classical $\tau$-functions. This would help us to
shed some light on the phenomenon. Besides, we discuss here some immediate
applications of the general approach. For example, we calculate the
$S$-matrix of the quantum Liouville theory.

In order to discuss how correlators of quantum integrable system are related
to classical hierarchies, one needs to deal with the quantum system in
terms suitable for the classical hierarchies. The key notion in these
hierarchies is the $\tau$-function, while the quantum integrable systems have
been formulated so far only in $R$-matrix terms \cite{FRT,FTQ}\footnote{See,
however, very interesting papers \cite{qt}.}.
Therefore, the main aim of this paper is to extend the notion of
$\tau$-function to the quantum case. To this end, we are going to start with
a detailed discussion of the usual $\tau$-function and, in particular, to
point out related group theory (symmetry) structures. Following this way, we
are naturally led to the notion of the {\it generalized} $\tau$-function that
can be still algebraically described and satisfies a set of Hirota-like
bilinear identities.

One can hardly over-estimate the necessity of the generalization of the
$\tau$-function notion. To begin with, the matrix models corresponding to the
standard hierarchies are just the simplest models. More involved and
physically important examples include the generalized
Kontsevich-Kazakov-Migdal model \cite{GKM2}, higher dimensional gauge
theories \cite{Wilact}, $2d$ gravity coupled to the conformal matter with
$c>1$ (i.e. with space-time dimension $d>2$ which again corresponds to the
higher-dimensional models) etc. In all these examples, the angular variables
of the matrices can not be factorized out and then integrated out. Thus, they
are described by more complicated "generalized" hierarchies.

Another important class of the generalized integrable systems has to
correspond to the WZWN theories with the level greater that 1. Further moving
into this direction should lead to the 2-loop, 3-loop etc. algebras within the
integrable framework. All the listed examples are associated with the so
called non-Cartanian hierarchies \cite{Mor3}.

Yet another important problem to be solved in the course of quantizing the
generalized $\tau$-functions is the careful regularization of the Liouville
theory \cite{Weight}. This problem can be also differently approached, as we
discuss later.

Let us remark here that we do not discuss why the
$q$-deformation described in the paper has something to do with
quantizing the system, i.e. with building quantum field theory (quantum
mechanics) in place of classical one. Instead, we address to the
papers \cite{FRT,FTQ} that discuss why turning to quantum groups
and quantum $R$-matrices corresponds to quantizing classical systems.

It was mentioned above that the main tool in studies of the generalized
$\tau$-functions is group theory. More concretely, one can use different
symmetry structures of the integrable hierarchies to formulate them in
algebraic terms. It is typical for the theoretical physics that the symmetry
approach allows one to derive the general laws that do not depend on dynamics.
In particular, it is the general symmetry structure that should explain the
integrability of effective actions in quantum field theories. This is why one
needs to reveal different algebraic structures and their role for the
quantization.

As a typical example, in this paper that is a review of the papers
\cite{MMV}-\cite{Mir96} we consider either the Toda integrable system, or
its particular case -- the Liouville system. In sections 3 and 6, it is
demonstrated that there are three different group structures in these
systems. Two of them are related with the group action in the spectral
parameter space: one is the group $GL(\infty)$ whose Grassmannian
parametrizes all solutions of the $2d$ Toda hierarchy of the general form,
the other one is the group $SL(n)$ that describes the reduction of this whole
large space of solutions to the solutions of the $SL(n)$ Toda. Thus, both
these groups act on solutions of the integrable hierarchy. It is these groups
that are substituted by their quantum counterparts in the course of
quantization.

Let us note, however, that the quantization procedure is not always a drill,
that is to say, not any group structure would be automatically substituted by
its quantum counterpart. Indeed, in the Liouville type system there is yet
another group acting in ``the space-time", i.e. on times of integrable
system. In the course of quantization, this group still remains classical
although is differently treated -- the orbit interpretation is involved.
In fact, if the classical system is obtained by the Hamiltonian reduction of
a free system given on the cotangent bundle to a simple real Lie group
\cite{OP2}, the quantum model is related to an irreducible
unitary representation of {\it the same} group. Therefore, the quantum system
is rather treated within the geometrical quantization framework \cite{Hart}.
This group structure is discussed in section 6. In the section, it is also
demonstrated that the group theory methods allows one not only to construct
the wave function of the quantum Toda system (that is called Whittaker
function) \cite{Jac}, but also to calculate its asymptotics
\cite{GinKor} and, therefore, $S$-matrix (i.e. two-point
correlation function).

Such calculations are known for a while \cite{STSh} but they have been
performed using the Iwasawa decomposition that is inconvenient for the
extension to the affine case. In section 6, we present alternative
calculation that is based on the Gauss decomposition, which is naturally
generalizable to the Liouville field theory described by the affine algebra
(in the affine case, one needs to use the ``point-wise" Gauss decomposition
analogous to that used when bosonizing the WZWN theories \cite{GMMOS}). All
the calculations of the wave function and its asymptotics are immediately
continued to this (affine) case giving the $S$-matrix of the $2d$ Liouville
theory. This $S$-matrix coincides with that obtained recently by
A. and Al.Zamolodchikovs \cite{ZZ} and H.Dorn and H.-J.Otto \cite{DO}.

The other group structures -- those $q$-deformed in the course of
quantization -- are discussed in the other sections. In particular, in
section 4, the generalized $\tau$-functions admitting immediate quantization
are discussed. It is natural to require for such $\tau$-functions
still to remain $c$-numbers. In s.4.5 such a deformation of the classical
differential hierarchy is, indeed, constructed so that the resulting
hierarchy becomes a difference one. Nevertheless, this hierarchy turns out to
be too trivial, since it can be obtained by a (quite complicated)
redefinition of times of the original hierarchy. Therefore, the problem of
constructing a non-trivial deformation of the $\tau$-function (that leads to
a non-commutative object) still persists.

This problem is solved, however, in the pure algebraic terms -- in the same
section 4, given a representation of any algebra, the $\tau$-function is
defined as the generating function of its matrix elements in some
(arbitrary) representation. Thus defined $\tau$-function satisfies a system
of bilinear identities (BI).  In the classical case, these BI are
differential ones. To quantize the system, it suffices to substitute the
classical group by its quantum counterpart.  This makes the differential BI
difference and the $\tau$-function itself -- operator.

Note that the $\tau$-function defined by an arbitrary representation of
algebra leads to BI given by {\it non-commutative} flows, that is to say,
the theory is no longer integrable in the usual (Liouville) sense. However,
it still preserves many important features of the integrable system and
corresponds to the non-Cartanian hierarchies \cite{Mor3} that we mentioned
above.

In order to obtain the $\tau$-function corresponding to the standard KP/Toda
hierarchy, one needs to consider the fundamental representations of the group
$SL(\infty)$. For these representations, there exists a special reduction of
the general $\tau$-function, which leads to the standard integrable
hierarchies.

The next problem is to construct some quantum counterpart of the KP hierarchy.
This problem is highly non-trivial, since the mentioned special reduction of
the general $\tau$-function admissible in the fundamental representations
does not endure the quantization. Nevertheless, one can construct the
$q$-deformation of some new classical hierarchy that is, by essence, the KP
hierarchy but with unusual evolution. This is discussed in section 5.

For the sake of continuity of the paper, we place more technical comments in
Appendices. Some additional details that are not included into this review
can be found in the original papers \cite{MMV}-\cite{Mir96} and, especially,
\cite{GKMMMO}, since the results of this latter are considered here very
briefly.

The author is grateful to A.Gerasimov, S.Kharchev, S.Khoroshkin, D.Lebedev,
A.Marshakov, A.Morozov, M.Olshanetsky and L.Vinet for the collaboration and
useful discussions and to G.Weight, A.Gorsky, V.Dobrev, A.Zabrodin,
Al.Zamolodchikov, V.Korepin, N.Nekrasov, N.Slavnov, I.Tyutin and J.Schnittger
for the useful discussions. The work is partially supported by grants
RFBR-96-01-01106 and INTAS-97-1038.

\section{Fermionic representation for integrable hierarchies}
\setcounter{equation}{0}
We start with a brief description of the fermionic approach to the standard
integrable hierarchies that is mostly due to the Kyoto school
\cite{J1}-\cite{J8} (see also \cite{UT,KMMM}) and describe in these terms the
important particular case of the semi-infinite (forced)
hierarchies\footnote{Right these hierarchies describe, among other, the
matrix models.}. Note that the fermionic approach is the most universal
language of those applied so far to describe the integrable systems, since it
corresponds to the formal (Sato \cite{Sato}) Grassmannian given in terms of
the formal series. Although alternative approaches (see, e.g.,
\cite{Kri3}) suit better, say, the finite-gap solutions
\cite{ZMNP}, it is the Sato Grassmannian that is necessary for dealing with
the singular points (describing, for instance, the matrix models, see
\cite{KSh}).

\subsection{General fermionic description of the standard integrable
hierarchies}
We begin with the description of the two-dimensional lattice Toda
(2TDL) hierarchy. Its $\tau$-function is defined through correlators in
the theory of free fermions $\psi, \tilde\psi$ ($b,c$-system of spin $1/2$).
Namely, the $\tau$-function is the ratio of the
correlators\footnote{Later on, we omit from this definition the
denominator correlator (which is nothing but a normalization factor) wherever
it is non-singular and, therefore, inessential.}
\beq\label{tau}
\tau _n(t,\bar t \mid G) \equiv \frac
{\langle n \mid e^H G e^{\bar H} \mid n \rangle}
{\langle n \mid  G  \mid n \rangle}
\eeq
in the theory of the free fermionic fields $\psi (z)$,  $\psi^{\ast} (z)$:
\ba\label{psi}
\psi (z) =
\sum _{\Z}
\psi _nz^n \\
\psi^{\ast} (z) =
\sum _{\Z}
\psi^{\ast} _nz^{-n-1}
\ea
\be\label{psicommrel}
\{\psi_k,\psi^{\ast}_m\}=\delta_{km},\ \ \ \{\psi_k,\psi_m\}=
\{\psi_k^{\ast},\psi^{\ast}_m\}=0
\ee
In the definition (\ref{tau}) we use the following notations
\beq\label{Ham}
H = \sum_{k > 0} t_kJ_k; \ \ \ \bar H = \sum_{k>0} \bar t_k J_{-k}
\eeq
where the currents are defined as
\beq\label{cur}
J(z) = \psi^{\ast} (z)\psi (z),\ \ \hbox{i.e.}\ \ J_k=\sum_{{\Z}} \psi_i
\psi_{i+k}^{\ast}
\eeq
The quantity $G$ is given by the formula
\beq\label{elGr}
G =\
:\exp \left\{ \sum_{m,n} {\cal G}_{mn}\psi^{\ast} _m\psi _n\right\}:
\eeq
and is an element of the group $GL(\infty)$ (or the Clifford group)
realized in the infinite dimensional Grassmannian (see also
equivalent descriptions in
\cite{Kri3}). The normal ordering is to be understood here as
ordering w.r.t. the vacuum $|0\rangle$ defined by the conditions
\beq\label{vac0}
\psi _m|0\rangle  = 0\ \ m < 0\hbox{ , }
\psi^{\ast} _m|0\rangle  = 0\ \ m \geq  0
\eeq
One can also consider a more general normal ordering that is defined w.r.t.
the vacuum $|k\rangle$ given by the conditions
\beq\label{vac}
\psi _m|k\rangle  = 0\ \ m < k\hbox{ , }  \psi^{\ast} _m|k\rangle  = 0\ \
m \geq  k
\eeq
Any concrete solution to the hierarchy ($\tau$-function (\ref{tau}))
depends only on the choice of element $G$
(or, equivalently, can be uniquely defined by the matrix
${\cal G}_{km}$). The remarkable result by the Kyoto school is the statement
that any solution to the 2TDL hierarchy is given by an element
$G$ of the form (\ref{elGr}), and, inversely, any such element
$G$ defines some solution to the 2TDL hierarchy.

Note that, within the fermionic approach \cite{J1}-\cite{J8}, 2TDL is the
most general one-component hierarchy.  Say, the KP hierarchy can be obtained
from the 2TDL hierarchy by canceling all the negative times and the zero
time. One can also arbitrarily fix these times and do not include the
corresponding evolution equations. In this respect, the KP hierarchy is a
sort of sub-hierarchy of 2TDL, not a reduction, and is described by the same
set of elements $G$ with less number of the evolution flows.

Another example is the Toda chain which is already a reduction that can be
obtained from 2TDL by special restricting element of the Grassmannian.
As a result, the $\tau$-function of the Toda chain depends not on the
positive and negative times separately but only on their sums
\cite{UT}\footnote{In most papers, the Toda chain depends only on {\it
difference} of times. It is due to the opposite sign of $\bar H$ used in
(\ref{tau}).} (therefore, one can consider only one set of the evolution
flows). This property of the Toda chain can be taken as its definition.

Now we derive some useful properties of the described system. One can easily
get from the commutation relations (\ref{psicommrel}) the transformation
law of the fermionic modes w.r.t. the action of the element of the
Grassmannian (\ref{elGr}):
\beq\label{A16}
G\psi _kG^{-1} = \psi _jR_{jk}\hbox{ , }  G\psi ^\ast _kG^{-1} =
\psi ^\ast _jR^{-1}_{kj}
\eeq
where the matrix $R_{jk} $ can be expressed through ${\cal G}_{jk}$ (see
\cite{HQM}). It will be clear later that $R_{jk}$ is an important building
block for determinant representations of the $\tau$-functions.

Let us introduce more notations. Using (\ref{psicommrel}), one can define
evolution of $\psi (z)$ and $\psi ^\ast (z)$ w.r.t. the times flows
$\{t_k\}$, $\{\bar t_k\}$:
\beq\label{A17}
\psi (z,t) \equiv  e^{H(t)}\psi (z)e^{-H(t)} = e^{\xi (t,z)}\psi (z)
\eeq
\beq
\label{A18}
\psi ^\ast (z,t) \equiv  e^{H(t)}\psi ^\ast (z)e^{-H(t)} =
e^{-\xi (t,z)}\psi ^\ast (z)
\eeq
\beq
\label{A19}
\psi (z,\bar t) \equiv  e^{\bar H(\bar t)}\psi (z)e^{-\bar H(\bar t)} =
e^{\xi (\bar t,z^{-1})}\psi (z)
\eeq
\beq
\label{A20}
\psi ^\ast (z,\bar t) \equiv  e^{\bar H(\bar t)}\psi ^\ast (z)e^{-\bar H(\bar
t)} = e^{-\xi (\bar t,z^{-1})}\psi (z)
\eeq
where
\beq\label{A21}
\xi (t,z) = \sum ^\infty _{k=1}t_kz^k
\eeq
Let us also define the Schur polynomials
\beq\label{Shur}
\exp\left\{\sum_{k>0}t_kx^k\right\} \equiv \sum _{k>0} P_k(t_k)x^k,\
\ \ \ P_k(t_k)=0,\ \ k<0
\eeq
Using this definition and (\ref{A17})-(\ref{A20}), one immediately gets
the evolution of the fermionic modes:
\beq\label{A22}
\psi _k(t) \equiv e^{H(t)}\psi _ke^{-H(t)}
=\sum ^{\infty}_{m=0} \psi _{k-m}P_m(t)
\eeq
\beq
\label{A23}
\psi ^{\ast }_k(t) \equiv e^{H(t)}\psi _k^{\ast }e^{-H(t)} =
\sum ^{\infty}_{m=0} \psi ^{\ast }_{k+m}P_m(-t)
\eeq
\beq\label{A24}
\psi _k(\bar t) \equiv e^{\bar H(\bar t)}\psi
_ke^{-\bar H(\bar t)} =\sum ^{\infty} _{m=0} \psi
_{k+m}P_m(\bar t)
\eeq
\beq
\label{A25}
\psi _k^{\ast }(\bar t) \equiv e^{\bar H(\bar t)}\psi
_k^{\ast }e^{-\bar H(\bar t)}=\sum ^{\infty} _{m=0} \psi
_{k-m}^{\ast } P_m(-\bar t)
\eeq
Now we can derive a determinant representation for the
$\tau$-function\footnote{The discussion below requires, certainly,
in real calculations some convergency properties, in particular,
of the determinants of infinite-dimensional matrices. It suffices to
consider the matrices with finite number of non-zero diagonals
\cite{DJKM}. However, sometimes this class of matrices turns out to be
too restrictive like being the case in matrix models \cite{GKM,KSh}.
On the other hand, the formulas below are direct implications of some
relations that can be checked out for formal series and, therefore, are
always correct once the convergency conditions are fulfilled. By the same
reason, we ignore the issue of existence of the completely filled vacuum
state below.}.  To this end, define the completely filled vacuum state
$|-\infty \rangle$ that satisfies the requirement
\beq\label{A26}
\psi ^\ast
_i|-\infty \rangle  = 0\hbox{  , }  i \in  {\Z}
\eeq
Then, any "shifted"
vacuum is obtained from the completely filled one as follows:
\beq\label{A27}
|n\rangle  = \psi _{n-1}\psi _{n-2} ...|-\infty \rangle
\eeq
Note that any element $G$ of the Clifford group (and, therefore,
$e^{\bar H(\bar t)}$) acts on $|-\infty \rangle $ trivially:
$G|-\infty \rangle  \sim  |-\infty \rangle $. Therefore, using
(\ref{A23}) and (\ref{A24}), one easily gets from (\ref{tau}):
\beq
\new
\begin{array}{c}
\tau _n(t,\bar t) =
\langle -\infty |...\psi ^\ast _{n-2}(-t)\psi ^\ast _{n-1}(-t)G\psi _{n-1}(
\bar t)\psi _{n-2}(\bar t) ...|-\infty \rangle  \sim \\
\sim \det [\langle -\infty |\psi ^\ast _i(-t)G\psi _j(\bar t)G^{-1}|-\infty
\rangle ]\left| _{i,j \leq  n-1}\right.
\end{array}\label{A28}
\eeq
Using (\ref{A16}), we now see that
\beq\label{A29}
G\psi _j(\bar t)G^{-1} =
\sum _{m,k} P_m(\bar t)\psi _kR_{k,j+m}
\eeq
and, therefore, the "manifest" solution to the 2TDL hierarchy reads
in the determinant form as
\beq\label{A30}
\tau _n(t,\bar t) \sim \det_{i,j<0} \ {C}_{i+n,j+n} (t,\bar t)
\eeq
where\footnote{Note that throughout we ignore the summation limits, since
they are automatically reproduced by the property (definition) of the
Schur polynomials $P_k=0$ at $k<0$ in (\ref{Shur}).}
\beq\label{A31}
{C}_{ij}(t,\bar t) = \sum _{k,m} R_{km}P_{k-i}(t)P_{m-j}(\bar t)
\eeq

Let us note that the time dependence of the matrix elements (\ref{A31}) in
the determinant (\ref{A30}) is given by the following
equations:
\be\label{h1}
\partial C_{ij}/\partial t_p = C_{i,j-p}\hbox{, }   j>p>0
\ee
\be
\label{h2}
\partial C_{ij}/\partial \bar t_{p} =  C_{i-p,j}\hbox{, }   i>p>0
\ee
which immediately follows from the corresponding property of the
Schur polynomials:
\beq\label{h3}
\partial P_k/\partial t_p = P_{k-p}
\eeq
that is a direct consequence of their definition.

It has been already remarked that the KP hierarchy is given by the
evolution only w.r.t. the positive times $\{t_k\}$, while the negative
times $\{\bar t_k\}$ are just parameters which give a set of points
in the Grassmannian and can be removed by a redefinition of the matrix
$R_{km}$. Then, the $\tau$-function of the (modified) KP hierarchy
reads as
\beq\label{A32}
\tau _n(t) = \langle n|e^{H(t)}G(\bar t)|n\rangle \sim \det_{i,j<0} \left[\sum
_k
R_{k,j+n}(\bar t)P_{k-i-n}(t)\right]
\eeq
where  $G(\bar t) \equiv  Ge^{\bar H(\bar t)}$ and
\beq\label{A33}
R_{kj}(\bar t) \equiv  \sum  _m R_{km}P_{m-j}(\bar t)
\eeq

\subsection{Bilinear identities for the classical 2TDL hierarchy}
For future use, we need two more {\it crucial} formulas in the theory of
integrable systems \cite{J3}:
\be
\langle n|\psi(z)\exp [H(t)]
=z^{n-1}\langle n-1|\exp [H(t- \epsilon (z^{-1}))]
\equiv z^{n-1}\hat X(z,t)\ \langle n-1|e^{{\displaystyle H(t)}}
\label{main1}
\ee
\be
\langle n|\psi^{\ast}(z)\exp [H(t)]
=z^{-n}\langle n+1|\exp [H(t+ \epsilon (z^{-1}))]
\equiv z^{-n}\hat X^{\ast}(z,t)\langle n+1|e^{{\displaystyle H(t)}}
\label{main2}
\ee
where
\be\label{A8}
\hat X(z,t) \equiv e^{\xi(z,t)} \; e^{-\xi(z,\tilde \partial_{t})}
\ee
\be\label{A9}
\hat X^{\ast}(z,t) \equiv e^{-\xi(z,t)} \; e^{\xi(z,\tilde \partial_{t})}
\ee
where as usual $\ti\d_t\equiv (\d_{t_1},\2\d_{t_2},\ldots)$.

Using these formulas, one can express the Baker-Akhiezer (BA)
function of the KP
hierarchy, defined via $\tau$-function by the formula \cite{J3}
\beq
\Psi (\mu |t_n) = e^{\displaystyle{{\Sigma t_n\mu^n}}} {\tau_0(t_n -
\left.{\mu ^{-n}\over n}\right|\bar t_n)\over \tau_0(t_n|\bar t_n)}
\eeq
through the fermionic correlator:
\be\label{BA1}
\Psi (\mu |t_n) ={\left<1\left|e^{H(t)}\psi(\mu)G\right|0\right>\over
\left<0\left|e^{H(t)}G\right|0\right>}
\ee
Analogously, one can consider the conjugated BA function given by
\be\label{BAC1}
\Psi^{\ast} (\mu |t_n)
={\left<1\left|e^{H(t)}\psi^{\ast}(\mu)G\right|0\right>\over
\left<0\left|e^{H(t)}G\right|0\right>}
\ee
Involving "non-zeroth" vacuums, one can also introduce more general pair of
the BA functions depending on the zero time $n$.

Considering the "whole" 2TDL hierarchy, one can reasonably introduce {\it
four} different BA functions, since, in this case, it makes sense to put
fermions both to the left and to the right of the Grassmannian element
$G$.

Now let us derive the bilinear identities (BI) satisfied by the
$\tau$-function of the 2TDL hierarchy (\ref{tau}) \cite{DJKM,UT}. To this
end, note that due to (\ref{A16}) the tensor product of the Grassmannian
elements $G\otimes G$ (\ref{elGr}) commutes with the tensor product
$\Gamma\equiv\sum_i\psi_i\otimes\psi_i^{\ast}$. To make our notations
consistent with the forthcoming sections, we denote hereafter
the fermionic modes $\psi_k$ through $\psi^+_k$ and $\psi^{\ast}_k$ through
$\psi^-_k$.

Look now at the matrix elements of the operator identity
\be\label{5s1}
\Gamma(G\otimes G)=(G\otimes G)\Gamma
\ee
taken in between the states
$\langle n+1|U(t)\otimes \langle m-1|U(t')$ and
$\bar U(\bar t)|n\rangle\otimes \bar U(\bar t')|m\rangle$, where we denote
$U(t)\equiv e^{H(t)}$ and $\bar U(\bar t)\equiv e^{\bar H(\bar t)}$:
\be\label{BIglinfty}
\sum_i \langle n+1| U(t) \psi_i^+ G \bar U(\bar t)
|n\rangle \cdot \langle m-1 | U(t')
\psi_i^{-} G \bar U(\bar t')
|m\rangle =\\
=\sum_i \langle n+1 | U(t) G\psi_i^+  \bar U(\bar t)
|n\rangle \cdot \langle m-1 | U(t')
G\psi_i^{-}  \bar U(\bar t')
|m\rangle
\ee
Then, one can rewrite (\ref{BIglinfty}) through the free fermion fields
(\ref{psi}):\footnote{The integration contour in this formula depends on
whether the integrand is expanded into positive or negative powers of
$z$.}
\be
\oint_{\infty} {dz}
\langle n+1 | U(t) \psi^+(z) G \bar U(\bar t)
|n\rangle \cdot \langle m-1 | U(t')
\psi^{-}(z) G \bar U(\bar t')
|m\rangle =\\=
\oint_0{dz}
\langle n+1 | U(t) G\psi^+(z)  \bar U(\bar t)
|n\rangle \cdot \langle m-1 | U(t')
G\psi^{-}(z)  \bar U(\bar t')
|m\rangle
\ee
Generalizing definitions of the BA functions (\ref{BA1})
and (\ref{BAC1}) to the non-zero vacuum and involving negative times
\be\label{BA}
\Psi^{\pm,i}_n\equiv \langle n\pm 1 | \hat U(t)
\psi_i^{\pm} G \hat{\bar U}(\bar t) |n\rangle
\ee
one can rewrite these formulas in a more compact
form (hereafter we use the notation $t$ both for positive and for
negative times all together when it can not mislead):
\be\label{BIA}
\sum_i\Psi^{+,i}_k(t)\Psi^{-,i}_l(t')=
\sum_j\bar\Psi^{+,j}_{k+1}(\bar t)\bar\Psi^{+,j}_{l-1}(\bar t')
\ee
Here $\bar\Psi$ gives the second pair of the BA functions defined by the
correlator with the fermion to the right of the element $G$.

These BA functions are generated by the action of the vertex operators
(\ref{A8})-(\ref{A9}) (see (\ref{main1})-(\ref{main2}):
\be
\sum_i\Psi^{+,i}_k(t)z^{i}\equiv \hat X^{+} (z,t)\tau_n(t),\ \ \
\sum_i\Psi^{-,i}_k(t)z^{-i-1}\equiv \hat X^{-} (z,t)\tau_n(t)
\ee
and analogously for ${\hat {\bar X}}^{\pm}(z,t)$.
Then, (\ref{BIglinfty}) can be rewritten as
\be
\oint_{\infty}{dz} \hat X^-(z,t)\tau_n(t,\bar t) \hat
X^{+}(z,t')\tau_m(t',\bar t')= \oint_0{dz}\hat
{\bar X}^-(z,\bar t)\tau_{n+1}(t,\bar t) \hat {\bar X}^{+}(z,\bar
t')\tau_{m-1}(t',\bar t')
\ee
This integral BI can be reduced to the infinite set of equations
(see also (\ref{Hirota2})) by expanding into
degrees of $t-t'$ and $\bar t-\bar t'$:
\be\label{KPHirota}
\sum ^\infty _{i=0} {P}_i(-2y){P}_{i+1}(\tilde D_t)
e^{[\sum _iy_iD_{t_i}]} \tau \cdot \tau  = 0\\
\tau_n\d_{ t_1}\d_{\bar t_{}}\tau_n- \d_{ t_1}\tau_n \d_{\bar t_{1}}\tau_n=
\tau_{n+1}\tau_{n-1}
\ee
where $D_t$ is the Hirota symbol given by its action onto the product of two
functions:\\
$\left.D^k_xf\cdot g\equiv \d_y[f(x+y)g(x-y)] \right|_{y=0}$, $\tilde D
\equiv (D_{t_{^1}}$, ${1\over 2} D_{t_{^2}},\ldots$), ${P}_i$'s are Schur
polynomials. Expanding these equations (\ref{KPHirota}) w.r.t. the arbitrary
set of parameters $\{y_k\}$, one finally arrives at the KP/Toda hierarchy
equations.

\subsection{Fermionic representation for the forced hierarchies}
Now we apply the described formalism to the case of the so-called forced
hierarchies and obtain, in particular, determinant formulas of
\cite{KMMOZ,KMMM}. The main difference of the formulas like
(\ref{A30}) with the determinant representations to be obtained
in this subsection is that the forced hierarchies corresponds to the
determinants of the {\it finite} matrices. This important difference is
due to the very specific boundary condition (from the point of view
of the particle chain)
\beq\label{forced}
\tau_0 = 1
\eeq
This condition may serve as a definition of the forced hierarchies
\cite{Kaup} investigated in \cite{KMMOZ,KMMM} in detail. It turns
out that the forced hierarchies are described by singular
elements of the Grassmannian, which depend on only positive (or
only negative) fermionic modes. It implies "quarter-infinite"
matrices ${\cal G}_{mn}$ in (\ref{elGr}) and, indeed, leads to the
determinants of finite matrices. We describe the forced hierarchies in this
subsection very briefly leaving some technical details till Appendix 1. More
details can be found in \cite{KMMOZ,KMMM}.

Thus, we are looking for the point of the Grassmannian in the form
\beq\label{elGrforced}
G = G_0P_+
\eeq
where $P_+$ is the projector onto the positive states:
\beq
P_+|n\rangle  = \theta (n)|n\rangle
\eeq
This projector admits a natural fermionic realization
\beq
P_+ = :\exp \left[\sum _{i<0}\psi _i\psi ^\ast _i\right]:
\eeq
and possesses a set of properties
\beq
P_+\psi ^\ast _{-k} = \psi _{-k}P_+ = 0\hbox{  , }  k > 0\
\eeq
\beq
\left[ P_+,\psi _k\right] = [P_+,\psi ^{\ast} _k] = 0\hbox{  , } k \geq  0\
\eeq
\beq
P^2_+ = P_+
\eeq
If the correlator in (\ref{tau}) contains this projector, we naturally
require for $G_0$ to depend only on $\psi _k$ and $\psi
^\ast _k$ with $k \geq 0$. We fix its form to be
\beq
G_0 = :\exp \left\{
\left(\int _\gamma  A(z,w)\psi _+(z) \psi ^\ast _+(w^{-1})dzdw \right) - \sum
_{i\geq 0} \psi _i\psi ^\ast _i \right\}:
\eeq
where $\gamma $ is an integration domain. The key formula (we leave its
proof till Appendix 1) gives the $\tau$-function corresponding to the
element of the Grassmannian choosen in the form (\ref{elGrforced}):
\be
\tau _n(t,\bar t) = \langle n|e^{H(t)}G_0P_+e^{\bar
H(\bar t)}|n\rangle  = \\ = {1\over n!} \int _\gamma \Delta (w)\Delta
(z)\prod ^n_{i=1}A(z_i,w_i)e^{\xi (t,z_i)+\xi (\bar t,w_i)}dz_idw_i
\ee
where $\Delta(z)\equiv \det_{i,j} z^{i-1}_j$ is the Van-der-Monde determinant.

Repeating the above calculation (see also Appendix 1), one can get
the determinant representation for the forced hierarchy:
\beq\label{detrep}
\tau _n(t,\bar t) = \left.\det \left[\partial ^i_{t_1}\partial
_{{\bar t}_1}^j\int _\gamma A(z,w)e^{\xi (t,z)+\xi (\bar t,w)}dzdw\right]
\right| _{i,j=0,...,n-1}
\eeq
This is, in fact, the determinant of the finite matrix, which provides, in
particular cases, the determinant representations for the matrix models
\cite{KMMM}.

Note that the matrix $A_{ij}$ defined by the expansion of
$A(z,w)$ into degrees of $z$ and $w$ coincides with the matrix $R_{ij}$ from
(\ref{A16}) (see Appendix 1). This simple relation of the fermion rotation
matrix and the matrix giving the element of the Grassmannian is
specific for the forced hierarchies.

The element of the Grassmannian determined in (\ref{elGrforced})
is in no way unique, since there are many possibilities to continue
the relation (\ref{forced}) to the negative values of $n$. The
choice we used corresponds to the condition
\be\label{forcedhaha}
\tau_n=0\ \ \ \hbox{at}\ \ \ n<0
\ee
This way of defining the $\tau$-function is advantageous when reducing
to the Toda chain, since, in this case, the $\tau$-function depends only on
the sum of times $t_k+\bar t_k$ that is the defining property of the
corresponding infinite hierarchy. Other important choices for the element
of the Grassmannian $G$ have been discussed in \cite{KMMOZ}. If considering
only dependence on times $t_k$, one of interesting choices corresponds to
the condition
\be
\tau_n(t_k)=\tau_{-n}\left((-)^kt_k\right)
\ee
and describes the (modified) CKP hierarchy \cite{J4}, i.e.
the corresponding element of the Grassmannian belongs to
$Sp(\infty)$. Another choice of this element is just
$G=G_0$ that corresponds to the condition
\be
\tau_n=1\ \ \ \hbox{at}\ \ \ n<0
\ee
etc.

\section{Group theory structures of the Toda molecule}
\setcounter{equation}{0}
Before generalizing the construction of the previous section, we discuss
group theory structures that are presented in integrable systems and take as
our main example the Liouville (Toda) theory\footnote{In fact, it has been
realized already quite a while ago \cite{OP,OP1} (see also, by the
appearance, the first paper on this stuff \cite{Bo}) that the group theory
approach is one of the most effective and elegant frameworks for classical
integrable systems.}. Note that it often happens for {\it different} group
structures to be presented in the {\it same} integrable system. In
particular, some groups can act in the solution space of integrable
hierarchy, while other ones -- just on the variables of integrable equations
(in "the space-time") \cite{Mir96}. We start with considering the structures
of the first type, leaving those of the second type till section 4.

As it was discussed in the previous section, all solutions of the
integrable 2DTL hierarchy (\ref{tau}) are described by the group elements of
$GL(\infty)$ (\ref{elGr}). This group acts in the space of solutions and
simultaneously in the space of the spectral parameter. This is the
group structure we are mainly interested in, since just this group
structure is in charge of the majority of defining properties of integrable
system and it is deformed in the course of quantizing. However, in the
example we are going to look at now, that is the Liouville system, there
exists yet another group.

We begin with a natural restriction of the forced hierarchies given by the
condition (\ref{forcedhaha}) onto finite systems. To this end, we impose
one additional constraint
\be\label{molecule}
\tau_n=0,\ \ \ n>p
\ee
for some $p$. This system corresponds to the two-dimensional
$(p-1)$-body Toda molecule
\cite{TM1}-\cite{LS} and is sometimes called non-periodic $SL(p)$
Toda chain \cite{OP}. Since it is just a very particular case of the
forced hierarchy, one can make use of the formulas of the previous
section and Appendix 1. Then it is immediate to check that
the $SL(p)$ Toda molecule is described by the arbitrary fermion rotation
matrix $R_{ij}$ only restricted to be of rank $p$ \cite{Mir94}.
This means that the matrix can be presented in the form
\be\label{ATM}
R_{ij}=\sum_k^p f^{(k)}_ig^{(k)}_j
\ee
$f^{(k)}_i$ and $g^{(k)}_j$ being arbitrary coefficients and
the kernel of the Grassmannian element
$A(z,w)$ (\ref{B8}) (whose power series expansion into powers
$z$ and $w$ gives the matrix $R_{ij}$) being of the form
\be\label{kernel}
A(z,w)=\sum_k^p f^{(k)}(z)g^{(k)}(w)
\ee
where $f^{(k)}(z)$ and $g^{(k)}(z)$ are arbitrary functions.

The simplest way to check this statement is to use the determinant
representation (\ref{B27}). In fact, from the first equation of the
2DTL hierarchy
\be\label{Hirota2}
\tau_n\d_{ t_1}\d_{\bar
t_{1}}\tau_n- \d_{ t_1}\tau_n \d_{ \bar t_{1}}\tau_n= \tau_{n+1}\tau_{n-1}
\ee
and condition (\ref{molecule}), one can establish that
$\log\tau_0$ and $\log\tau_p$ satisfy the free wave equation
\be\label{wave}
\partial_{ t_1}\partial_{\bar t_{1}}\log \tau_0=
 \partial_{ t_1}\partial_{\bar t_{1}}\log \tau_{p}=0
\ee
Since the relative normalization of $\tau_n$ is not fixed yet,
one can always choose $\tau_0=1$. Then,
\be \label{mol}
\tau_{0}(t)=1\;\;,\;\;\;\;\;
\tau_{p}(t) = \chi(t_1)\bar \chi(\bar t_{1})
\ee
where $\chi(t_1)$ and $\bar \chi(\bar t_{1})$ are arbitrary functions.
The two-dimensional Toda system with the boundary conditions
(\ref{mol}) has been considered in \cite{TM1}.
Solution to the equation (\ref{Hirota2}) is, in this case,
\cite{LS}:
\be\label{C1}
\tau_{n}(t)\;=\;
        \det\;\partial_{t_1}^{i-1} (-\partial_{\bar t_{1}})^{j-1}\tau_{1}(t)
\ee
with
\be\label{C2}
\tau_{1}(t) = \sum_{k=1}^{p} a^{(k)}(t)\bar a^{(k)}(\bar t_{1})
\ee
where functions
$a^{(k)}(t)$ and $\bar a^{(k)}(\bar t_{1})$ satisfy the conditions
\be
\det\partial^{i-1}_{t_1}a^{(k)}(t)=\chi(t)\; , \;\;
\det (-\partial_{
\bar t_{1}})^{i-1}\bar a^{(k)}(\bar t_{1})=\bar \chi(\bar t_{1})
\ee
This result corresponds to formula (\ref{B27}) with the kernel $A(z,w)$ of
the form (\ref{kernel}). Let us draw that, although the Toda molecule looks
like a forced hierarchy with one additional projector, it is
generically described by the infinite number of the fermion modes because of
infinitely many non-vanishing elements of the matrix (\ref{ATM}).

Note that the function $\tau_p$ can be put equal to unity. This is
exactly the choice corresponding to the group $SL(p)$, while the general
solution is associated with the group $GL(p)$.

Let us discuss now how one can observe this group structure in the Toda
molecule. For the sake of simplicity, we consider below the instance
of the Liouville system ($p=2$), the extension to the higher groups
being quite immediate. The simplest way of doing is to restore the manifest
group structure in (\ref{C1})-(\ref{C2}). This has been done in
\cite{LS1}. Moreover, this representation admits the deformation so
that the quantum Liouville system is merely described by the quantum group
$SL_q(2)$ \cite{Babelon1,J}. Still, it is more convenient for
us to employ here a little bit different approach to the problem \cite{Mir94}.

Namely, we are going to look at the Liouville theory as at a reduction of the
general two-dimensional Toda system associated with the group $GL(\infty)$.
In this case, the functions $f^{(k)}(z)$ and $g^{(k)}(z)$ (either the set of
the corresponding coefficients $f^{(k)}_i$ and $g^{(k)}_i$) describe
the way for the group $SL(2)$ associated with the Liouville system to
embed into the group $GL(\infty)$ in the course of reduction. Thus, the
different embeddings are related by external $GL(\infty)$-automorphisms of
the group $SL(2)$.

Let us describe it in more explicit terms. The Liouville equation is
of the form:
\be\label{Liouv}
\partial \tau_{1} \bar\partial \tau_{1} -
\tau_{1}\partial\bar\partial\tau_{1} =\tau_0\tau_2 = 1
\ee
or
\be
\partial\bar\partial \phi = 2e^{\phi}, \ \
\tau_{1} = e^{-\phi/2}
\ee
Its general solution
\be\label{Louvsol}
\tau_{1}(t,\bar t|g)=(1+A(t)B(\bar t))\left[{\partial
A\over\partial t}{\partial B\over\partial\bar t}\right]^{-\2}
\ee
is parametrized by two arbitrary functions $A(t)$ and $B(\bar t)$.
These functions are related with the fermion rotation matrix via the
following formulas that can be obtained by comparing (\ref{B27})
and (\ref{Louvsol}):
\be\label{at}
G\psi_i G^{-1}=\left(\int\int dt{e^{-tx}\over \sqrt{\partial A}} dxx^{i-1}
\right)\cdot\sum_k\left(\int \int
d\bar t{e^{-\bar ty}\over\sqrt{\bar\partial B}}
dyy^{k-1}\right)\psi_k+\\+\left(\int \int dt{e^{-tx}A\over \sqrt{\partial A}}
dxx^{i-1}
\right)\cdot\sum_k\left(\int \int
d\bar t{e^{-\bar ty}B\over\sqrt{\bar\partial B}}
dyy^{k-1}\right)\psi_k
\equiv f_i\Psi^{(1)}+g_i\Psi^{(2)}
\ee
Thus, the element $G$ of the Grassmannian, which is associated
with the Liouville system rotates the fermionic modes in the
two-dimensional space given by the "dressed" fermions
$\Psi^{(1,2)}_i$ that evidently depend on the concrete choice of
$G$, i.e. on the functions $A$ and $B$. Each choice of
$G$ fixes such a pair of the fermions so that they give the embedding of the
group $SL(2)$ into the group $GL(\infty)$ and all the variety of solutions to
the Liouville equation is given by all possible choices of the functions $A$
and $B$, i.e. by all possible embeddings of this type.

\section{Generalized $\tau$-function and bilinear identities}
\setcounter{equation}{0}
Our previous discussion of $\tau$-functions was mainly concentrated on their
determinant representations. However, the $\tau$-function can be also
introduced by the hierarchy of defining equations, namely, by an (infinite)
set of bilinear identities of the Hirota type (\ref{Hirota2}) and
(\ref{KPHirota}). We show in this section how the $\tau$-function can be
generalized still to satisfy an hierarchy of bilinear identities, and in the
next section we discuss its determinant representations that turn out to be
more particular property.

The main idea of our approach is to formulate any integrable hierarchy in
pure algebraic terms so that it could be constructed for any (highest weight)
representation of any group. In particular, quantizing means substituting the
group by its quantum counterpart.

\subsection{$\tau$-function and bilinear identities}
\paragraph{\b-function.}
Thus, for a given universal enveloping algebra (UEA) $U({\cal G})$
(over some ring) and a Verma module $V$ of this algebra, we define
the \b-function as a generating function of all the matrix elements
$\langle m | g | \bar m \rangle_V$ \cite{GKLMM}:
\be\label{gtau}
\tau_V(t,\bar t|g) \equiv
\sum_{m,\bar m\in V}
s^V_{m,\bar m}(t,\bar t)
\left.\langle m | g | \bar m
\rangle\right._V
\ee
The main ambiguity of this definition is the choice of the function
$s^{V}_{m,\bar m}(t,\bar t)$, which is to be done in a "clever" way
to have $\tau$-function with good properties. The ambiguity is partially
eliminated for the highest weight representations, when one can require
\be
\tau_V (t,\bar t|g)  =
<0_V| U(t) g \bar U(\bar t)| 0_V >
\label{etau}
\ee
with some operators $U$ and $\bar U$ that does not depend on
$V$, $|0_V>$ being the highest weight vector\footnote{Evolution operators
$U(t)$ ($\bar U(\bar
t)$) for the general representation depend on all raising
(lowering) generators of algebra, which generally do not commute. Therefore,
the time flows in the hierarchy equations described below should not commute
in contrast to the case of the standard integrable hierarchies. Still these
evolution operators can commute in the specially choosen representations, for
instance, in the fundamental representations of the classical groups -- see
section 6 for more details.}.

This requirement naturally arises if one considers the $\tau$-function of the
2DTL hierarchy that corresponds to $\G = SL(p)$ with $V$ being one of
the $p-1$ fundamental representations. We will return to this case in
section 5.

However, the set of bilinear identities that we are going to derive now does
not depend on the concrete choice of the functions
$s^{V}_{m,\bar m}(t,\bar t)$, thus, we will return to their explicit
form later.

\paragraph{Intertwining operators and bilinear identities.}
Let us generalize now the derivation of the BI of section 2 to the case of
the general $\tau$-function (\ref{gtau}). Its part that can be formulated
in purely algebraic terms is the derivation of the counterpart of
(\ref{BIglinfty}). The procedure of derivation consists now of
some steps.

1. The starting point is to embed a Verma module
$\widehat V$ into the tensor product $V\otimes W$, where
$W$ is some (arbitrary) finite-dimensional representation of
${\cal G}$.\footnote{In the affine case, one should use finite-dimensional
(evaluation) representations as well. These are representations with
zero central charge, and they are not highest weight representations
-- cf. the definition of vertex operators in
\cite{FR,JM}.} With the fixed choice of $V$ and $W$, there exist only
finite number of possible $\widehat V$.

Define now some generalization of the fermions -- intertwining operators.
The right \iw of the type $W$ is defined to be a homomorphism of the
${\cal G}$-modules:
\be\label{inttw}
E_R: \ \ \widehat V \longrightarrow V\otimes W
\ee
This \iw can be constructed for the highest weight vectors
\be\label{vacuum2}
|{\bf 0} \rangle_{\widehat V} =
\left( \sum_{\{p_{\bf\alpha},i_{\bf\alpha}\}}
 A\{p_{\bf\alpha},i_{\bf\alpha}\}
 \left(\prod_{{\bf\alpha}>0}
 \left. (T_{-{\bf\alpha}}\right)^{p_{\bf\alpha}}\otimes
 \left. (T_{-{\bf\alpha}}\right)^{i_{\bf\alpha}}\right) \right)
| {\bf 0} \rangle_V \otimes | {\bf 0} \rangle_W
\ee
and then continued explicitly to the whole representation in accordance
with formula:
\be\label{vacuum}
\widehat V = \left\{ | {\bf n_{\bf\alpha}} \rangle_{\widehat V} =
 \prod_{{\bf\alpha}>0}
 \left(\Delta (T_{-{\bf\alpha}}\right)^{n_{\bf\alpha}}
 | {\bf 0} \rangle_{\widehat V} \right\}
\ee
where $T_{-{\bf\alpha}}$ are (lowering) generators of algebra, which belong
to its maximal negative (right) nilpotent subalgebra
$\bar N({\cal G})$, $\balpha$ are the positive roots, the vacuum state is
annihilated by all the raising operators $T_{{\bf\alpha}}$ from the positive
(left) nilpotent subalgebra $N({\cal G})$, the Verma module is built by the
action of all generators $T_{-{\bf\alpha}}$ on the highest weight state:
$V = \left\{|{\bf n}_{\bf\alpha} \rangle_V=
\prod_{{\bf\alpha} > 0}\left(T_{-{\bf\alpha}}\right)^{n_{\bf\alpha}}| {\bf
0}\rangle_V\right\}$. Formula (\ref{vacuum}) reflects the fact that
the action of ${\cal G}$ on the tensor product of representations (Verma
modules) is defined by the co-multiplication $\Delta$ and, for finite
$W$, allows one to present every $| {\bf
n_{\bf\alpha}} \rangle_{\widehat V}$ as a {\it finite} sum of states
$| {\bf m_{\bf\alpha}} \rangle_{V}$ with the coefficients taking values in
$W$.

2. As the next step, we consider another triple of modules that define
the left \iw
\be\label{inttw2}
E'_L: \ \ \widehat V' \longrightarrow W' \otimes V'
\ee
so that the product $W\otimes W'$ contains the {\it
unit} representation of ${\cal G}$.

\paragraph{BI in terms of UEA.}
We turn now to the immediate derivation of the BI. They can be obtained
in two forms -- in the operator form that is analogous to
(\ref{5s1}) and as a relation for the matrix elements (i.e. in terms of
the algebra of functions on the group). The first form can be obtained
in two steps. As the first step, one needs to consider the projection
of the product $W\otimes W'$ onto the unit representation
\be\label{proj}
\pi: \ \ W\otimes W' \longrightarrow I
\ee
that can be manifestly provided by multiplying any element from
$W\otimes W'$ and
\be\label{explproj}
\pi = \left._W \langle {\bf 0} | \otimes \left._{W'} \langle {\bf 0} |
\left( \sum_{\{i_{\bf\alpha},i'_{\bf\alpha}\}}
A_{\pi}\{i_{\bf\alpha},i'_{\bf\alpha}\}
 \left(\prod_{{\bf\alpha}>0}
 \left. (T_{+{\bf\alpha}}\right)^{i_{\bf\alpha}}\otimes
 \left. (T_{+{\bf\alpha}}\right)^{i'_{\bf\alpha}}\right) \right)
\right.\right.
\ee
Using this projection, one can construct the new \iw
\be\label{Gamma}
\Gamma: \ \
\widehat V \otimes \widehat V' \stackrel{E_R\otimes E_L'}{\longrightarrow}
V \otimes W \otimes W' \otimes V'
\stackrel{I\otimes \pi \otimes I}{\longrightarrow} V \otimes V'
\ee
possessing the property (\ref{5s1})
\be\label{CRGamma}
\Gamma (g\otimes g) = (g\otimes g) \Gamma
\label{Ggg=ggG}
\ee
for any group element $g$ such that
\be\label{grel}
\Delta(g)=g\otimes g
\ee

Put it differently, the space $W\otimes W'$ contains the canonical
element of pairing $w_i\otimes w^i$ commuting with the action of
$\Delta (g)$. This means that the operator
$\sum_i E_i\otimes E^i :  V\otimes V' \longrightarrow \widehat V\otimes
\widehat V'$ ($E_i\equiv E(w_i),\ E^i\equiv E(w^i)$) commutes with $\Delta
(g)$.

The identity (\ref{Ggg=ggG}) is the algebraic formulation of the BI. In order
to obtain differential (or difference) identities, one needs (like it was in
the above example of the 2DTL hierarchy) to use the definition (\ref{etau})
and consider the matrix element of the identity (\ref{Ggg=ggG}) along with
the evolution operators between the highest weight states.  This gives rise
to identities (there is a lot of equivalent BI in accordance with a lot of
possible choices of the representations $V$ and $\widehat V$) for the objects
like FBA, i.e. for the averages (\ref{etau}) with additional insertions of
the \iws $E_i$ and $E^i$. Then, using the commutation relations of the \iws
with algebra generators, one can push $E_i$ to the highest weight vector. The
result of this procedure can be imitated by the action of some differential
or difference operators that leads to the differential or difference BI. This
latter step is, however, not always possible and requires the correct choice
of the evolution operators $U(t)$, $\bar U(\bar t)$. Up to now, there is no
general recipe for such a choice. However, in the next two subsections we
construct some explicit examples that illustrate how the choice is usually
made. Certainly, the concrete differential (difference) form of the BI
depends on the concrete algebra and its concrete representation and can not
be obtained in a general algebraic form. Moreover, we will show soon that the
same system can be described by different differential hierarchies -- or even
by  difference and differential hierarchies -- depending on the choice of
the evolution operators.

\paragraph{BI in terms of the algebra of functions.}
To conclude the subsection, we obtain the BI in terms of matrix elements,
i.e. in terms of the algebra of functions on the group. This language is dual
to that of UEA.

Thus, let us write down the matrix element (\ref{Ggg=ggG}) taken in between
four states
\be\label{CRGamma2}
\left._{V'} \langle k' | \left._V \langle k |
(g\otimes g) \Gamma | n \rangle_{\widehat V}
 |n' \rangle_{\widehat V'} \right.\right. =
\left._{V'} \langle k' | \left._V \langle k |
\Gamma (g\otimes g) | n \rangle_{\widehat V}
|n' \rangle_{\widehat V'} \right.\right.
\ee

Action of the operator $\Gamma$ can be presented by the formula
\be\label{me2}
\Gamma | n \rangle_{\hat\lambda} | n' \rangle_{\hat\lambda'}
= \sum_{l,l'} | l \rangle_\lambda | l' \rangle_{\lambda '}
\Gamma(l,l'| n,n')
\ee
i.e. (\ref{CRGamma2}) goes into
\be\label{me3}
\sum_{m,m'} \Gamma (k,k'| m,m')
\frac{|| k ||^2_\lambda
 || k' ||^2_{\lambda '}}
 {|| m ||^2_{\hat\lambda}
 || m' ||^2_{\hat\lambda '}}
\langle m | g | n \rangle_{\hat\lambda}
\langle m' | g | n' \rangle_{\hat\lambda '}
= \sum_{l,l'}
\langle k | g | l \rangle_\lambda
 \langle k' | g | l' \rangle_{\lambda '}
\Gamma(l,l'| n,n')
\label{bilimat}
\ee
To rewrite this expression as some differential or difference equation, one
needs to use formula (\ref{gtau}) for the \b-function. Then, one can get the
generating equation for the identities (\ref{bilimat}) making use of the
explicit formulas for the matrix elements
$\Gamma(l,l'| n,n')$ that can be calculated in the group theory framework
(they are nothing but the corresponding Clebsh-Gordon coefficients). Again
to present the generating equation in the differential form, one needs to
choose correctly the coefficients $s^{R}_{m,\bar m}(t,\bar t)$ in formula
(\ref{gtau}). Technically, the simplest way to make this choice is
often to use concrete representation, i.e. BI formulated in terms
of the algebra of functions. We illustrate this method in the next
section.

To conclude this subsection, note that, in order to reproduce literally
the derivation of the BI for the 2DTL hierarchy, one would introduce
for {\it the same} triple of representations $V$, $V'$, $W$ the
conjugated pair of the \iws
\be
\Phi:\ \hat V\otimes W\to V,\ \ \ \ \Phi^{\star}:\ V\to
W\otimes \hat V
\ee
that, by definition, satisfy the conditions
\be
\Delta (g)\Phi=\Phi g,\ \ \ \Phi^{\star}\Delta (g)=g\Phi^{\star}
\ee
since being homomorphisms, and the algebra action in the tensor product is
given by the co-multiplication. These conditions are the literal
generalization of formula (\ref{A16}), and the canonical element commutes
again with the group element $g$ giving rise to the equation (\ref{Ggg=ggG}).
Still, the proposed derivation of the BI using the pair of
representation triples looks more general, and, hence, we use it throughout
this section.

\subsection{Example -- $SL_q(2)$}
To illustrate the above described rather abstract construction of the
generalized $\tau$-function and BI, consider now as an example the
case of the quantum group $SL_q(2)$ \cite{GKLMM}. This example, at the same
time, illustrates what are the $\tau$-function and BI in quantum integrable
systems. It was already noted that, in the framework under consideration,
quantizing just implies substituting group by its quantum counterpart. Note
that, from now on, by $q$-deformation we assume the real
$q>1$, although other values of $q$ can be also easily treated by the
methods developed.

\paragraph{Bilinear identities.}
Algebra $U_q(SL(2))$ \footnote{All the necessary facts of the quantum groups
that we use below can be found in the review \cite{Dem}.} is given
by the generators $T_+$, $T_-$ and $T_0$ subject to the commutation
relations
\be
q^{T_0} T_\pm q^{-T_0} = q^{\pm 1}T_\pm, \ \
\phantom. [T_+,T_-] = \frac{q^{2T_0}-q^{-2T_0}}{q - q^{-1}}
\ee
and by the co-multiplication
\be\label{coprod}
\Delta(T_\pm) = q^{T_0} \otimes T_\pm + T_\pm \otimes q^{-T_0}, \ \
\Delta(q^{T_0}) = q^{T_0}\otimes q^{T_0}
\ee
The Verma module $V_\lambda$ of the highest weight
$\lambda$ (not obligatory half-integer) consists of the elements
\be\label{notation1}
| n \rangle_\lambda \equiv T_-^n|0 \rangle_\lambda, \ \ n\geq0
\ee
such that
\be\label{notation}
T_- | n \rangle_\lambda = | n+1 \rangle_\lambda, \ \
T_0 | n \rangle_\lambda =
 (\lambda - n) | n \rangle_\lambda , \ \
T_+ | n \rangle_\lambda \equiv
b_n(\lambda) | n-1 \rangle_\lambda\\
b_n(\lambda) = [n]_q[2\lambda + 1 - n]_q, \ \ \
 [x]_q \equiv \frac{q^x - q^{-x}}{q - q^{-1}},\ \ \ [n]_q!\equiv
[1]_q[2]_q\ldots[n]_q \\
|| n ||^2_\lambda \equiv \left._\lambda\langle n | n \rangle\right._\lambda
= \frac{[n]_q!\ \Gamma_q(2\lambda +1)}{\Gamma_q(2\lambda +1-n)}
\stackrel{\lambda \in {\bf Z}/2}{=}
\frac{[2\lambda]_q![n]_q!}{[2\lambda -n]_q!}
\ee
Now we obtain BI working in terms of the algebra of functions on the quantum
group, i.e. manifestly calculating matrix elements of the operator $\Gamma$.

As the module $W$ we choose the spin $\2$ irrep of $U_q(SL(2))$. Then,
$\widehat V=V_{\lambda\pm {\f 2}}$, $V=V_\lambda$.

To calculate matrix elements of $\Gamma$, one needs to project
the tensor product of two different $W$ onto the singlet state $S =
|+\rangle |-\rangle - q|-\rangle|+\rangle$:
\be\label{qdethaha}
(A| + \rangle + B| - \rangle)\otimes
(| + \rangle C + | - \rangle D) \longrightarrow
AD-qBC
\ee

With our choice of $W$, one needs to consider two different cases:

($A$) $\widehat
V=V_{\lambda-{\f 2}}$ and $\widehat
V'=V_{\lambda'-{\f 2}}$, or

($B$) $\widehat
V=V_{\lambda-{\f 2}}$ and $\widehat
V'=V_{\lambda'+{\f 2}}$.

Using formulas (\ref{vacuum}) and (\ref{coprod}), one now easily
obtains matrix elements of the projection operator\footnote{In order to
simplify the formulas, hereafter we omit the sign of tensor product from the
notation of the states $|+\rangle
\otimes|0\rangle_\lambda$ etc.}:

\underline{Case A}:
\be\label{projA}
| n \rangle_{\lambda - \frac{1}{2}}
| n' \rangle_{\lambda ' - \frac{1}{2}}
\longrightarrow q^{\frac{n'-n-1}{2}} \left(
[n'-2\lambda']_qq^{\lambda '}
| n+1 \rangle_\lambda | n' \rangle_{\lambda '} -
[n-2\lambda]_qq^{-\lambda}
| n \rangle_\lambda | n'+1 \rangle_{\lambda '}\right)
\ee

\underline{Case B}:
\be\label{projB}
| n \rangle_{\lambda + \frac{1}{2}}
| n' \rangle_{\lambda ' - \frac{1}{2}}
\longrightarrow q^{\frac{n'-n-1}{2}} \left(
[n'-2\lambda']_qq^{\lambda '}
| n \rangle_\lambda | n' \rangle_{\lambda '} -
[n]_qq^{+\lambda +1}
| n-1 \rangle_\lambda | n'+1 \rangle_{\lambda '}\right)
\ee

Using explicit formulas (\ref{projA})-(\ref{projB}) for the matrix
elements of $\Gamma(l,l'| n,n')$, one can find for the identities
(\ref{bilimat}) a simple generating equation
(see the end of the previous subsection) provided the evolution operators are
choosen to be
$q$-exponential\footnote{$e_q(x)\equiv \sum_{n\ge 0} {x^n\over [n]_q!}$.}
of the generators $U(t)=e_q(tT_+)$, $\bar U(\bar t)= e_q(\bar t T_-)$. This
generating equation is of the form

\underline{Case A}:
\be\label{qeqA}
\sqrt{M_{\bar t}^- M_{\bar t'}^+}
\left( q^{\lambda '} D_{\bar t}^{(0)}
                   \bar t'D_{\bar t'}^{(2\lambda ')} -
  q^{-\lambda } \bar t D_{\bar t}^{(2\lambda)} D_{\bar t'}^{(0)}\right)
\tau_\lambda (t,\bar t | g)\tau_{\lambda '}(t',\bar t'| g) = \nn \\
= [2\lambda]_q[2\lambda']_q \sqrt{M_{t}^- M_{t'}^+}
\left( q^{-(\lambda+\frac{1}{2})}t' - q^{(\lambda' + \frac{1}{2})}t\right)
\tau_{\lambda - \frac{1}{2}}(t,\bar t | g)\tau_{\lambda ' - \frac{1}{2}}
(t',\bar t'| g)
\ee
Here $D_t^{(\alpha)} \equiv
\frac{q^{-\alpha}M^+_t - q^\alpha M^-_t}{(q - q^{-1})t}$, and
$M^{\pm}$ are the multiplicative shift operators:
$M^{\pm}_tf(t)=f(q^{\pm 1}t)$.

\underline{Case B}:
\be\label{qeqB}
\sqrt{M_{\bar t}^- M_{\bar t'}^+}
\left( q^{\lambda '} \bar t' D_{\bar t'}^{(2\lambda ')} -
  q^{(\lambda +1)} \bar t D_{\bar t'}^{(0)}\right)
\tau_\lambda (t,\bar t | g)\tau_{\lambda '}(t',\bar t'| g) = \nn \\
= \frac{[2\lambda ']_q}{[2\lambda +1]_q} \sqrt{M_{t}^- M_{t'}^+}
\left( q^{\lambda '} tD_{t}^{(2\lambda +1)} -
  q^{\lambda} t' D_{t}^{(0)}\right)
\tau_{\lambda + \frac{1}{2}}(t,\bar t | g)\tau_{\lambda ' -
\frac{1}{2}}(t',\bar t'|g)
\ee

The classical limit of these equations is

\underline{Case A}:
\be\label{eqA}
\left(2\lambda \frac{\partial}{\partial \bar t'} -
    2\lambda' \frac{\partial}{\partial \bar t} +
 (\bar t' - \bar t)\frac{\partial^2}{\partial \bar t\partial \bar t'}
\right)
 \tau_\lambda (t,\bar t | g)\tau_{\lambda '}(t',\bar t'| g)
= 4\lambda\lambda '  (t' - t)
\tau_{\lambda - \frac{1}{2}}(t,\bar t | g)\tau_{\lambda ' - \frac{1}{2}}
(t',\bar t'| g)
\ee

\underline{Case B}:
\be\label{eqB}
\phantom{fhg}\left[(\bar t'-\bar t){\partial\over\partial\bar t'}-2\lambda'
\right]
\tau_\lambda(t,\bar t|g)\tau_{\lambda '}(t',\bar t'|g)=\\=
{2\lambda '\over 2\lambda +1}\left[(t-t'){\partial\over\partial t}-
2\lambda-1\right]\tau_{\lambda+{\f 2}}(t,\bar t|g)\tau_{\lambda '-
{\f 2}}(t',\bar t'|g)
\ee

Thus, we obtain {\it different} BI that are satisfied by {\it the same}
$\tau$-function. In fact, it suffices to use, say, the first equation (case
$A$) to fix the $\tau$-function completely.

\paragraph{Classical limit.}
Before going into further details, consider the case of the classical
group $SL(2)$ with $\tau$-function satisfying the BI (\ref{eqA})-(\ref{eqB}).
One can construct the general solution to these equations, however, it is
even simpler to obtain $\tau$-function immediately from the definition
(\ref{etau}). Taking the representation of arbitrary spin
$\lambda$, one gets the result:
\be\label{tau-sl2-clas}
\tau_\lambda=\left._\lambda\langle 0|e^{tT_-}ge^{\bar tT_+}
|0\rangle\right._\lambda = (a+b\bar t+ct+dt\bar t)^{2\lambda}
\ee
the group element $g$ being parametrized by three parameters:
\be\label{g-clas}
g=e^{x_+T_+}e^{x_0T_0}e^{x_-T_-}
\ee
and
\be\label{abcd-clas}
a\equiv e^{{\f 2}x_0}+x_+x_-e^{-{\f 2}x_0},\ \ b\equiv x_+e^{-{\f 2}x_0},
\ \ c\equiv e^{-{\f 2}x_0}x_-,\ \ d\equiv e^{-{\f 2}x_0}
\ee
i.e.
\be\label{det-clas}
ad-bc=1
\ee

Let us now return to the BI. It was already mentioned that it is
sufficient to look at the equation (\ref{eqA}), since every solution
to this equation satisfies all other BI, say, the equation
(\ref{eqB}), or other equations obtained with other choices of
$V$, $\widehat V$ and $W$. The general solution (\ref{eqA}) is
3-parametric and, certainly, coincides with (\ref{tau-sl2-clas}).

Let us note that (\ref{tau-sl2-clas}) at $\lambda={1\over 2}$
(fundamental representation) has little to do with the general solution
of the Liouville equation (\ref{Louvsol}) that is associated with the
$SL(2)$-reduction (see section 3). However, there is a connection
between them: although the set of solutions of the Liouville equation
(\ref{Liouv}) is far more rich, solutions (\ref{eqA}) are contained in this
set.

In order to understand how the set of solutions of the Liouville equation
should be restricted, one needs to rewrite
(\ref{eqA}) (like it is usually done for Hirota equations) as
a system of differential equations obtained by expanding into
powers of $\epsilon = \frac{1}{2}(t - t')$ and $\bar\epsilon =
\frac{1}{2}(\bar t - \bar t')$. For instance,
at $\lambda = \lambda'$ one gets from
(\ref{eqA}):
\be\label{expeq}
\hbox{coefficient in front of }\epsilon:  \ \ \
\partial \tau_\lambda \bar\partial \tau_\lambda -
\tau_\lambda\partial\bar\partial\tau_\lambda =
2\lambda \tau_{\lambda - \2}^2 \\
\hbox{coefficient in front of }\bar\epsilon: \ \ \
2\lambda\tau_\lambda \bar\partial^2\tau_\lambda =
(2\lambda -1)(\bar\partial\tau_\lambda)^2 \\
\ldots
\ee
If $\lambda = \frac{1}{2}$, the first of these two equations is nothing
but the Liouville equation (\ref{Liouv}):
\be
\partial \tau_{\2} \bar\partial \tau_{\2} -
\tau_{\2}\partial\bar\partial\tau_{\2} =\tau_0^2 = 1
\ee
while the second one
\be
\bar\partial^2\tau_{\2} = 0
\ee
is the condition that strongly restricts the solutions of the Liouville
equation. It implies that the two arbitrary functions $A(t)$ and  $B(\bar t)$
parametrizing these solutions are reduced to linear functions of
times. In terms of Grassmannian, this means that the
$SL(2)$-\b-function considered here corresponds to the embedding of the
$SL(2)$ matrix into the left upper corner of the $GL(\infty)$ matrix -- see
section 3. Certainly, it had to be expected, since
$SL(2)$ \b-function "knows" nothing of the $GL(\infty)$ group. In
particular, in order to quantize the "full-fledged" Liouville equation,
at first, one needs to study the $GL_q(\infty)$ system, and then its
corresponding reduction.

\paragraph{Solutions to the quantum BI.}
Let us discuss now solutions to the quantum BI. The general $c$-number
solution to the equation (\ref{qeqA}) can be manifestly constructed.
Indeed, one can easily check that
\be\label{qsolA}
\tau_\lambda = [\alpha + {\f \alpha}t\bar t]_q^{2\lambda} \equiv
\sum_{i\geq 0} \frac{\Gamma_q(2\lambda +1)}
{\Gamma_q(2\lambda +1 - i)}\frac{\alpha^{2\lambda-2i}(t\bar t)^i}{[i]_q!}
\ee
satisfies the equation (\ref{qeqA}), since
\be
D_{t}^{(0)} [\alpha + {\f \alpha}t\bar t]_q^{2\lambda} =
{\f \alpha}[2\lambda]_q
   [\alpha+ {\f \alpha}t\bar t]_q^{2\lambda -1}\bar t, \nn \\
tD_t^{(2\lambda)} [\alpha + {\f \alpha}t\bar t]_q^{2\lambda} = -\alpha
[2\lambda]_q
   [\alpha + {\f \alpha}t \bar t]_q^{2\lambda - 1}
\ee
However, this gives only the one-parametric set of solutions in contrast to
the classical case. The explanation of this phenomenon is due to the fact
that, of all the elements of $U_q(SL(2))$, only the Cartan one has the
"correct" co-multiplication law (\ref{grel}), while, in the classical case,
there is the three-parametric family of the group elements
(\ref{g-clas}) having such a co-multiplication.

\paragraph{Quantum \b-function.}
However, in the quantum case, there still exists the way to construct
a three-parametric family of solutions. To this end, it suffices to
consider the non-commutative $\tau$-function\footnote{The idea of
introducing non-commutative $\tau$-function has been also proposed
in \cite{J}. Note that there were also other
very interesting attempts in the same direction \cite{qt},
although I do not quite understand their relation
to our approach.}. Indeed, (\ref{gtau}) implies that \b-function takes
its values in the algebra of functions on the quantum group
$SL_q(2)$, i.e. is {\it non-commutative} quantity. For instance,
in the fundamental representation, it is equal to
\be
\tau_{1\over 2}= \langle + | g | + \rangle  + \bar t \langle +
| g | - \rangle
+ t \langle - | g | + \rangle + t\bar t \langle - | g | - \rangle=
 a+b\bar t+ct+dt\bar t
\ee
where the generators $a,b,c,d$ of the algebra of functions
$A(SL_q(2))$ are the elements of the matrix
\beq
{\cal T} = \left(
\begin{array}{cc}
 a& b \\ c& d
\end{array}
\right),\ \ \ \ \ ad-qbc=1
\eeq
with commutation relations given by the equations
${\cal T}{\cal T}{\cal R}= {\cal R}{\cal T}{\cal T}$ \cite{FRT}
\be
ab = qba, \\
ac = qca, \\
bd = qdb, \\
cd = qdc, \\
bc = cb, \\
ad - da = (q - q^{-1})bc
\label{core}
\ee
In order to obtain such a non-commutative $\tau$-function from
(\ref{etau}), one needs to consider $g$ as an element of UEA given over
some {\it non-commutative} ring instead of the complex number field.
This increases the number of the group elements satisfying
(\ref{grel}). This ring is exactly $A_q(SL(2))$ (see the next
subsection).

In order to construct non-commutative $\tau$-function at any
representation of spin $\lambda$, one can decompose this representation
into the representations of spins
$\lambda-{1\over 2}$ and ${1\over 2}$ \cite{GKLMM}:
\be
\phantom._{\lambda }\langle k | g | n \rangle_{\lambda} =
q^{-\frac{k+n}{2}}\left[
\phantom._{\lambda-{\f 2}} \langle k | g | n \rangle_{\lambda-{\f 2}}
\phantom. \langle + | g | + \rangle +
q^{\lambda }[n]_q\phantom._{\lambda-{\f 2}}\langle k | g | n-1
\rangle_{\lambda-{\f 2}} \langle + | g | - \rangle +\right.\\+ \left.
q^{\lambda }[k]_q\phantom._{\lambda-{\f 2}}\langle k-1 | g | n
\rangle_{\lambda-{\f 2}} \langle - | g | + \rangle +
q^{2\lambda}[k]_q[n]_q\phantom._{\lambda-{\f 2}}\langle k-1 | g | n-1
\rangle_{\lambda-{\f 2}} \langle - | g | - \rangle \right]_q
\ee
Recurrently applying this procedure, one gets
\be
\tau_\lambda(t,\bar t|g) =
\tau_{\lambda - \2}(q^{-\2}t,q^{-\2}\bar t | g)
\tau_{\2}(q^{\lambda - \2}t,q^{\lambda - \2}\bar t | g) =  \\
\stackrel{{\hbox{if}}\ \lambda \in \hbox{\bf Z}/2}{=}
\tau_{\2}(q^{\2 - \lambda}t, q^{\2 - \lambda}\bar t|g)
\tau_{\2}(q^{{3\over 2}-\lambda}t, q^{{3\over 2} - \lambda}\bar t | g) \ldots
\tau_{\2}(q^{\lambda - \2}t, q^{\lambda - \2}\bar t | g)
\ee

\subsection{Universal T-operator (group element)}
\paragraph{General construction.}
Let us construct now in more explicit terms the group element given over
the non-commutative ring -- the algebra of functions on the quantum group.
That is, we construct such an element
$g\in U_q({\cal G})\otimes A(\G)$ of the tensor product of
UEA $U_q({\cal G})$ and its dual algebra of functions $A(\G)$ that
\be\label{Tcoprod}
\Delta_U(g)=g\otimes_{U} g \in A(\G)\otimes U_q({\cal G})\otimes U_q({\cal
G})
\ee
To construct this element \cite{FRT,FG,Mir94,R5}, we fix some basis
$T^{(\alpha)}$ in $U_q({\cal G})$. There exists a non-degenerated
pairing between $U_q({\cal G})$ and
$A(\G)$, which we denote $<...>$. We also fix the basis
$X^{(\beta)}$ in $A(\G)$ orthogonal to $T^{(\alpha)}$ w.r.t. this
pairing. Then, the sum
\be\label{grouel}
\hbox{{\bf T}}\equiv\sum_{\alpha}X^{(\alpha)}\otimes T^{(\alpha)}\in A(\G)
\otimes U_q({\cal G})
\ee
is exactly the group element we are looking for. It is called the universal
{\bf T}-matrix (as it is intertwined by the universal
${\cal R}$-matrix) or the universal group element.

In order to prove that (\ref{grouel}) satisfies formula
(\ref{Tcoprod}) one should note that the matrices
$M^{\alpha\beta}_{\gamma}$ and $D^{\alpha}_{\beta\gamma}$ giving
respectively the multiplication and co-multiplication in
$U_q({\cal G})$
\be
T^{(\alpha)}\cdot T^{(\beta)}\equiv M^{\alpha\beta}_{\gamma}T^{(\gamma)},\ \
\Delta(T^{(\alpha)})\equiv D^{\alpha}_{\beta\gamma}T^{(\beta)}\otimes
T^{(\gamma)}
\ee
give rise to, inversely, co-multiplication and multiplication in the
dual algebra $A(\G)$:
\be\label{DM}
D^{\alpha}_{\beta\gamma}=\left<\Delta(T^{(\alpha)}),X^{(\beta)}\otimes
X^{(\gamma)}\right>\equiv \left<T^{(\alpha)},X^{(\beta)}\cdot
X^{(\gamma)}\right>\\
M^{\alpha\beta}_{\gamma}=\left<T^{(\alpha)}T^{(\beta)},X^{(\gamma)}\right>=
\left<T^{(\alpha)}\otimes T^{(\beta)},\Delta(X^{(\gamma)})\right>
\ee
Then,
\be
\Delta_U(\hbox{{\bf T}})=\sum_{\alpha}X^{(\alpha)}\otimes
\Delta_U(T^{(\alpha)})=
\sum_{\alpha,\beta,\gamma}D^{\alpha}_{\beta\gamma}X^{(\alpha)}\otimes
T^{(\beta)}\otimes T^{(\gamma)}=\\=\sum_{\beta,\gamma}X^{(\beta)}X^{(\gamma)}
\otimes T^{(\beta)}\otimes T^{(\gamma)}=\hbox{{\bf T}}\otimes_U\hbox{{\bf T}}
\ee
This is the first defining property of the universal {\bf T}-operator,
which coincides with the classical one. The second property that allows
one to consider {\bf T} as an element of the "true" group is the group
multiplication law $g\cdot g'=g''$ given by the map:
\be\label{grouplaw}
g\cdot g'\equiv \hbox{{\bf T}}\otimes_A
\hbox{{\bf T}}\in A(\G)\otimes A(\G)\otimes U_q({\cal G})
\longrightarrow g''\in A(\G)\otimes U_q({\cal G})
\ee
This map is canonically  given by the co-multiplication and is
again the universal {\bf T}-operator:
\be
\hbox{{\bf T}}\otimes_A\hbox{{\bf T}}=\sum_{\alpha,\beta}X^{(\alpha)}\otimes
X^{(\beta)}\otimes T^{(\alpha)}T^{(\beta)}=\sum_{\alpha,\beta,\gamma}
M^{\gamma}_{\alpha,\beta}X^{(\alpha)}\otimes X^{(\beta)}\otimes T^{(\gamma)}=
\sum_{\alpha}\Delta(X^{(\alpha)})\otimes T^{(\alpha)}
\ee
i.e.
\be
g\equiv\hbox{{\bf T}}(X,T),\ \ g'\equiv \hbox{{\bf T}}(X',T),\ \
g''\equiv\hbox{{\bf T}}(X'',T)\\
\ \ X\equiv\{X^{(\alpha)}\otimes I\}\in A(\G)\otimes I,
\ \ X'\equiv\{I\otimes X^{(\alpha)}\}\in I\otimes A(\G)\\
X''\equiv\{\Delta(X^{(\alpha)})\}\in A(\G)\otimes A(\G)
\ee

\paragraph{T-operator for $SL_q(2)$.}
In order to get more compact formulas, let us redefine the generators
of $U_q(SL(2))$ so that the co-multiplication becomes non-symmetric:
\be
T_+\longrightarrow T_+q^{-T_0},\ \ T_-\longrightarrow q^{T_0}T_-\\
\Delta(T_+)=I\otimes T_++T_+\otimes q^{-2T_0},\ \
\Delta(T_-)=T_-\otimes I+q^{2T_0}\otimes T_-
\ee
This replace results in substituting everywhere
the $q$-number $[n]_q$ for the $q$-number
$(n)_q\equiv {1-q^{2n}\over 1-q^2}$, and the
$q$-exponential $e_q(x)$ -- for ${\cal E}_q(x)\equiv {\f e_q(-x)}=\sum_{k\ge
0} {x^k\over [k]_q!}q^{-k(k-1)/2}=\sum_{k\ge 0}{x^k\over (k)_q!}$. Besides,
the difference operators $D_t^{(\alpha)}= \frac{q^{-\alpha}M^+_t - q^\alpha
M^-_t}{(q - q^{-1})t}$ are replaced by $D_t^{(\alpha)}=
\frac{q^{2\alpha}M^+_t - 1}{(q^2 - 1)t}$.

Fix now the basis $T^{(\alpha)}=T_+^iT_0^jT_-^k$ in $U_q(SL(2))$.
Then, using the co-multiplication for $T^{(\alpha)}$, one can
calculate the matrix $D^{\alpha}_{\beta\gamma}$ (\ref{DM}) and
construct manifest the orthonormal basis $X^{(\alpha)}$:
\be
X^{(\alpha)}={x_+^i\over (i)_{q}!}{x_0^j\over j!}{x_-^k\over (k)_{q^{-1}}!}
\ee
where generating elements $x_{\pm},\ x_0$ produce the Borel Lie algebra
\be\label{xalgebra}
\phantom{.}[x_0,x_{\pm}]=(\log q) x_{\pm},\ \ [x_+,x_-]=0
\ee
Thus,
\be\label{grouelsl2}
\hbox{{\bf T}}={\cal E}_{q}^{x_+T_+}e^{x_0T_0}{\cal E}_{q^{-1}}^{x_-T_-}
\ee

\paragraph{General formula for T-operator.}
In the case of general (simple) Lie algebra, there are different explicit
representations of the {\bf T}-operator. For instance, in \cite{FG},
in accordance with the general construction
(\ref{grouel}), there was considered the so-called PBW-basis
(Poincare-Birkgoff-Witt basis) \cite{Dem}, presented by the ordered monomials
of the algebra generators corresponding to all roots.
In this basis, the group element is represented by the product of the
$q$-exponentials like (\ref{grouelsl2}) given by the Gauss decomposition, and
the generating elements of the algebra of functions satisfy the relations
like (\ref{xalgebra}). In fact, these relations give again a Borel Lie
algebra. It is related to the general structure of the algebra of functions
as a double (see \cite{Mir94} and reference [11] therein).

Another interesting representation for the group element of arbitrary
simple Lie algebra has been found in \cite{R5}. This representation is
remarkable, since it is constructed in the Chevalle basis of UEA, i.e. in
the generators corresponding to the simple roots. To construct the linear
basis in UEA in terms of Chevalle generators for arbitrary simple Lie algebra
is a non-trivial problem, since they satisfy some additional constraints --
Serre relations \cite{Zhel}, and, therefore, are not free generating
elements. However, in the paper \cite{R5}, the group element was constructed
in these terms not basing on expressions like (\ref{grouel}), but by the
immediate check of the relation (\ref{Tcoprod}) taking into account the Serre
relations.

Namely, in \cite{R5} the group element has been presented in the form:
\be
g = g_Ug_Dg_L \\
g_U =\left.\prod_s\right.^<\ {\cal E}_q^{\theta_sT_{i(s)}},   \ \ \
g_L = \left.\prod_s\right.^>\ {\cal E}_{q^{-1}}^{\chi_sT_{-i(s)}}, \ \ \
g_D = \prod_{i=1}^{r_{\G}} e^{\vec\phi\vec H}
\label{genpar}
\ee
where generators $T_i$ correspond only to the {\it simple} roots
$\pm \vec\alpha_i$, $i = 1,\ldots,r_{\G}$. Every root
$\vec\alpha_i$ may appear some times in this product so that there are
different parametrizations of the group element, which depend on
the choice of the set $\{s\}$ and of the map
$i(s)$ of the choosen set to the set of the simple roots. Technically,
the calculations in terms of Chevalle generators are easier because of
the specifically simple co-multiplication law:
\be
\Delta(T_i) = T_i\otimes q^{-2H_i} + I\otimes T_i
\\ \Delta(T_{-i}) = T_{-i}\otimes I + q^{2H_i}\otimes T_{-i} \label{quco}
\ee
Note that formula (\ref{genpar}) is written down for the simple-laced
algebras. Generally, one needs to consider the $q$-exponential with parameter
$q^{||\vec\alpha_i||^2/2}$, not $q$.

The generating elements of the algebra of functions
$\theta, \chi, \vec\phi$ in this parametrization satisfy the
quadratic algebra that is the exponential of the Heisenberg algebra:
\be
\theta_s\theta_{s'} = q^{-\vec\alpha_{i(s)}\vec\alpha_{i(s')}}
\theta_{s'}\theta_s,\ \ \ s<s'\\
\chi_s\chi_{s'} = q^{-\vec\alpha_{i(s)}\vec\alpha_{i(s')}}
\chi_{s'}\chi_s,\ \ \ s<s'\\
e^{\vec\beta\vec\phi}\theta_s = q^{\vec\beta\vec\alpha_{i(s)}}
\theta_s e^{\vec\beta\vec\phi}\\
e^{\vec\beta\vec\phi}\chi_s = q^{\vec\beta\vec\alpha_{i(s)}}
\chi_s e^{\vec\beta\vec\phi}
\label{comrel}
\ee
These relations can be read off from the parametrization (\ref{genpar}) and
formula (\ref{Tcoprod}).

\subsection{\b-function and representations of the algebra of functions}
In accordance with our definition, the non-commutative $\tau$-function is an
element of the algebra of functions on the quantum group. Therefore, the
natural problem is, having some fixed representation of this algebra, to
determine the "value" of the $\tau$-function in this representation.
At the same time, each representation has to be associated with some solution
of the integrable hierarchy. In fact, this is just invariant description for
the group acting in the space of the spectral parameter \cite{Mir96}. Indeed,
any concrete solution of the classical 2DTL hierarchy is described by the
$\tau$-function (\ref{tau}) with some concrete matrix
${\cal G}_{mn}$ that gives the element of the Grassmannian (\ref{elGr}).
In terms of the group element, this means that we fix some trivial
representation of the algebra of functions merely given by $c$-numbers. These
representations exhaust the representations of the algebra on the classical
group, however, in the quantum case, there are non-trivial representations.
Nevertheless, one still should identify every representation of
$A(\G)$ with some solution of the BI hierarchy. Certainly, any reduction of
the integrable system selects some subspace in the solution space, i.e.
restricts the class of representations under consideration, and can be
described, as a rule, by some additional group structure (similar to the Toda
molecule -- see section 3).

Thus, we come to the following general algebraic scheme of constructing
integrable hierarchy:

{\bf For any given UEA $U({\cal G})$, using the procedure described in
the present section, one can introduce the $\tau$-function that satisfies
BI and takes its values in the algebra $A(\G)$ of functions on the group.
Besides, a natural action of UEA on the algebra of functions is
fixed, and any concrete solution of the BI, in the case of the Lie algebra
corresponding to the solution of the classical hierarchy,
is given by fixing the representation of the algebra of functions.}

Let us consider now the quantum $\tau$-function that is an operator. Then,
the natural question is whether one is able to make some $c$-number
quantities of it. This is important, say, to expose the connections between
$\tau$-function and the generating function of correlators in quantum
system. The simplest $c$-number quantity is the "double" generating function
that generates matrix elements of both the representation and
co-representation of UEA (any co-representation of UEA is equivalent, by
duality, to some representation of the algebra of functions). This generating
function depends on four sets of times and should satisfy BI w.r.t. the
indices of both the representation and the co-representation, i.e. describe
some {\it four-dimensional} equation system.

Another $c$-number function is the $\tau$-function itself taken in the
trivial co-representation. This case is closest to the classical case
and, hence, is especially interesting. We discuss this case in detail
later, however, first, we describe briefly the structure of the
co-representations of quantum groups. We consider the simplest case of
$SL_q(2)$ (more general consideration -- see
\cite{Mir94} and reference [11] therein).

Co-representations of UEA $SL_q(2)$ are given by the representations
of the algebra (\ref{xalgebra}). This is Borel algebra, therefore, it
has no non-trivial finite-dimensional irreducible representations
\cite{Zhel}. Its finite-dimensional representations are reducible
but not completely reducible. The irreducible representations have
been considered first in
\cite{WSoi}\footnote{Striktly speaking, in \cite{WSoi} there have
been discussed the $*$-representations, i.e. the representations with
an additional involution.}. To compare with the results of this paper,
rewrite the standard generators of $A(SL_q(2))$ (\ref{core})
in terms of the
algebra (\ref{xalgebra}) (this is a sort of bosonization
of the algebra of functions):
\be
a=e^{{\f 2}x_0}+x_+x_-e^{-{\f 2}x_0},\ \ b= x_+e^{-{\f 2}x_0},
\ \ c= e^{-{\f 2}x_0}x_-,\ \ d= e^{-{\f 2}x_0}
\ee
It is remarkable that these expressions coincide with
(\ref{abcd-clas}), but different $x$ are not commuting.

In these terms, the only two irreducible representations are
the trivial one given by $x_+=x_-=0$, i.e. $ad=1,\ b=c=0$, and
the infinite-dimensional representation that can be given
explicitly by the action on the basis
$\{e_k\}_{k\ge 0}$ \cite{WSoi}
\be
ae_k=(1-q^{2k})^{\2}e_{k-1}\ (ae_0=0),\ \ de_k=(1-q^{2k+2})^{\2}e_{k+1},\ \
ce_k=\theta q^ke_k,\ \ be_k= -\theta^{-1}q^{k+1}e_k
\ee

This structure of representations can be easily generalized to other quantum
groups (of rank $r$), since the algebras of the generating elements $x$ are
always Borel ones. Thus, the whole set of irreducible representations can be
again exhausted by the trivial and infinite-dimensional representations (all
of them being $r$-parametric). The explicit formulas for them can be found in
\cite{Soi} (see also \cite{Mir94} and reference [11] therein).

\paragraph{\b-function in the trivial representation
and difference KOS hierarchy.}
Thus, $\tau$-function in the trivial representation gives some pattern of the
$c$-number function. This function is, certainly, too simple and
flat. However, now we obtain the equation satisfied by the
$\tau$-function in the trivial representation, which has as its solutions
quite non-trivial $c$-number functions.

First note that, even in the case of $SL_q(2)$, there is no naive
determinant representation similar to (\ref{detrep}). Indeed, introduce
(for the sake of brevity, we denote $D\equiv D^{(1)}$):
\be
C^1_1 = \tau_F = a + b\bar t + ct + dt\bar t
\ee
Then,
\be
C^1_2 = D_{\bar t}C^1_1 = b+dt, \ \ \
C^2_1 = D_t C^1_1 = c+d\bar t, \ \ \
C^2_2 = D_{\bar t}D_t C^1_1 = d
\ee
and $C^a_b$ can not be identified with the generating elements
of the coordinate ring $SL_q(2)$ $a$, $b$, $c$, $d$ (\ref{core})
(for instance, $C^1_2 C^2_1 \neq
C^2_1 C^1_2$). Thus, the determinant ${\det_q}C$ is not the adequate
object and, hence, it (either the definition
of $C$) should be modified. Namely, in the case of
$SL_q(2)$, the appropriate formula is of the form
\be\label{sl2hirota}
\tau_{F^{(2)}} = {\rm det}_q g = 1=
C^1_1 C^2_2 - q C^1_2 M^-_{\bar t} C^2_1 =
\tau_F D_tD_{\bar t}\tau_F - q D_{\bar t}\tau_F M^-_tD_t\tau_F
\ee
We will return to the determinant formulas for higher groups
in the next section. Now let us note that the $\tau$-function
in the trivial representation
$A(SL_q(2))$ -- $ad=1,\ b=c=0$ satisfies the simpler equation
(it has been first proposed in \cite{Sat}, and we abbreviate
it as KOS, in accordance with the first letters of the names of
the authors)
\be\label{hirotadif}
\tau_F D_tD_{\bar t}\tau_F - D_{\bar t}\tau_F D_t\tau_F=1
\ee
This equation is of the form close to the equation (\ref{Hirota2}), which is
easily extendable to the $SL_q(n)$-case, in contrast to (\ref{sl2hirota}).
This is not surprising, since the equation (\ref{hirotadif}) is solved by the
$\tau$-function in the trivial representation.

One more argument in favor of the equation (\ref{hirotadif}) is that
it naturally leads to the determinant of the Schur $q$-polynomials
(\ref{qSchur}). Indeed, consider the trivial element $g$
(corresponding to the trivial representation of $A(\G)$).
Then, for the case of $SL_q(n)$, one can get by the direct calculation
(see \cite{Mir94} and reference [5] therein)
\be
\langle k,0,\ldots,0|
{\cal E}_{q^{-1}}^{s_2T_{12}}{\cal E}_{q^{-1}}^{s_3T_{13}}\ldots
{\cal E}_{q^{-1}}^{s_nT_{1n}}\times {\cal E}_q^{\bar s_2T_{21}}
{\cal E}_q^{\bar s_3T_{31}}\ldots
{\cal E}_q^{\bar s_nT_{n1}}|k,0,\ldots,0\rangle=P^{(q)}_k(s\bar s)
\ee
where $\langle k,0,\ldots,0|$ is the symmetric product of the $k$
simplest fundamental representations.

\subsection{Difference KOS hierarchy}
\paragraph{Difference KOS hierarchy from the 2DTL hierarchy.}
Consider the equation (\ref{hirotadif}) more carefully and, in
particular, demonstrate that it can be obtained as the equation
of the classical 2DTL hierarchy with a non-trivial evolution
\cite{MMV}, i.e. that the {\it difference} equation
(KOS) (\ref{hirotadif}) can be obtained in the framework of the
standard {\it differential} ($GL(\infty)$) 2DTL hierarchy.
To this end, let us consider the $\tau$-function\footnote{In the
next section, we show that this $\tau$-function is the generating
function of the matrix elements in fundamental representations.},
defined just by (\ref{A30}) with redefined time flows
\be
C^k_l(s,\bar s) \longrightarrow
{\cal C}^k_l(s,\bar s) = \sum_{i,j} P^{(q)}_{i-k}(s)\
R_{ij}\ P^{(q)}_{j-l}(\bar s)
\ee
where Schur $q$-polynomials are defined by the formula
\be\label{qSchur}
\prod_i {\cal E}_{q^i}(s_iz^i) = \sum_j P^{(q)}_j(s)z^j
\ee
and satisfy the conditions
\be
D_{s_i} P^{(q)}_j(s) = (D_{s_1})^i P^{(q)}_j(s) = P^{(q)}_{j-i}(s)
\ee
Thus, we obtain
\be
D_{s_i}{\cal C}^k_l = {\cal C}^{k +i}_l, \ \
D_{\bar s_i}{\cal C}^k_l = {\cal C}^k_{l+i}
\ee
and
\be
\tau^{(P^{(q)})}_{n} (s,\bar s | g) =
\det_{1\leq k,l \leq n} D_{s_1}^{k -1}
D_{\bar s_1}^{l -1} {\cal C}^1_1(s,\bar s)
\ee
Thus defined $\tau$-function, indeed, satisfies the equations
(\ref{hirotadif}) \cite{Sat,MMV}:
\be
\tau_{k}\cdot D_{s_1}D_{\bar s_1}\tau_{k}-
D_{s_1}\tau_{k}\cdot D_{\bar
s_1}\tau_{k}=\tau_{k-1}\cdot M^+_{s_1}M^+_{\bar s_1}\tau_{k+1}\\
\ldots\ \ \
\ee
where ellipses means the rest of equations of the KOS hierarchy.
The simplest way to prove this equation is to rewrite the $\tau$-function
using the formula
\be
\det D_{s_1}^i D_{\bar s_1}^j C=
q^{-(n-1)(n-2)}(1-q)^{n(n-1)}(t\bar t)^{{n(n-1)\over 2}}
\det_{0\le i,j < n}  (M^+_{s_1})^i (M^+_{\bar s_1})^j C
\label{taudetdif}
\ee
and then apply the Jacobi identity\footnote{The Jacobi identity
is the particular ($p=2$) case of the general identity for the minors
of arbitrary matrix
$$
\sum_{i_p}C_{ri_p}\hat C_{i_1\ldots i_p|j_1\ldots j_p} = \frac{1}{p!}\sum_P
(-)^P
\hat C_{i_1\ldots i_{p-1}|j_{P(1)}\ldots j_{P(p-1)}} \delta_{rj_{P(p)}}
$$
where the sum in the r.h.s. goes over all permutations of $p$ indices,
and $\hat C_{i_1\ldots i_p|j_1\ldots j_p}$ denote the determinant
(minor) of the matrix obtained from $C_{ij}$ by removing rows
$i_1\ldots i_p$ and columns $j_1\ldots j_p$. Using that
$(C^{-1})_{ij} = \hat C_{i|j}/\hat C$, one can rewrite this identity as
$$
\hat C \hat C_{i_1\ldots i_p|j_1\ldots j_p} =
\left(\frac{1}{p!}\right)^2\sum_{P,P'} (-)^P(-)^{P'}
\hat C_{i_{P'(1)}\ldots i_{P'(p-1)}|j_{P(1)}\ldots j_{P(p-1)}}
\delta_{i_{P'(p)}|j_{P(p)}}
$$}.

\paragraph{Fermionic language for the KOS hierarchy.}
To conclude this section, we show how the KOS hierarchy can be
rewritten in the fermionic language. Since the $\tau$-function is
obtained from the Toda hierarchy by the redefinition of times, one
can simply substitute the new times into the $\tau$-function of the
2DTL hierarchy. Indeed, the relation
\be
\prod_{k=1}^\infty {\cal E}_{q^k}(s_k z^k) = \prod_{k=1}^\infty e^{t_kz^k}
\ee
allows one to express times $t$ through $s$
\be
\sum_{k=1}^\infty t_kz^k = \sum_{n,k=1}^\infty
\frac{s_k^n(1-q_k)^n}{n(1-q_k^n)} z^{nk}
\label{tverT}
\ee
so that
\be
P^{(q)}_k(s) = P_k(t)
\ee
Thus, the $\tau$-function can be presented in the form
\be\label{48}
\tau_{n}(s,\bar s|g) = \tau_{n}(t,\bar t|g)\ \stackrel{(\ref{tau})}{\sim} \
\langle n | e^{H(t)} g e^{\bar H(\bar t)}
|n\rangle
\ee
where
\be\label{49}
H(t) = \sum_{n>0} t_nJ_{+n} \ \stackrel{(\ref{tverT})}{=} \
\sum_{n,k=0}^\infty \frac{s_k^n(1-q_k)^n}{n(1-q_k^n)} J_{+nk} \\
\bar H(\bar t) = \sum_{n>0} \bar t_n J_{-n} =
\sum_{n,k=0}^\infty \frac{\bar s_k^n(1-q_k)^n}{n(1-q_k^n)} J_{-nk}
\ee

$\tau$-function of the KOS hierarchy can be also considered as the
$\tau$-function of the 2DTL hierarchy in Miwa variables
\cite{MMV}. To this end, one needs to consider the general Miwa
transform in some special points. That is, using formulas like
\be
t_k = \frac{1}{k}\frac{((1-q)s_1)^k}{1-q^k} =
\frac{1}{k}\sum_{l \ge 0} \left( (1-q)q^l s_1\right)^k
\label{semperMitr}
\ee
one easily gets that the $\tau$-function (\ref{48}) is described by the
following set of the Miwa variables
\be\label{Miwa100}
\left\{ \left.e^{2\pi i a/k}\mu_k q_k^{-l/k}   \right| a=0,\ldots k-1;
\ \ l\geq 0\right\}, \ \ \
\mu_k = \left( (1-q_k)s_k \right)^{-1/k}
\ee

This means that the KOS hierarchy can be considered as the 2DTL hierarchy in
Miwa variables with the special choice (\ref{Miwa100}) of these latter.

\section{Quantum and classical KP hierarchy with different
evolutions}
\setcounter{equation}{0}
\subsection{Structure of fundamental representations}
In this section, we investigate a particular case of the construction
considered above, that is, the integrable hierarchies associated with the
fundamental representations of the groups $SL(p)$ and $SL_q(p)$
\cite{KMM1,KMM2,GKLMM,Mir94}. They are of great importance, since
it is these cases that correspond to the standard 2DTL hierarchy and its
quantum counterpart and, besides, just in these cases there are some
determinant representations for the $\tau$-functions. Note that, in order to
obtain determinant formulas in quantum case, we have to extend slightly the
definition of the $\tau$-function (\ref{etau}) by introducing
non-commutative times.

Strictly speaking, the 2DTL hierarchy is described by the group
$SL(\infty)$. However, we consider the case of general $p$, although
the results almost do not depend on $p$.

We start with describing the fundamental representations of these groups
\cite{Zhel} (see also Appendix 2). The group $SL(p)$ has
$r\equiv\hbox{rank}\ = p-1$ fundamental representations, the simplest one
$F\equiv F_1$ being the $p$-plet containing the states
\be\label{simplfrep}
\psi_i = T_-^{i-1}| 0 \rangle, \ \ i = 1,\ldots,p
\ee
Here the generator $T_-$ is the sum of all $r$ {\it simple} roots
of $SL(p)$:
$T_- = \sum_{i=1}^r T_{-{\bf\alpha}_i}$.  All other fundamental
representations $F_{k}$ can be now constructed as skew degrees of
$F = F_{1}$:
\be\label{frep}
F^{(k)} = \left\{ \Psi^{(k)}_{i_1\ldots i_k} \sim
 \psi_{[i_1}\ldots \psi_{i_k]} \right\}
\ee
$F_{k}$ is given by action of the operators
\be\label{coprodFR}
R_k(T_-^i) \equiv T_-^i\otimes I \otimes \ldots \otimes I +
I\otimes T_-^i \otimes \ldots \otimes I +
I\otimes I \otimes \ldots \otimes T_-^i
\ee
on the highest weight vector. These operators commute with each other.
Given an integer $k$, evidently, there are exactly $k$ independent ones among
them (with $i = 1,\ldots,k$). The fact that all the fundamental
representations are generated by the same generator $T_-$ is a remarkable
property that can serve as a definition of the fundamental representations
and is in charge of all essential features of the classical integrable
hierarchies, in particular, of many commuting Hamiltonians.

Now one can manifestly construct either a pair of the \iws
analogous to those in s.4.2:
\be\label{inttwfrep}
I_{(k)}:\ \ F_{k+1}
\longrightarrow F_{k} \otimes F, \ \ \ I^{\ast}_{(k)}:\ \ F_{k-1}
\longrightarrow F^*\otimes F_{k}\\ {\hbox{so that}} \ \ \ \Gamma_{k|k'}:\ \
F_{k+1}\otimes F_{k'-1} \longrightarrow F_{k} \otimes F_{k'}
\ee
or a pair of the fermionic \iw
\be\label{fermioniop}
\psi^+:\ \ F_1\otimes
F_k\longrightarrow F_{k+1},\ \ \ \psi^-\equiv I_{(k)}
\ee
Here
\be\label{frep2}
F^* = F^{(r)} = \left\{\psi^i \sim \epsilon^{ii_1\ldots i_r}
  \psi_{[i_1}\ldots \psi_{i_r]} \right\}\\
I_{(k)}:\ \ \Psi^{(k+1)}_{i_1\ldots i_{k+1}} =
            \Psi^{(k)}_{[i_1\ldots i_k}\psi^{\phantom{fgh}}_{i_{k+1}]}, \ \ \
I^*_{(k)}:\ \ \Psi^{(k-1)}_{i_1\ldots i_{k-1}} =
            \Psi^{(k)}_{i_1\ldots i_{k-1}i}\psi^i
\ee
and $\Gamma_{k|k'}$ is constructed through the embedding
$I \longrightarrow F \otimes F^*$ induced by the pairing
$\psi_i \psi^i$: the basis in the linear space
$F^{(k+1)}\otimes F^{(k'-1)}$ induced by $\Gamma_{k|k'}$ from the basis in
the space $F^{(k)}\otimes F^{(k')}$ is
$\Psi^{(k)}_{[i_1\ldots i_k} \Psi^{(k')}_{i_{k+1}]i'_1\ldots i'_{k'-1}}$.

Now one can rewrite the operator $\Gamma$ in terms of matrix elements
\be\label{gkdet}
g^{(k)}\left({{i_1\ldots i_k}\atop{j_1\ldots
j_k}}\right) \equiv \langle \Psi_{i_1\ldots i_k} | g | \Psi_{j_1\ldots j_k}
\rangle= \det_{1\leq a,b\leq k} g^{i_a}_{j_b}
\ee
in the following way
\be\label{gkgk}
g^{(k)}\left({{i_1\ldots i_k}\atop{[j_1\ldots j_k}}\right)
g^{(k')}\left({{i'_1\ldots i'_k}\atop
{j_{k+1}]j'_1\ldots j'_{k'-1}}}\right) =
g^{(k+1)}\left({{i_1\ldots i_k[i'_{k'}}
\atop{j_1\ldots j_{k+1}}}\right)
g^{(k'-1)}\left({{i'_1\ldots i'_{k'-1}]}\atop
{j'_1\ldots j'_{k'-1}}}\right)
\ee
This is the manifest form of BI (\ref{Ggg=ggG}) for matrix elements in the
case of fundamental representations, which is identically hold for any
$g^{(k)}$ of the form (\ref{gkdet}). Making use of ``dressing" the matrix
elements of $g^{(k)}$ by the Schur polynomials depending on times,
one can easily get (\ref{A30}) from (\ref{gkdet}) and the 2DTL equation
(\ref{Hirota2}) -- from (\ref{gkgk}). However, later on we always work with
BI in UEA and use the fermionic pair of the \iws (\ref{fermioniop}).

Note that the fundamental representations for the quantum group
$SL_q(p)$ have the same structure with antisymmetrization replaced by
$q$-antisymmetrization. Some details can be found in Appendix 2 and
\cite{GKLMM}.

Let us now turn to more invariant presentation. Namely, note that, as
described above, the fundamental representations are completely generated by
the operators
$\displaystyle{T_{\pm}^{(k)}
=}$\\$=\displaystyle{\sum_{\vec\alpha :  h(\vec\alpha) = k}
T_{\pm\vec\alpha}}$ that are sums of all the generators of $SL(p)$,
associated with the positive/negative roots of ``weight" $k$ (in the first
fundamental representation $F$, these $T_{\pm}^{(k)}$ are the
$p\times p$ matrices with units on the $k$-th upper/lower diagonal and zeroes
wherever else.

It was already mentioned in s.4.1 that the evolution operators
$U(t)$, $\bar U(\bar t)$ in general representation should contain
{\it all} raising, respectively, lowering
operators\footnote{Since $\tau$-function is, by definition, the
generating function of {\it all} the matrix elements.}
and correspond to non-commutative flows. At the same time,
for the fundamental representations $F_n$ they can be choosen in the form
\be
U(t) = \exp \left(\sum_{k\geq 1} t_kT_+^{(k)}\right), \\
\bar U(\bar t) = \exp \left(\sum_{k\geq 1} \bar t_kT_-^{(k)}\right)
\label{Uop}
\ee
The essential property of the operators $T_{\pm}^{(k)}$ is that they
are commuting:
\be
\left[ T_+^{(k)},\ T_+^{(l)} \right] = 0,\ \ \
\left[ T_-^{(k)},\ T_-^{(l)} \right] = 0
\ee
hence, $U(t)$, $\bar U(\bar t)$ (\ref{Uop}) correspond to commuting flows.
Certainly, general evolution operators depend on additional times and
mutually non-commutative generators of the group.

The operators defined in (\ref{Uop}) celebrate the following properties:

(i) $U,\ \bar U\ \in\  SL(p)$;

(ii) more concretely, $U$, $\bar U$ belong to the nilpotent
subgroup $NSL(p)$ of the group $SL(p)$. In fact,
$NSL(p)$ is a subgroup of the Borel subgroup:
$NSL(p) \subset BSL(p) \subset SL(p)$ (in the fundamental representation
$F_1$, the subgroup $BSL(p)$ consists of all the upper-triangle
matrices with unit determinant, while the matrices from
$NSL(p)$ are additionally constrained to have only units on the main
diagonal);

(iii) since co-multiplication is
\be
\Delta(T_{\pm\vec\alpha}) = T_{\pm\vec\alpha}\otimes I +
         I \otimes T_{\pm\vec\alpha}
\label{claco}
\ee
then
\be
\Delta U(t) = U(t) \otimes U(t) =
\left(U(t)\otimes I\right) \left(I\otimes U(t)\right)
\label{clacoU}
\ee
In other words, the evolution operators $U$, $\bar U$ are the group
elements of the subgroup $NSL(p)$.

These properties are rather appealing and it is natural to try to preserve
them in quantization. However, there are two immediate things to be taken
into account. First, there is nothing similar to the operators
$T_{\pm}^{(k)}$ at $q\neq 1$ (at least, formulated independently of
a specific representation $R$). This means that the manifest expressions for
$U(t)$ and $\bar U(\bar t)$ should be very different from (\ref{Uop}).

Second, there is no reasonable notion of the nilpotent subgroup
$N\G_q$ in quantum case: only the quantum deformation of the Borel
subgroup $B\G_q \subset \G_q$ is well defined (see s.4.3). Indeed,
if one chooses as $U(t)$ an object like $g_U$
and as $\bar U(\bar t)$ -- like $g_L$ (see (\ref{genpar})), due to the
absence of the factor $g_D$, $\Delta(g_U) \neq g_U
\otimes g_U, \ \ \ \Delta(g_L) \neq g_L \otimes g_L$, and the nilpotent
subgroup $N\G_q$ is really absent. At the same time, the Borel subgroup
$B\G_q$ exists, since $\Delta(g_Ug_D) = (g_Ug_D)\otimes(g_Ug_D)$).

Despite this problem, we consider $U$ and $\bar U$ as quantities like
$g_U$ and $g_L$ respectively (however, see the end of this section).
By this reason, formula (\ref{clacoU}) in the quantum case requires
some modification. More precisely, instead of (\ref{clacoU}) we obtain
\be
\Delta(U(\xi)) = U^{(2)}_L(\xi) \cdot
U^{(2)}_R(\xi)
\ee
where
\be\label{1.11n}
U(\xi) = \left.\prod_s\right.^<\ {\cal
E}_q\left(\xi_sT_{i(s)}\right)
\ee
\be
U^{(2)}_L = \left.\prod_s\right.^<\ {\cal
E}_q\left(\xi_sT_{i(s)}\otimes q^{-2H_{i(s)}}\right) \neq I\otimes U(\xi)\\
U^{(2)}_R = \left.\prod_s\right.^<\
{\cal E}_q\left(\xi_s I\otimes T_{i(s)}\right) = I\otimes U(\xi)
\ee
This expression is later used to derive the determinant formulas for the
quantum $\tau$-function.

Note that, taking the evolution operators as elements of the decomposition
(\ref{genpar}), we admit the new approach to the generalized
$\tau$-function when not only the group element but also the operators
$U$ and $\bar U$ are elements  of
$U({\cal G})\otimes A(\G)$, i.e. times are {\it non-commutative}
parameters. Exactly this approach allows one to construct determinant
representations for the $\tau$-functions.

Let us point out that, hereafter within this framework, we always understand
by multiplication of the evolution operators
$U$ and $\bar U$ and group element $g$ in the definition of the
$\tau$-function (\ref{etau}) the group multiplication law
(\ref{grouplaw}), that is to say, the elements $\theta,\
\phi$ and $\chi$ of the algebra (\ref{comrel}) in evolution operators commute
with elements of the corresponding algebra in $g$, see s.4.3.

\subsection{Group element parametrization in fundamental
representations}
Now let us return to the general parametrization (\ref{genpar})
of the group element given in terms of Chevalle generators and
consider the group $SL(p)$. Then, the simplest way to choose the map
$i(s)$ is
$$
i(s):\ \ \ 1,2,\ldots,r-1,r;\  1,2,\ldots,r-1;\  1,2,3;\  1,2;\  1
$$
with $s = 1,\ldots,\frac{p(p-1)}{2}$ (dimension of the group $SL(p)$), i.e.
\be
U(\xi) = \prod_{1\leq i \leq p}\prod_{i<j\le p}\exp\left(
\xi_{ij} T_{j-i}\right)
\label{i(s)}
\ee

However, since we are going to discuss the fundamental representation, it
suffices, how one could see in the previous subsection, to consider the
orbits of $SL(p)$ parametrized by only $r$ variables.
Therefore, out next goal is to find an adequate parametrization of these
orbits. In the classical case ($q=1$), there exist, at least, three such
parametrizations considered below \cite{KMM1,KMM2}. However, only one of
them admits a simple quantum deformation but, instead, does not
correspond to (\ref{Uop}). The problem with constructing these
orbits of smaller dimensions is due to necessity of the reduction consistent
with the commutation relations (\ref{comrel}), i.e. due to necessity of
choosing a subclass of rather specific representations of the algebra of
functions on quantum group.

\paragraph{Parametrization A.} The simplest possibility is to
restrict the set $\{s\}$ onto $s = 1,\ldots,r$ and
choose $i(s) = s$, i.e.
\be
U^{(A)}(\xi) = \left.\prod_{i=1}^{r_{\G}}\right.^<\
\exp\left(\xi_iT_i\right)
\label{UA}
\ee
This is sufficient to generate all the states of any fundamental
representation from the corresponding vacuum vector (highest weight vector)
but
$<0_{F_n}|\ U^{(A)}(\xi) $ has little to do with
$<0_{F_n}|\ U(t)$ (with $U(t)$ given by formula (\ref{Uop})).
Maybe it is better to say that the identification
$ <0_{F_n}|\ U^{(A)}(\xi) = <0_{F_n}|\ U(t)$ gives rise to a quite
complicated map $\xi_i(t)$ that manifestly depends on $p$.

One can certainly construct the KP/Toda hierarchy in terms
of variables $\xi$ instead of the usual times $t$ (see below) but
it is {\it not} obtained just as a replace of variables -- the whole
construction looks absolutely different. This is the price for easy
generalization of this construction to the case of $q\neq 1$:
one needs just to write instead of (\ref{UA})
\be
U^{(A)}(\xi) = \left.\prod_{i=1}^{r_{\G}}\right.^<\
{\cal E}_q\left(\xi_iT_i\right)
\label{UAq}
\ee
where $\xi$ are non-commuting quantities,
\be
\xi_i\xi_j = q^{-\vec\alpha_i\vec\alpha_j}\xi_j\xi_i,
\ \ \ i<j
\ee
It is shown below that now any statement at $q=1$ in
$\xi$-parametrization is easily generalized to the case of
$q\neq 1$.

\paragraph{Parametrization B (conventional).} The other choice of
parametrization is related to the conventional times of the KP/Toda
hierarchy:
\be
U^{(B)}(\xi) = \prod_s \exp\left(\xi_sT_{i(s)}\right) =
U(t) = \exp \left(\sum_k t_kT_+^{(k)}\right)
\label{2.16n}
\ee
This implies that $ <0_{F_n}|\ U^{(B)}(\xi) = <0_{F_n}|\ U(t)$
with some $p$-independent functions $\xi_s(t)$. The key difference
between the two sides of the equality (\ref{2.16n}) is that the r.h.s.
contains mutually-commuting combinations of root generators, while the l.h.s.
contains only mutually non-commuting Chevalle generators. Such a
reparametrization actually exists, but the set $\{s\}$ should contain,
at least, $\frac{p(p-1)}{2}$ elements and one can take
$i(s)$ just as in (\ref{i(s)}). However, now not all of the
$\xi_s$'s are independent: instead they are expressed through
$r$ times $t_k$. For instance, the $t_1$-dependence of
$\xi_{ij}$ is given by
\be
\xi_{ij} = \frac{t_1}{p+i-j} + {\cal O}(t_2, t_3,
\ldots)  \label{xivert}
\ee

However, in order to construct some reasonable quantum
deformation of parametrization B, one needs to reproduce
the proper commutation relations
\be
\xi_s\xi_{s'} = q^{-\vec\alpha_{i(s)}\vec\alpha_{i(s')}}
\xi_{s'}\xi_s,\ \ \ s<s'
\ee
for $\frac{p(p-1)}{2}$ variables $\xi_s$ as a corollary of
{\it some} relations between $r$ variables $t_k$
(which, of course, do not commute when $q\neq 1$). To make this possible, one
should also somehow deform formulas (\ref{xivert}) at $q\neq 1$.
Solution to this problem is unknown so far.

\paragraph{Parametrization C (Miwa variables).} Yet another option is to make
the (representation-independent) Miwa transform
$t_k = \frac{1}{k}\sum_a \lambda_a^k$. This transformation is
perfectly consistent with the simple-root decomposition:
\be
U(t) = \prod_a \exp\left(\sum_{k=1}^{r_{\G}}\frac{\lambda_a^k}{k}
T_+^{(k)}\right) =
\prod_a \left(\prod_{i=1}^{r_{\G}} e^{\lambda_aT_i}\right)
\label{claMipa}
\ee
The set $\{s\}$ and mapping $i(s)$ are not of the "most economic" type
(\ref{i(s)}), but the general rule (\ref{comrel}) of the quantum deformation
is, of course, applicable.

In this parametrization, however, the problem is that
(\ref{i(s)}) implies the quantum formula in the
form
\be
\prod_a \left(\prod_{i=1}^{r_{\G}} {\cal E}_q\left(\lambda_{ai}T_i
\right)\right)
\ee
where $\lambda_{ai}$ with different $i$ and the same $a$ do not commute.
At the same time, the constraint
$\lambda_{ai} = \lambda_{aj}$ for $i\ne
j$ is of crucial importance for the classical
($q=1$) formula (\ref{claMipa}). What is the proper quantum deformation of
this constraint remains unclear.

\subsection{Classical ($q=1$) KP/Toda hierarchy in different
parametrizations}
\paragraph{Determinant formulas and system of equations.}
To understand better the scheme developed here, let us rederive now the
standard features of classical integrable hierarchies within the group theory
framework under consideration. These features include determinant
representations for the $\tau$-functions and the set of differential BI. We
consider the fundamental representations of the group $SL(p)$, since
just in this case there exist determinant representations for the
$\tau$-functions. For their derivation, one does not need to specify any
parametrization and can consider some arbitrary evolution
$U(t)$.

Thus, we consider the $\tau$-functions
$\tau_n \equiv \tau_{F_n}(t,\bar
t|g)$. Let us begin with the simplest one
$$
\tau_1 = <0_{F_1}|U(t) g \bar U(\bar t)|0_{F_1}>
$$
Note that the specific feature of $F=F_1$ is
$$
<0_{F}|U(t) =
\sum_k {\cal P}_k(t) <0_{F} | T_+^k =
\sum_k {\cal P}_k(t) <k_{F}|
$$
where the r.h.s. is re-expanded in terms of "generalized Schur polynomials"
(the first equality in this formula defines these polynomials)
and $p$ states of $F = F_1$ are denoted as $<k_F| = <0_F|
T_+^k$, $k = 0,...,r=p-1$. Thus,
\be
\tau_1(t,\bar t|g) =
\sum_{k,\bar k} {\cal P}_k(t){\cal P}_{\bar k}(\bar t)
<0_{F}|T_+^kg T_-^{\bar k}|0_{F}> =\\=
\sum_{k,\bar k} {\cal P}_k(t){\cal P}_{\bar k}(\bar t)
<k_{F}| g |\bar k_{F}>=\sum_{k,\bar k} {\cal P}_k(t)g_{k,\bar k}{\cal P}_{\bar
k}(\bar t)  \label{clatau1}
\ee
One can also define
\be
\tau_1^{m\bar m}
\equiv <m_F|U(t) g \bar U(\bar t)|\bar m_F> = \sum_{k,\bar k} {\cal P}_k(t)
g_{m+k,\bar m+\bar k}{\cal P}_{\bar k}(\bar t)
\label{clataumm1}
\ee

Now we return to the generic fundamental representation
$F_n$. Since
\be
\left.<m_1\ldots m_n\right._{F_n}| =
<{m_1}_F|\otimes<{m_2}_F|\otimes\ldots\otimes<{m_n}_F| +\\
+
{\hbox{antisymmetrization over}\ }m_1,\ldots,m_n =\\
= \sum_P (-)^P <m_{P(1)}|\otimes <m_{P(2)}| \otimes\ldots\otimes
<m_{P(n-1)}|
\label{as}
\ee
the vacuum (highest weight) state $F_n$ can be written as
\be
<0_{F_n}| = <0,1,\ldots,n-1_{F_n}|=\\
= \sum_P (-)^P <P(0)_F|\otimes <P(1)_F| \otimes\ldots\otimes
<P(n-1)_F|
\label{ass}
\ee
Since, for the classical group, (see (\ref{clacoU}))
\be\label{coprodcl}
U(t)|_{F_n} = \Delta^{n-1}U(t) = U(t)^{\otimes n}, \ \ \
g|_{F_n} = \Delta^{n-1}(g) = g^{\otimes n}
\ee
one finally gets
\be
\tau_{n+1}(t,\bar t|g) \equiv <0_{F_n}| U(t) g \bar U(\bar t)
|0_{F_n}> =  \\
= \sum_{P,\bar P} (-)^P(-)^{\bar P}
\prod_{k=0}^{n} <P(k)_F|U(t) g \bar U(\bar t)|\bar P(k)_F> = \\
= \det_{0\leq m,\bar m<n}\tau_1^{m\bar m} = \det_{0\leq m,\bar m<n}
\sum_{l,\bar l} {\cal P}_{l-m}(t)g_{l_F,\bar l_F}{\cal P}_{\bar l-\bar
m}(\bar t)=\\ = \sum_{{1<m_1<m_2<...}\atop{1<\bar m_1<\bar m_2<...}}
\det_{ji}{\cal P}_{m_j-i}(t) \det_{ji}g_{m_{j}\bar m_{i}} \det_{ij}
{\cal P}_{\bar m_i-\bar j}(\bar t)
\label{det}
\ee
This result is to be compared with formula
(\ref{A30}) as $p\to\infty$, with the corresponding Dynkin diagram of the
group $A_{p-1}$ being infinite in both directions in this limit
(semi-infinite for the forced hierarchy).

The next step is to obtain the differential BI. The determinant formulas
(\ref{det}) are, certainly, not the best starting point to this end.
It is far more convenient to apply the methods based on the \iws, which we
develop in this review. In particular, for the fundamental representations
one should consider the fermionic pair $\psi_i^{\pm}$
($i=1\ldots p$) of \iws
$\psi^\pm :\ \ F_1\otimes F_n$~{\ams \symbol{29}}~$F_{n+1}$.
Then, the derivation of BI for arbitrary $p$ almost literally coincides
with that for $p\to\infty$, which was described in s.4.1.
One only needs to introduce the fermions containing the finite number
of degrees of the spectral parameter
$\psi^+ (z)\equiv \sum_{i=1}^p
\psi_i^+ z^i$ and $\psi^{-}(z)\equiv\sum_{i=1}^p
\psi_i^{-}z^{p-i+1}$, and, respectively, to define vertex operators through
the formulas
\be
\sum_i\Psi^{+,i}_k(t)z^{i}\equiv \hat X^{+} (z,t)\tau_n(t),\ \ \
\sum_i\Psi^{-,i}_k(t)z^{p-i+1}\equiv \hat X^{-} (z,t)\tau_n(t)
\label{VO}
\ee
and similarly for ${\hat {\bar X}}^{\pm}(z,t)$. Then BI acquires the form
\cite{KMM1,KMM2}
\be\label{IntHir}
\oint{dz\over z^{p+2}} \hat X^-(z,t)\tau_n(t,\bar t) \hat
X^{+}(z,t')\tau_m(t',\bar t')= \oint{dz\over z^{p+2}}\hat
{\bar X}^-(z,\bar t)\tau_{n+1}(t,\bar t) \hat {\bar X}^{+}(z,\bar
t')\tau_{m-1}(t',\bar t')
\ee

In order to obtain concrete systems of differential equations, one now needs
to consider particular choices of parametrizations of
$U(t)$, $\bar U(\bar t)$.

\paragraph{Conventional parametrization (B).}
This parametrization leads to the standard KP/Toda hierarchy
with determinant representation
(\ref{A30})-(\ref{A31}) for the $\tau$-function. Indeed,
using (\ref{Uop}), one gets
$$
<0_{F}|U(t) = <0_{F} | \exp\left(\sum_k t_kT_+^{(k)}\right)=
\sum_k P_k(t) <k_{F}|
$$
with the conventional Schur polynomials $P_k(t)$ (\ref{Shur}).
The main peculiarity of this evolution is the property
\be
\tau_1^{m\bar m} = \frac{\partial}{\partial t_m}
\frac{\partial}{\partial\bar t_{\bar m}}\tau_1 =
\left(\frac{\partial}{\partial t_1}\right)^m
\left(\frac{\partial}{\partial\bar t_1}\right)^{\bar m}\tau_1
\label{clat-derB}
\ee
In order to obtain the system of equations in this parametrization, one
can note that the vertex operators (\ref{A8})-(\ref{A9}) are
\be
\hat X^{+}(z,t) = \hbox{Pr}_p\left[e^{\xi(z,t)}
\hbox{Pr}\left[z^ne^{-\xi(z^{-1},\tilde
\partial_{t})}\right]\right]\\
\hat X^{-}(z,t) = \hbox{Pr}_p\left[e^{-\xi(z,t)}\hbox{Pr}\left[z^{p-n+1}
e^{\xi(z^{-1},\tilde
\partial_{t})}\right]\right]
\ee
(and similarly for the other pair of the vertex operators), where
$\hbox{Pr}[f(z)]$ projects onto the polynomial part of the function
$f(z)$ and
$\hbox{Pr}_l[f(z)]$ projects onto the polynomial part of degree
$l$.

Substituting these formulas into (\ref{IntHir}) and expanding this latter
into degrees of $t_i-t_i'$, we arrive at the set of equations that gives
equations of the standard 2DTL hierarchy as
$p\longrightarrow \infty$ \cite{DJKM}.

\paragraph{KP/Toda hierarchy in parametrization A.}
Consider now the same hierarchy with a different evolution given by the
parametrization A. From now on, we denote, for brevity,
$\hat U(\xi) \equiv U^{(A)}(\xi)$, and the corresponding $\tau$-function
is $\hat\tau(\xi,\bar\xi|g)$. This $\tau$-function is linear
in each time $\xi_i$ and, therefore, has simpler determinant representation
and satisfies simpler hierarchy. Indeed,
(\ref{clatau1}) now turns into
\be
\hat\tau_1(\xi,\bar\xi|g) \equiv <0_{F_1}|\hat U(\xi) g \hat{\bar U}(\bar
\xi)|0_{F_1}> = \sum_{k,\bar k \geq 0} s_k\bar s_{\bar k} <k|g|\bar k>
\label{clatau1s}
\ee
where $s_k = \xi_1\xi_2\ldots \xi_{k}$, $s_0=1$, and (\ref{clataumm1})
is replaced by
\be
\hat\tau_1^{m\bar m}(\xi,\bar\xi|g)
\equiv <m_{F_1}|\hat U(\xi) g \hat{\bar U}(\bar \xi)|\bar m_{F_1}> =
\frac{1}{s_m\bar s_{\bar m}} \sum_{{k \geq m}\atop{\bar k \geq \bar m}} s_k\bar
s_{\bar k} <k|g|\bar k>   = \\ =\frac{1}{s_{m}\bar s_{\bar m}}\sum_{{k \geq
m}\atop{\bar k \geq \bar m}} \frac{\partial}{\partial\log s_k}
\frac{\partial}{\partial\log \bar s_{\bar k}}\tau_1(\xi,\bar\xi|g)=
\frac{1}{s_{m-1}\bar s_{\bar
m-1}}\frac{\partial}{\partial \xi_m} \frac{\partial}{\partial \bar
\xi_{\bar m}}\tau_1(\xi,\bar\xi|g)
\label{clataumm1s}
\ee
Therefore,
\be
\hat\tau_{n+1} = \det_{0\leq m,\bar m \leq n}\hat\tau_1^{m\bar m}
= \left(\prod_{m=1}^n s_m\bar s_{\bar m}\right)^{-1}
\det_{(m,\bar m)} \left(
\sum_{{k \geq m}\atop{\bar k \geq \bar m}} s_k\bar s_{\bar k}
<k|g|\bar k>\right) = \\
=  \frac{1}{s_n\bar s_n}
\sum_{k,\bar k \geq n} s_k\bar s_{\bar k}
\det_{0 \leq m,\bar m \leq n-1}\left(
\begin{array}{cc}
g_{m\bar m} & g_{m\bar k} \\
g_{k\bar m} & g_{m\bar m}
\end{array} \right)
\equiv \frac{1}{s_n\bar s_n}
\sum_{k,\bar k \geq n} s_k\bar s_{\bar k}
{\cal D}_{k\bar k}^{(n)}
\label{detA}
\ee
One can compare the determinant representations
(\ref{A30})-(\ref{A31}) and
(\ref{detA}) to get the connection between different coordinates
$t$ and $\xi$. This connection possesses the structure
$s_k\sim$ some functions of $P_j(t)$.  \footnote{In the simplest case of the
first fundamental representation, one needs to identify
$\tau_1(t|g)
=\hat\tau_1(\xi|g)$, i.e. $ s_k = P_k(t),\ \ \ {\partial\over\partial
t_k}=\sum_i s_{i-k}{\partial\over \partial s_i}$. However, identification
of $\tau_n(t)$ and $\hat\tau_n(\xi)$ with $n\neq 1$ leads to
different relations between $\xi$ and $t$.}

One can also easily obtain the differential equations for the
$\tau$-function in parametrization A. Indeed, it is straightforward to find
FBA (\ref{BA}) and substitute this into (\ref{BIA})
\be
\Psi^{+,n+k+1}_n(\xi)={s_{n+k}\over
s_{n}}\left(\tau_n(\xi)-\xi_n{\partial\tau_n(\xi)\over\partial\xi_n}\right)\\
\Psi^{+,n+1}_n(\xi)=
\left(\tau_n(\xi)-\xi_n{\partial\tau_n(\xi)\over\partial\xi_n}\right),\ \
\Psi^{+,n}_n(\xi)=-{\partial\tau_n(\xi)\over\partial \xi_n}\\
\Psi^{-,k}_n(\xi)-\xi_{n-1}{\partial\Psi^{-,k}_n(\xi)\over\xi_{n-1}}
={s_{n-1}\over s_{k-1}}{\partial\tau_n(\xi)\over\partial
\log\xi_{k}}+{s_{n-1}\over s_{k-2}}{\partial\tau_n(\xi)\over\partial\xi_{k-1}}
\ \hbox{for}\ k>n\\
\Psi^{-,n}_n(\xi)-\xi_{n-1}{\partial\Psi^{-,n}_n(\xi)\over\xi_{n-1}}
=\tau_n(\xi)+ {\partial\tau_n(\xi)\over \log\xi_{n}}\\
\Psi^{-,n-1}_n(\xi)=\xi_{n-1}\tau_n(\xi),\ \ \Psi^{-,k}_n(\xi)=0\ \hbox{for}\
k<n-1\label{exBAA}
\ee
As for the values of $\Psi^{+,k}_n(\xi)$ for $k<n$, they are constants that
can hardly be expressed as an action of a differential operator on
$\tau_n(\xi)$. This means that formula
(\ref{BIA}) (where manifest expressions for
$\bar\Psi$ analogous to (\ref{exBAA}) can be also easily written down) does
not lead to differential equations when $k$ and $l$ are arbitrary choosen.
If, however, one chooses
$k\le l-1$, because of multiple cancellations due to
(\ref{exBAA}), formula(\ref{BIA}) is almost a differential equation. It can
be easily transformed to a differential equation by putting
$\xi_{n-1}=0$ (see (\ref{exBAA})). One can easily check that the number of
independent equations obtained in this way is sufficient to determine the
$\tau$-function in full. This means that the whole hierarchy can be still
presented in the differential form.

\subsection{Quantum case ($q\neq 1$)}
\paragraph{$q$-determinant-like representation.}
Let us demonstrate now how the technique developed in the previous subsection
is deformed to the quantum case, i.e. the group $SL_q(p)$ and, in
particular, obtain $q$-determinant-like representations for the
$\tau$-functions similar to (\ref{det}). We also demonstrate  that, in
parametrization A, formula
(\ref{clataumm1s}) expressing $\tau_1^{m\bar m}$ through
derivatives of $\tau_1$ is still correct for $q\ne 1$, with all the
derivatives replaced by difference operators.

In order to obtain $q$-determinant-like representations, it suffices to
study
{\it any} $U(\xi)$ of the form (\ref{1.11n}), without reference to the
particular parametrization A like it was in the case of classical
groups.\footnote{Actually, we require the $U(\xi)$ to be an element from
$N\G_q$ and to be expressed only through the generators associated with {\it
simple positive} roots:
$U(\xi)=U\{\xi_s|T_i\}$. Formula (\ref{1.11n}) is a possible but not unique
realization of these requirements.}

As a result of absence of the diagonal factor $g_D$, the classical
co-multiplication law (\ref{coprodcl}) is replaced in quantum case
by the following co-multiplication rule:
\be
\Delta^{n-1}(U\{T_i\})
= \prod_{m=1}^n U^{(m)}
\ee
where
\be
U^{(m)} = U\left\{\ I\otimes\ldots
\ldots \otimes I\otimes T_i \otimes q^{-2H_i}\otimes \ldots\otimes
q^{-2H_i}\right\}
\ee
($T_i$ appears at the $m$-th place in the tensor product). Similarly,
\be
\bar U^{(m)} = \bar U
\left\{\ q^{2H_i}\otimes\ldots \ldots \otimes q^{2H_i}\otimes T_{-i}
\otimes I\otimes \ldots\otimes I\right\}
\ee
Let
$$
H_i|\bar j_{F_1}> = h_{i,\bar j}|\bar
j_{F_1}>, \ \ \ <j_{F_1}|H_i = h_{i,j}<j_{F_1}|
$$
(in fact, for $SL(p)$ $\
2h_{i,i-1} = +1,\ 2h_{i,i} = -1$, and all the rest are vanishing). Then
\be
\tau_n^{j_1\ldots j_n\bar j_1\ldots\bar j_n}(\xi_s,\bar \xi_s|g)
\equiv
\left(\otimes_{m=1}^n<j_m|\right)
\Delta^{n-1}(U)\ g^{\otimes n}\ \Delta^{n-1}(\bar U)
\left(\otimes_{m=1}^n|\bar j_m>\right)=
\\= \prod_{m=1}^n \ <j_m|
U\left\{T_i
q^{-2\sum_{l = m+1}^nh_{i,j_l}}\right\}\ g\
\bar U\left\{T_{-i}
q^{2\sum_{l=1}^{m-1}h_{i,\bar j_l}}\right\}\ |\bar j_m>=\\
=\prod_{m=1}^n \tau_1^{j_m\bar j_m}\left(\xi_s
q^{-2\sum_{l = m+1}^nh_{i(s),j_l}},\bar \xi_s
q^{2\sum_{l=1}^{m-1}h_{i(s),\bar j_l}}\right)
\label{4.42n}
\ee

In order to get a $q$-determinant-like counterpart of (\ref{det}),
one should replace antisymmetrization
by $q$-antisymmetrization in formulas (\ref{as})-(\ref{ass}),
since, in the quantum case, the fundamental representations are
described by $q$-antisymmetrized vectors
(see details in s.5.2 of \cite{GKLMM}
and in Appendix 2). We define the $q$-antisymmetrization
as a sum over all permutations
\be
\left( [1,\ldots,k]_q\right) = \sum_P (-q)^{{\rm deg}\ P}
\left(P(1),\ldots,P(k)\right)
\ee
where
\be
{\hbox{deg}}\ P = \ \hbox{\ \ number\ of\ inversions\ in}\ P
\ee
Then, $q$-antisymmetrizing (\ref{4.42n}) with $j_k=k-1,\ \bar j_{\bar
k}=\bar k-1$, one finally gets
\be
\tau_n(\xi,\bar \xi|g) = \sum_{P,P'} (-q)^{{\rm deg}\ P + {\rm deg}\ P'}
\prod_{m=0}^{n-1}
\tau_1^{P(m)P'(\bar m)}\left(\xi_s
q^{-2\sum_{l = m+1}^{n-1}h_{i(s),P(l)}},\bar \xi_s
q^{2\sum_{\bar l=0}^{m-1}h_{i(s),P'(\bar l)}}\right)
\label{detq}
\ee
If there would be no $q$-factors
twisting the time variables, this would be just a $q$-determinant
which is defined by the formula
\be
{\det}_q A \sim  A^{[1}_{[1}\ldots A^{n]_q}_{n]_q}
= \sum_{P,P'} (-q)^{{\hbox{deg}}\ P + {\hbox{deg}}\ P'}
\prod_{a} A^{P(a)}_{P'(a)}
\label{qdet}
\ee
The concrete example of formula (\ref{detq}) in the simplest non-trivial
case of the second fundamental representation, which makes its structure more
transparent is described in Appendix 3.

Note that the $q$-factors in all the expressions obtained above can be
trivially reproduce by action of the operators
$$
M_j^{\pm}:\ \ M_j^{\pm}\xi_{s} =
q^{\pm\delta_{j,i(s)}}\xi_{s}
$$
$$
\bar M_j^{\pm}:\ \ \bar M_j^{\pm}\bar \xi_{s} =
q^{\pm\delta_{j,i(s)}}\bar \xi_{s}
$$
We already observed some times that these operators typically emerge in
deformations of integrable hierarchies.

Let us now briefly discuss the equations satisfied by the quantum
$\tau$-function. For their derivation, one needs to follow the same
line as in the classical case and introduce \iws. In quantum case, one
should distinguish between the right and left \iws:
$\Phi^{\pm,R} :\ \ F_n\otimes F_1$~{\ams
\symbol{29}}~$F_{n+1}$ and $\Phi^{\pm,L} :\ \ F_1\otimes F_n$~{\ams
\symbol{29}}~$F_{n+1}$. These operators $\Phi^{\pm,R,L}$ can be expressed
through the classical \iws (fermions)\footnote{This is a
corollary of similarity of the irreducible representations of UEA's
in the classical and quantum cases (at $|q|\ne 1$) \cite{Ros}.}:
\be
\Phi^{\pm,R}_i=q^{-\sum_{j=1}^{i-1}\psi^+_j\psi^-_j}\psi^{\pm}_i,\ \ \
\Phi^{\pm,L}_i=q^{\sum_{j=1}^{i-1}\psi^+_j\psi^-_j}\psi^{\pm}_i
\label{phipsi}
\ee
In analogy with the classical case, one needs to consider now
the operator $\Gamma=\sum_i\Phi^{+,L}_i \otimes\Phi^{-,R}_i$
commuting with
$g\otimes g$.
Then, introducing quantum FBA as averages of the quantum \iws
$\Phi^{\pm}$ (properly labeled by indices $L$ and $R$) and corresponding
vertex operators, one obtains the same equations
(\ref{BIglinfty})-(\ref{IntHir}) but with re-defined entries.
Technically, the action of vertex operators can be calculated with the help
of formula (\ref{phipsi}) and leads to quantum counterparts of formulas
(\ref{main1})-(\ref{main2}).

Now we fix the concrete parametrization A and show how all these formulas
work.

\paragraph{Parametrization A.}
Most of expressions of the previous subsection remain almost the same
in this parametrization (although, in the quantum case, one needs
to take care of normal ordering of different objects). In particular,
\be
\hat\tau_1(\xi,\bar\xi|g) \equiv\ <0_{F_1}|\hat U(\xi) g \hat{\bar U}(\bar
\xi)|0_{F_1}>\ = \sum_{k,\bar k \geq 0} s_k\bar s_{\bar k} <k|g|\bar k>
\label{clatau1sq}
\ee
where again $s_k = \xi_1\xi_2\ldots \xi_{k}$, $s_0=1$, while
$\bar s_k=\bar\xi_k\ldots\bar\xi_2\bar\xi_1$, $\bar s_0=1$ and
\be
\hat\tau_1^{m\bar m}(\xi,\bar\xi|g)=
s_m^{-1}
\left(\sum_{{k \geq m}\atop{\bar k \geq \bar m}}
s_k\bar s_{\bar k} <k|g|\bar k>\right) \bar s_{\bar m}^{-1}=
s_{m-1}^{-1}\left( D_{\xi_m} \bar D_{\bar \xi_{\bar m}}\tau_1(\xi,\bar\xi|g)
\right)\bar s_{\bar m-1}^{-1}
\label{clataumm1sq}
\ee
Here\footnote{There is an ambiguity in the choice of these operators
as the $\tau$-function is a linear function of times and, therefore,
any linear operator which makes unity from $\xi$ is suitable. We fix
them to act naturally on the $q$-exponential and correspond to the
difference operators that we used throughout this review.}
$ D_{\xi_i} f(\xi)\equiv
\frac{1}{\xi_i}\frac{M_i^{+2}-1}{q^2-1}f(\xi)$, $ \bar
D_{\bar\xi_i}f(\bar\xi) \equiv
\left[\frac{M_i^{-2}-1}{q^{-2}-1}f(\bar\xi)\right]\frac{1}{\bar\xi_i}$ (in
these operators, the order is also crucial!).
Now one can manifestly express $\tau_n$ through $\tau_1$, using
formulas (\ref{detq}) and (\ref{clataumm1sq}). In Appendix 3, we describe
how it works for the simplest case of the second fundamental
representation.

FBA for the $\tau$-function in parametrization A is given by
the following expressions
\be
\Psi^{+,n+k+1}_n(\xi)=q^{n+1}s_n^{-1}s_{n+k}M^+_{n+1}\ldots M^+_{n+k-1}
\left(\tau_n(\xi)-\xi_nD_n\tau_n(\xi)\right)\\
\Psi^{+,n+1}_n(\xi)=q^{n+1}
\left(\tau_n(\xi)-\xi_nD_n\tau_n(\xi)\right),\ \
\Psi^{+,n}_n(\xi)=-q^nD_n\tau_n(\xi)\\
\Psi^{-,k}_n(\xi)-\xi_{n-1}D_{n-1}\Psi^{-,k}_n(\xi)=\\=q^{n-2}
s_{k-2}^{-1}s_{n-1}D_{i-1}\tau_n(\xi)+q^{n-2}
s_{k-1}^{-1}s_{n-1}\xi_iD_i\tau_n(\xi)
\ \hbox{for}\ k>n\\
\Psi^{-,n}_n(\xi)-\xi_{n-1}D_{n-1}\Psi^{-,n}_n(\xi)=q^{n-2}
\tau_n(\xi)+ \xi_nD_n\tau_n(\xi)\\
\Psi^{-,n-1}_n(\xi)=q^{n-2}\xi_{n-1}\tau_n(\xi),\ \ \Psi^{-,k}_n(\xi)=0\
\hbox{for}\ k<n-1
\ee
Substituting these expressions in (\ref{BIA}), one obtains a set of
equations which is a quantum counterpart of the KP/Toda hierarchy in
parametrization A.

To conclude this section, point out once more that the determinant
representations of the classical $\tau$-functions do not turn into exact
$q$-determinant representations at $q\ne 1$. The reason for this is that, in
the quantum case, the evolution operator is no longer the group element.
This happens because no nilpotent subgroup $N\G_q$ exists in the quantum
group. To avoid this problem, one could begin from a slightly different
parametrization of the $\tau$-function such that the evolution operator lies
in the Borel (not just nilpotent) subgroup $B\G_q$. In the classical limit,
additional Cartan generators can be removed by redefinition of the element
$g$ (labelling the point of the Grassmannian). However, in the quantum case
the Cartan part of the evolution would essentially change the result: the
evolution operator becomes a group element (for $B\G_q$) and all "twists" of
times in formula (\ref{detq}) disappear.
Thus defined $\tau$-function is just the $q$-determinant. In order to
fulfil the whole scheme, one still needs to find an appropriate
parametrization of (a set of) group elements
$B\G_q$ by exactly $r_{\G}$ ``times".

\section{Wave function and $S$-matrix in quantum Liouville (Toda) theory}
\setcounter{equation}{0}
As it has been demonstrated above, the group acting in the space of
the spectral parameter, i.e. on solutions of integrable system, in the course
of quantization should be merely replaced by the corresponding quantum group.
In essence, it is the very sense of the quantization procedure. Now we are
going to consider an example of some different group structure that can be
also manifested in the Toda theory. It is the group acting in "the
space-time", i.e. on variables of equation. It is the remarkable property of
this group that it survives in the course of quantization, but the way of
dealing with it as well as the objects under consideration do completely
change.

Namely, consider again the Toda molecule additionally reduced to the Toda
chain, i.e. depending only on the sum of times. Let us demonstrate that this
system can be obtained via the Hamiltonian reduction of a free system given
on the cotangent bundle for a Lie group \cite{OP2,Perelomov}.

\subsection{Classical Toda system as Hamiltonian reduction}
Indeed, consider the cotangent bundle $T^*\G$ for the real split group
$G$. By the group shifts, it can be reduced to the pair
$(Y,g),~Y\in{\cal
G}^*,~g\in \G$, where ${\cal G}$ is the algebra corresponding to the
group $\G$. There is the canonical bi-invariant symplectic form on
$T^*\G$
\beql{r1}
\om=\de Y(\de gg^{-1})
\eq
and the set of invariant commuting Hamiltonians
\beql{r2}
\frac{<Y^{d_k}>}{d_k}, k=1\ldots,r
\eq
where $d_k=2,\ldots$ are invariants of ${\cal G}$ and $Y^{d_k}$ are
$ad^{\ast}\G$-invariant polynomials on ${\cal G}^*$. It is the upstairs
Hamiltonian system.

Now consider the symplectic reduction on $T^*\G$ (i.e. the reduction
preserving the symplectic form $\om$ (\ref{r1}))
w.r.t. the action
of the left and right nilpotent subgroups $\bar{N}\oplus N$ of
the group $G$
\beql{r8}
g\ms vg, ~Y\ms vYv^{-1},~~v\in
\bar{N}, \\ g\ms gn,~Y\ms Y,~~n\in N
\eq
It gives the two moment maps
\beql{r9}
\mu_v=Pr_{\bar{\cal N}^*} ~Y,~~~\mu_n=Pr_{{\cal N}^*}~
g^{-1}Yg
\eq
Making use of the Gauss decomposition, one can transform $g$ by
(\ref{r8}) to the Cartan subgroup $A$. Let $g=\exp\phi,~\phi\in{\cal A}$.
Assume that
\beql{r10}
\mu_v=Pr_{\bar{\cal N}^*} ~Y=\mu^L,~~~\mu_n=Pr_{{\cal N}^*} ~g^{-1}Yg=\mu^R
\eq
where $\mu^{L}=\sum_{\al\in\Pi}\mu^{L}_{\al}{\cal G}_{\al}$,
$\mu^{R}=\sum_{\al\in\Pi}\mu^{R}_{\al}{\cal G}_{-\al}$, $\mu^{R,L}_{\al}$ are
arbitrary constants and $\Pi$ denotes the set of all positive roots.
After ``diagonalizing" $g$ by the Gauss decomposition at the point
$g=\exp\phi,~\phi\in{\cal A}=h_{\G}$, one can manifestly solve the
constraint (\ref{r10})
\beql{r11}
Y=p+\sum_{\al\in\Pi}\left(\mu^R_{\al}e^{\al(\phi)}{\cal G}_{\al}
+\mu^L_{\al}{\cal G}_{-\al}\right)
\eq
Then, the reduced symplectic form acquires the canonical form
\be
\om^{red}=\de p\de\phi=\sum_{k=1}^r\de p_k\de\phi_k
\ee
and the reduced phase space $\bar N\setminus \setminus T^*\G//N$
is of dimension
\be
2\dim \G-2\dim \bar N -2\dim N=2r
\ee
Combination $W(Y;\tau_1\ldots,\tau_r)=\sum_{k=1}^r\frac{\tau_k}{d_k}
<Y^{d_k}>$ defines the commuting Hamiltonians of the classical hierarchy
given by the Lax operator (\ref{r11}). In particular,
\be\label{r13}
<Y^2>=\frac{1}{2}\sum_{k=1}^rp_k^2+
\sum_{\al\in\Pi}\mu_{\al}e^{2\al(\phi)}
\ee
is the conventional Toda chain Hamiltonian \cite{FT}. The classical
action in this representation takes the form
\beql{r12}
S^I=\int Y(\d_tgg^{-1})-<Y^2>+<B_L,Y-\mu^L>+<B_R,g^{-1}Yg-\mu^R>
\eq
where $B_L\in \bar{\cal N},~B_R\in {\cal N}$ are the Lagrange multipliers.

In order to compare with the standard Toda chain formulas, it suffices
to choose as $\G$ the group $SL(p)$ and consider its fundamental
representation. Then, $Y$ (\ref{r11}) coincides with the Toda chain Lax
operator etc.

The above construction is essentially based on the Gauss decomposition
\cite{GKMMMO}. In fact, the same system can be obtained via Iwasawa
decomposition. Just this latter decomposition is the standard one used so far
in studies of these system. In particular, in the next subsections we
consider quantization of the Toda system, which has been originally performed
in the Iwasawa decomposition framework and led to the so-called Whittaker
wave function (WF). However, we use everywhere the Gauss decomposition that
is more applicable to the affine algebras \cite{GKMMMO}. One can
find in \cite{GKMMMO} further details and some discussion of the
connection of the two approaches -- the Gauss and Iwasawa ones.

\subsection{Solution to the Liouville quantum mechanics -- general scheme}
A non-trivial property of the ``space-time" group manifested
above\footnote{This space-time group is not to be mixed with the group
$SL(p)$ describing the reduction to the Toda molecule (see s.3). In order to
obtain that these two are absolutely different, one can remark that the
two-dimensional Toda molecule is described by the same reduction group
$SL(p)$ but by the affine space-time group $\widehat SL(p)$.}
in the $SL(p)$ Toda molecule is that, in the course of quantization, this
group remains classical, while its interpretation changes: quantization
implies an orbit interpretation. Indeed, while the classical
system is obtained via the Hamiltonian reduction of a free system
given on the cotangent bundle for a simple real Lie group
\cite{OP2}, the quantum model is related to irreducible unitary
representations of {\it the same} group. Thus, the quantum system should be
rather interpreted within the geometrical quantization approach
\cite{Hart}.

Let us now briefly describe the general approach that allows one to solve
differential equations for the eigenvalues (EV) and, in particular, the
Schroedinger equation, in group theory terms. After this, we consider
concrete examples. Technically one needs to make some steps.

{\bf 1)} Let $g(\xi|T) \in A_{\G}(\xi)\otimes U_{\cal G}$ be
``the universal group element" (see s.4.3) of the Lie group
$\G$,\footnote{This can also be a quantum group, which leads, however, to the
difference Liouville equation -- see \cite{ORog}.} where $\xi$ somehow
parametrizes the group manifold and
$T$ are generators of $\G$ in some (not obligatory irreducible)
representation, their only property being
$[T^a, T^b] = f^{abc}T^c$.

{\bf 2)} For any given parametrization $\{\xi\}$, one can introduce
two sets of differential operators ${\cal D}_{R,L}(\xi)$ such that
\be
{\cal D}_L^a(\xi) g(\xi|T) = T^a g(\xi|T)\\
{\cal D}_R^a(\xi) g(\xi|T) = g(\xi|T) T^a
\ee
These operators satisfy the obvious commutation relations
\be
\left[{\cal D}_L^a, {\cal D}_L^b\right] = -f^{abc}{\cal D}_L^c\\
\left[{\cal D}_R^a, {\cal D}_R^b\right] = f^{abc}{\cal D}_R^c\\
\left[{\cal D}_L^a, {\cal D}_R^b\right] = 0
\label{coreDD}
\ee

{\bf 3)} For a given representation ${\cal R}$,
scalar product $<\ |\ >$
and two elements of the representation
$<\psi_L|$ and $|\psi_R>$ construct the matrix element
\be
F_{{\cal R}}(\xi | \psi_L,\psi_R) =
<\psi_L | g(\xi|T) |\psi_R>
\ee
The action of any combination of the differential operators
${\cal D}_R$ on $F$ inserts the same combination of generators $T$ to the
right of $g(\xi|T)$. If $|\psi_R>$ happens to be an eigenvector of
this combination of generators, $F$ provides a solution to the
corresponding differential equation (as the eigenvalue problem).

{\bf 4)} Of special interest are the Casimir operators, since
$|\psi_R>$ which are their eigenvectors are just elements of irreducible
representations.

{\bf 5)} In particular, the quadratic Casimir operator, when expressed
through ${\cal
D}_R$, is the Laplace operator and gives rise to some important
Schroedinger equations, of which the equation for the Toda system is a
typical example.

In fact, in order to get the Schroedinger equation for the Toda-type system,
one needs to impose additional constraints on the states
$<\psi_L|$ and $|\psi_R>$, which correspond to the Hamiltonian reduction
(\ref{r8}) of the free motion on the group manifold of
$\G$ (in the classical case). Such a reduction, as we saw in the previous
subsection, is usually associated with some decomposition of
$\G$ into a product of subgroups, and different decompositions can lead to
equivalent reductions. Essentially, reduction allows one to get rid of the
dependence of
$F(\xi)$ on some of the coordinates $\xi$. In the case of finite dimensional
Lie groups, the remaining coordinates can be choosen to correspond to the
Cartan torus, while, for affine algebras, it is more natural to preserve the
dependence on all the diagonal matrices. The matrix element we discuss here,
WF $F$, is sometimes called Whittaker function.

Let us remark that the differential operators ${\cal D}_L$
and ${\cal D}_R$ realize respectively the left and the right regular
representations of the algebra of $\G$
(in fact, the left one is antirepresentation), which can be
invariantly given by the action on the algebra $A_{\G}$
of functions on the group:
\be
\pi_{reg}^L(h)f(g)=f(hg),\ \ \ \pi_{reg}^R(h)f(g)=f(gh),\ \ \ g,h\in G
\ee
Manifestly these operators can be constructed in the following way.
Let us consider the group element $g$ and the (formal) differential
operator $\hbox{d}$ acting as the full derivative on functions of
$\xi_i$, i.e. ${\hbox{d}}\equiv
\sum_id\xi_i {\p\xi_i}$. Then, one may calculate
$g^{-1}\cdot{\hbox{d}}g$
(Maurer-Cartan form) and expand it into the generators of the algebra:
\be\label{d1hren}
g^{-1}\cdot{\hbox{d}}g=\sum_{a,i}c_{a,i}T_ad\xi_i
\ee
i.e.
\be\label{d2hren}
{\hbox{d}}g=\sum_{a,i}c_{a,i}\left({\cal D}_R^ag\right)d\xi_i
\ee
Now one reads off the manifest form of the differential operators
${\cal D}_R$ from this expression. Analogously, one can calculate ${\cal
D}_L$.

Note that it suffices to consider $g$ only in the fundamental representation,
since matrix elements in this representation generate the whole algebra
$A_{\G}$ and the group action can be extended onto the whole algebra using
the co-multiplication. At the same time, calculating the coefficients
$c_{a,i}$ in the fundamental representation is quite easy.

As an instance of the general scheme, we briefly consider below the example of
the Schroedinger equation for the Toda-type system, while the details can be
found in \cite{GKMMMO}. We begin with the simplest case of the group
$SL(2)$, i.e. the Liouville quantum mechanics. The problem we address to is
to find an integral representation for the WF and its asymptotics
that define the scattering $S$-matrix.

\subsection{Liouville quantum mechanics -- $SL(2,\R)$}
\paragraph{Liouville system for the group $SL(2,\R)$.}
The corresponding Lie algebra is given by the relations
\be
\phantom{dsd}[T_+,T_-]=T_0,\ \ [T_{\pm},T_0]=\mp 2T_{\pm}
\ee
In the fundamental representation
\be
T_0=\left(\begin{array}{cc}1&0\\0&-1\end{array}\right),\
T_+=\left(\begin{array}{cc}0&1\\0&0\end{array}\right),\
T_-=\left(\begin{array}{cc}0&0\\1&0\end{array}\right)
\ee
The quadratic Casimir operator is equal
\be
C=(T_-T_++T_+T_-)+\2 T_0^2=2T_-T_++T_0+\2 T_0^2
\ee

Let us now perform the first step of our general scheme and
parametrize the group element as follows\footnote{This parametrization
differs from that used in the previous sections by the permutation
of the nilpotent subgroups, which is given by action of the inner
group automorphism -- antipode.}
\be
g(\psi,\phi,\chi|T) = e^{\psi T_-}e^{\phi T_0}e^{\chi T_+}
\ee
In the fundamental representation, the group element has the form
\be\label{groupelG}
g = \left( \begin{array}{cc} 1 & 0 \\ \psi & 1 \end{array}\right)
 \left( \begin{array}{cc} e^{\phi} & 0 \\ 0 & e^{-\phi} \end{array}\right)
   \left( \begin{array}{cc} 1 & \chi \\ 0 & 1 \end{array}\right) =
  \left( \begin{array}{cc} e^{\phi} & \chi e^{\phi} \\
   \psi e^{\phi} & \psi\chi e^{\phi} + e^{-\phi} \end{array}\right)
\ee
As the next step, we need to construct differential operators
realizing the right and the left regular representations. We already remarked
that, to this end, one should calculate ``the currents"
$g^{-1}\cdot {\hbox{d}}g$ and ${\hbox{d}}g\cdot g^{-1}$:
\be\label{gdgG}
g^{-1}dg = \left(\begin{array}{cc}
      - e^{2\phi}\chi d\psi + d\phi &
              - e^{2\phi}\chi^2 d\psi + 2\chi d\phi + d\chi \\
      e^{2\phi}d\psi & e^{2\phi}\chi d\psi - d\phi \end{array}\right)\\
dg\cdot g^{-1} =  \left(\begin{array}{cc}
      - e^{2\phi}\psi d\chi + d\phi &  e^{2\phi}d\chi \\
 -e^{2\phi}\psi^2 d\chi + 2\psi d\phi + d\psi & e^{2\phi}\psi d\chi - d\phi
         \end{array} \right)
\ee
Using this formula and (\ref{d1hren}), one easily gets
\be\label{rrG}
\frac{\partial g}{\partial\phi} = g(T^0 + 2\chi T^+) =
                (T^0 + 2\psi T^-)g\\
\frac{\partial g}{\partial\chi} = gT^+ =
   e^{2\phi}(T^+ - \psi T^0 - \psi^2 T^-)g\\
\frac{\partial g}{\partial\psi} =
   g(-\chi^2 T^+ - \chi T^0 + T^-)e^{2\phi} = T^-g
\ee
and, using (\ref{d2hren}), --
\be\label{rrrG}
{\cal D}_R^+ = \frac{\partial}{\partial\chi},\ {\cal D}_R^0=-2\chi{\p\chi}+
{\p\phi},\ {\cal D}_R^-=e^{-2\phi}{\p\psi}+\chi{\p\phi}-\chi^2{\p\chi}
\ee
\be\label{lrrG}
{\cal D}_L^- = \frac{\partial}{\partial\psi},\ {\cal D}_L^0=-2\psi{\p\psi}
+{\p\phi},\ {\cal D}_L^+=-\psi^2{\p\psi}+\psi{\p\phi}+e^{-2\phi}{\p\chi}
\ee

At the third step, we need to choose some representation. Consider
the principal (spherical) series of representations induced by
one-dimensional representations of the Borel subalgebra
\cite{GGPS}. The space of the representation consists of functions of one
variable and matrix elements (scalar product) are given by
integrals with the flat measure. Action of the algebra is given
by the differential operators
\be\label{j-rep}
T_+ = \frac{\partial}{\partial x},\ T_0=-2x{\p x}
+2j,\ T_-=-x^2{\p x}+2jx
\ee
where $j$ is spin of the representation.
We consider here only unitary representations, since, for these
representations only, the scalar product is given by a convergent integral.
The unitarity condition implies that $j+\2$ is pure
imaginary\footnote{Similarly, for the general $SL(p)$, all
$j_i+\2$ (i.e. $\bj+\brho$) should be pure imaginary.}.

The fourth step implies that we consider the quadratic Casimir
operator (it can be calculated both in the left and in the right
regular representations):
\be
C=\2 {\partial^2\over\partial\phi^2}+{\p\phi}+2e^{-2\phi}
{\partial^2\over\partial \psi\partial\chi}
\ee
Its EV in the spin $j$ representation is
$2(j^2+j)$. Therefore, the matrix element
$F=<\psi_L|g(\theta,\phi,\chi)|\psi_R>$ with $<\psi_L|$ and $\psi_R|$
belonging to the spin $j$ representation satisfies the equation
\be
CF=2(j^2+j)F
\ee

As the last, fifth step we need to fix the reduction conditions.
We choose them to be
\be {\p\chi}F_G=i\mu_RF_G,\ {\p\psi}F_G=i\mu_LF_G
\ee
i.e. (see (\ref{rrrG}) and (\ref{lrrG}))
\be\label{statecondG}
T_+|\psi_R>=i\mu_R|\psi_R>,\ \ <\psi_L|T_-=i\mu_L<\psi_L|
\ee
Under this reduction, the Hamiltonian (quadratic Casimir operator)
is equal to (cf. (\ref{r13}))
\be\label{HG2}
H=\2 {\partial^2\over\partial\phi^2}+{\p\phi}-2\mu_R\mu_Le^{-2\phi}
\ee
and the function $\Psi(\phi)=e^{\phi}F$ satisfies the following
Schroedinger equation
\be\label{SchG}
\left[\2
{\partial^2\over\partial\phi^2}-2\mu_R\mu_Le^{-2\phi}\right]
\Psi(\phi)=2\left(j+{\f 2}\right)^2\Psi(\phi)\equiv
\lambda^2\Psi(\phi)
\ee
which is actually the Liouville Schroedinger equation.

Note that, using the representation (\ref{groupelG}) and the
conditions
(\ref{statecondG}), one can immediately obtain the equation
(\ref{HG2}) through the following chain of formulas\footnote{
From now on, we consider only the matrix element of
$e^{\phi_I T_0}$, instead of
$g$, in order to exclude trivial $\psi$- and $\chi$-dependencies of
$F$.}:
\be
2j(j+1)F_G^{(j)}\equiv 2j(j+1)<\psi_L|e^{\phi T_0}|\psi_R>_j=
<\psi_L|e^{\psi T_-}e^{\phi T_0}e^{\chi T_+}\hat C|\psi_R>_j=\\=
<\psi_L|e^{\phi T_0}(2T_-T_++T_0+\2 T_0^2)|\psi_R>_j=
\left(-2\mu_R\mu_Le^{-2\phi}+{\p\phi}+\2 {\partial^2\over\partial
\phi^2}\right)<\psi_L|e^{\phi T_0}|\psi_R>_j
\ee

\paragraph{Solving the Liouville Schroedinger equation.}
Now we should only solve the reduction conditions (\ref{statecondG}) and
find out some explicit integral representation for the WF. We use
formulas (\ref{j-rep}). Then, the conditions (\ref{statecondG})
take the form
\be
T_+|\psi_R>={\p x}\psi_R(x)=i\mu_R\psi_R(x),\\
<\psi_L|T_-=(2x+x^2{\p x}+2jx)\psi_L(x)=i\mu_L\psi_L(x)
\ee
Their solutions are
\be
\psi_R(x)=e^{i\mu_Rx},\ \ \psi_L(x)=x^{-2(j+1)}e^{-{i\mu_L\over x}}
\ee
This finally provides us with the solution to the Schroedinger equation
(\ref{SchG}):
\be\label{LWFG}
e^{\phi}F^{(j)}(\phi)=e^{\phi}<\psi_L|e^{\phi T_0}|\psi_R>_j
=e^{\phi}\int x^{-2(j+1)}e^{-{i\mu_L\over x}}e^{\phi(2j-2x{\p x})}
e^{i\mu_Rx}\ dx=\\= \left({i\over \mu_R}\right)^{-(2j+1)}e^{-(2j+1)\phi}
\int_0^{\infty} x^{-2(j+1)}e^{-{\mu_L\mu_Re^{-2\phi}
\over x}- x}dx=\\=2\left(i\sqrt{{\mu_L\over\mu_R}}\right)^{-(2j+1)}
K_{2j+1}(2\sqrt{\mu_L\mu_R}e^{-\phi})\ dx
\ee
where $K_{\nu}(z)$ is the Macdonald function
\cite{GraRyz}.

\paragraph{Harish-Chandra function and its asymptotics.}
Let us once more return to the Schroedinger equation (\ref{SchG}). The
condition of unitarity of some representation implies pure imaginary
$\lambda\equiv \sqrt{2}(j+\2)$, i.e. continuous spectrum of the
Schroedinger equation (and real energy). Solutions of this equation
obviously oscillate at only one infinity, while exponentially decreasing at
the other one. We can consider the scattering problem. Then, the reflection
$S$-matrix is equal to the ratio of the coefficients of falling and reflected
waves. We really need to find the exponential asymptotics
$e^{\pm\lambda\phi}$ at {\it the same} infinity. Technically, it can be done
as follows: one considers a small real shift of $\lambda$ so that, depending
on its sign, one can extract one of the two asymptotics.

We perform this procedure for the Macdonald function:
\be\label{Mac2}
K_\nu(z)={\f 2}\left({z\over 2}\right)^\nu\int_{0}^\infty
e^{-t-{z^2\over 4t}}t^{-\nu-1}dt
\ee
If $\hbox{Re}\nu\le 0$, one can through away the asymptotically
(i.e. at small $z$) small term ${z^2\over 4t}$ so that the integral
is equal to
\be
K_\nu(z)\stackreb{z\to 0}{\sim}{\f 2}\left({z\over
2}\right)^\nu\int_0^\infty e^{-t}t^{-\nu-1}dt=
-{\pi\over 2}{1\over \sin\pi\nu \Gamma(1+\nu)}\left({z\over 2}\right)^\nu
\ee
On the other hand, at $\hbox{Re}\nu\ge 0$ one needs to
make the replace $t\to
z^2t$ and the result reads as
\be
K_\nu(z)\stackreb{z\to 0}{\sim}{\f 2}(2z)^{-\nu}\int_0^\infty
e^{-{1\over 4t}}t^{-\nu-1}dt=
{\pi\over 2}{1\over \sin\pi\nu \Gamma(1-\nu)}\left({z\over 2}\right)^{-\nu}
\ee

Note that, though the ratio of the two calculated asymptotics is unambiguously
defined, each of them separately is defined up to a common normalization
factor of the WF. This factor is used to fix partially requiring the absence
of poles at finite values of momentum. The issue of the entirely
fixed normalization is discussed in \cite{GKMMMO}.

In this concrete case, the condition of absence of the poles implies that
the WF (\ref{LWFG}) is to be multiplied by the function
$\sin\pi\lambda$ that, cancelling all the poles, introduce new additional
zeroes. The two asymptotics of the so-normalized WF
are\footnote{Hereafter, we omit the trivial factor depending on the
cosmological constant $\mu_L\mu_R$.}
\be\label{hchsl2}
c_{\pm}={\f \Gamma (1\pm \lambda)}
\ee
These functions are called Harish-Chandra functions \cite{HC} and play an
important role in the group theory.

\subsection{Toda quantum mechanics}
\paragraph{Quantum mechanics for the Toda molecule -- $SL(p,\R)$.}
We apply now the above developed general scheme to more general case of the
group $SL(p,R)$. Calculations in this case are quite similar to the
Liouville case, however, they require some additional information
on the structure of the unitary representations of
$SL(p,R)$. Therefore, we write down here only the results and refer to
Appendix 4 and, especially, the reference \cite{GKMMMO} which contain the
detailed calculations.

Namely, like the previous subsection, we consider the matrix element
$F$ given in some unitary representation of the group and act on it
by some Casimir operator. Since the group $SL(p,R)$ is of rank $p-1$,
i.e. has $p-1$ independent Casimir operators, the matrix element $F$
satisfies now $p-1$ differential equations for the EV simultaneously.
Thus, this matrix element solves the quantum problem with the Hamiltonians
of s.6.1 corresponding to the Toda system. Among these Hamiltonians,
there is the quadratic Hamiltonian (\ref{r13}) that gives rises to
the second order equation, that is to say Schroedinger equation. It is
of the form
\be
\left( {\partial^2\over\partial\bphi^2}
-2\sum_i\mu^L_i\mu_i^Re^{\balpha_i\bphi}
\right)\Psi^{(\blambda)}(\bphi)=\blambda^2
\Psi^{(\blambda)}(\bphi)
\ee
the function $\Psi^{(\blambda)}(\bphi)$ being
the rescaled matrix element
$e^{-\brho\bphi} F^{(\blambda)}(\bphi)$ and
$\mu^{R,L}_i$ being a set of the cosmological constants.
Manifestly calculating this matrix element, one gets the integral solution
to this equation:
\be\label{res}
\Psi (\phi)=e^{-\blambda\bphi}
\int \prod_{i<j}dx_{ij}\prod_{i=1}^{p-1}
\Delta_i^{-(\blambda\balpha_i+1)}(xS^{-1})\times
e^{ i\mu_{i}^Rx_{i,i+1}e^{\balpha_i\bphi} -i\mu_{p-i}^L
\frac{\Delta _{i,i+1}(xS^{-1})}{\Delta _{i}(xS^{-1})}}
\ee
Here $\balpha_i$ are the positive simple roots, $S$ is the inner automorphism
of the group $SL(p)$, which maps the upper-triangle matrices to
the lower-triangle ones,
$S_{ij}\equiv\delta_{i+j,p+1}$, $\Delta_i(A)$ is the upper main $i\times
i$ minor of the matrix $A$ and $\Delta_{i,i+1}(A)$ is the same minor of the
matrix with $i$ and $i+1$ columns exchanged.

Now looking at the asymptotics of this WF, one can obtain the
Harish-Chandra functions \cite{GKMMMO}. First of all, note that
the number of different asymptotics coincides with the number of elements
of the Weyl group for $SL(p)$, while the WF (\ref{res}) is expressed
through only the simple roots. Thus, it is natural that the Harish-Chandra
functions are labeled by elements of the Weyl group and are related
by action of this group. In order to calculate these asymptotics, one
may use the standard technique developed for the Iwasawa decomposition
case \cite{GinKor} and the result reads as
\be\label{gfh}
c_s(\blambda)=\prod_{\balpha\in\De^+}{\f \Gamma(1+s\blambda\cdot\balpha)}
\ee
where $s$ denotes elements of the Weyl group and the product runs
over all the positive roots. The ratio of these Harish-Chandra
functions gives us the reflection $S$-matrix of the theory.

\paragraph{Quantum field Liouville system.}
As the next step, we extend our construction to the affine case.
Indeed, considering the group $\widehat{SL(2)}$ and the ``point-wise"
Gauss decomposition (used for bosonization of affine algebras -- see
\cite{GMMOS}), one obtains the Schroedinger equation for the
two-dimensional {\it field} Liouville system. This means that,
calculating the Harish-Chandra functions, we get in this case
the $S$-matrix (or the two-point correlation function) for the
two-dimensional Liouville field theory. Indeed, choose the system of positive
roots in this case as follows
\be
\balpha_0 + n(\balpha_0+\balpha_1),\
n(\balpha_0+\balpha_1),\
\balpha_1 + n(\balpha_0+\balpha_1),\
n = 0,1,2,\ldots\
\ee
Denote
\be
\blambda\cdot\balpha_0 = \frac{1}{2} - p + \tau,\ \ \
\blambda\cdot\balpha_1 = -\frac{1}{2} +p
\ee
Then, one can use formula (\ref{gfh}) for the Harish-Chandra functions
only shifting the argument of the $\Gamma$-function by 1/2 because of the
affine situation:
\be\label{HCHA}
c(\blambda) = \prod_{n\geq 0}
\Gamma^{-1}(p + n\tau)
\prod_{n\geq 1}\Gamma^{-1}(n\tau)\Gamma^{-1}(1-p + n\tau)
\ee
This expression certainly requires a careful regularization, but
all the infinite products cancel from the corresponding reflection
$S$-matrix (two-point function)
\be
S(p) = \frac{c(-p)}{c(p)} =
\frac{\Gamma(1+p)\Gamma(1+{p\over\tau})}{\Gamma(1-p)\Gamma(1-{p\over\tau})}
\ee
This expression coincides with the two-point functions obtained in the papers
\cite{DO,ZZ} in a very different way\footnote{In notations of
\cite{ZZ}, $p=2iP/b$ and $\tau=b^2$.}.

Now we introduce a fundamental building block -- the function
\be
i(\blambda)=i(p,\tau)\sim\prod_{m,n\ge 0}(p+m+n\tau)
\ee
where the product goes over only one (positive) quadrant in the plane
$m,n$. One can construct from this block both the Harish-Chandra function
(\ref{HCHA}) and the elliptic theta-functions
\be
\theta(p+\2+{\tau\over
2},\tau)\sim i(p,\tau)i(-p,\tau)i(p,-\tau)i(-p,-\tau)
\ee
and also the $q$-exponential
\be
e_q(e^{2\pi ip})\sim {\f i(p,\tau)i(-p,-\tau)},\ \ \ 1=e^{i\pi\tau}
\ee
just by taking products over different quadrants. All these expressions
certainly require some accurate regularization, since the unrestricted
products diverge. The regularized expressions can be found in
\cite{Mam} and \cite{ZZ}.

\app{Forced hierarchies}
In this Appendix, we reproduce some technical calculations
omitted from s.2.3 \cite{KMMOZ,KMMM}. Namely, we consider the forced
hierarchy given by the condition
\beq
\tau _n = 0\hbox{  , }  n < 0
\eeq
This corresponds to choosing the element of the Grassmannian in the form
\beq
G = G_0P_+
\eeq
where $P_+$ is the projector onto the positive states:
\beq
P_+|n\rangle  = \theta (n)|n\rangle
\eeq
which is realized as the fermionic operator
\beq\label{B4}
P_+ = :\exp \left[\sum _{i<0}\psi _i\psi ^\ast _i\right]:
\eeq
and enjoys the properties
\beq\label{B5}
P_+\psi ^\ast _{-k} = \psi _{-k}P_+ = 0\hbox{  , }  k > 0\
\eeq
\beq\label{B6}
\left[ P_+,\psi _k\right] = [P_+,\psi ^{\ast} _k] = 0\hbox{  , } k \geq  0\
\eeq
\beq\label{B7}
P^2_+ = P_+
\eeq
We also fix $G_0$ to include only positive fermionic modes
$\psi _k$ and $\psi ^\ast _k$ at $k \geq 0$:
\beq\label{B8}
G_0 = :\exp \left\{ \left(\int _\gamma  A(z,w)\psi _+(z)
\psi ^\ast _+(w^{-1})dzdw \right) - \sum _{i\geq 0} \psi _i\psi ^\ast _i
\right\}:\ \
\eeq
where $\psi _+(z) =\sum _{k\geq 0}\psi _kz^k$ , $\psi ^\ast _+(z)
=\sum _{k\geq 0}\psi ^\ast _kz^{-k}$ and $\gamma $ is some integration
domain. Besides, we need the projector onto the negative states
\beq\label{B9}
P_- = :\exp \left[- \sum _{i\geq 0}\psi _i\psi ^\ast _i\right]:\
\eeq
enjoying the following properties
\beq\label{B10}
P_-\psi _k = \psi _k^\ast P_- = 0\hbox{  , }  k \geq 0\
\eeq
\beq\label{B11}
\left[ P_-,\psi _{-k}\right] = [P_-,\psi _{-k}{^\ast} ] = 0\hbox{  , }
k > 0\
\eeq
\beq\label{B12}
P^2_- = P_-
\eeq
Let us calculate the state
\beq
G_0P_+e^{\bar H(y)}|n\rangle \
\eeq
One can easily check that it vanishes at $n<0$. Indeed, using
(\ref{vac}), (\ref{A19}) and (\ref{B6}), one can obtain that, at $n < 0$,
the state
$$
e^{\bar H(\bar t)}|n\rangle  = \psi ^\ast _{-n}(\bar t)\hbox{ ... }
\psi ^\ast _{-1}(\bar t)e^{\bar H(\bar t)}|0\rangle
$$
contains only negative modes $\psi ^\ast _{-m}$ ( $m > 0$).
Hence, action of $P_+$ annihilates this state due to formula
(\ref{B5}). At $n \geq  0$, using (\ref{vac}), (\ref{A19}) and
(\ref{B6}), one gets
\beq
P_+e^{\bar H(\bar t)}|n\rangle  = \psi _{n-1}(\bar t)\hbox{ ... }
\psi _0(\bar t)P_+e^{\bar H(\bar t)}|0\rangle
\eeq
One should use that
\beq\label{B15}
P_+e^{\bar H(\bar t)}|0\rangle  = |0\rangle
\eeq

\paragraph{Proof of (\ref{B15}).} Denote
$$
|\bar t\rangle  = P_+e^{\bar H(\bar t)}|0\rangle
$$
Then
$$
{\partial \over \partial \bar t_k} |\bar t\rangle  = P_+e^{\bar H(\bar t)}
\sum ^{k-1}_{i=0} \psi ^{\ast} _{i-k}\psi _i|0\rangle  = 0
$$
due to (\ref{vac0}) and (\ref{B5}). Since
$\displaystyle{|\bar t\rangle \left| _{\bar t_k=0} \right.= |0\rangle}$,
formula (\ref{B15}) is proved.

Thus, using formulas (\ref{B8}), (\ref{B9}) one obtains
\beq\label{B16}
\new
\begin{array}{c}
G_0P_+e^{\bar H(\bar t)}|n\rangle  = G\ \psi (\bar t)\ldots
\psi (\bar t)|0\rangle  = \\
= \sum     {1\over m!} \int _\gamma  \prod ^m_{i=1}A(z_i,w_i)dz_idw_i
\psi _+(z_1)\ldots\psi _+(z_m) \times \\
\times P_-
\psi ^\ast _+(w^{-1}_m)\ldots\psi ^\ast _+(w^{-1}_1)\psi _{n-1}(\bar t)
\ldots \psi _0(\bar t)|0\rangle
\end{array}
\eeq
Let us now check that only the $m=n$ term contributes into the
infinite sum (\ref{B16}). Indeed, at
$m > n$ the state $\psi ^\ast _+(w^{-1}_m)\ldots\psi
^\ast _+(w^{-1}_1)\psi _{n-1}(\bar t) \ldots \psi _0(\bar t)|0\rangle $
vanishes, since, in this case, some of the positive modes in
$\psi^\ast _+(w^{-1}_i)$ can be pushed to the vacuum
$|0\rangle $ cancelling it. Inversely, at $m<n$ some of the positive
modes of $\psi_k(-\bar y)$ can be pushed to the projector
$P_-$ and, due to (\ref{B10}), cancel it. Therefore,
\beq\label{B17}
\new
\begin{array}{c}
G_0P_+e^{\bar H(\bar t)}|n\rangle
= {1\over n!} \int _\gamma  \prod ^n_{i=1}A(z_i,w_i)dz_idw_i
\psi _+(z_1)\ldots \psi _+(z_n) \times \\
\times P_-
\psi ^\ast _+(w^{-1}_n)\ldots \psi ^\ast _+(w^{-1}_1)\psi _{n-1}(\bar t)
\ldots  \psi _0(\bar t)|0\rangle
\end{array}
\eeq
Now we use the following statement:
\beq\label{B18}
\new
\begin{array}{c}
\psi ^\ast _+(w^{-1}_n)\ldots \psi ^\ast _+(w^{-1}_1)\psi _{n-1}(\bar t)
\ldots \psi _0(\bar t)|0\rangle
= \Delta (w)\exp
\left[\sum ^n_{j=1}\xi (\bar t,w_j)\right]|0\rangle
\end{array}
\eeq

\paragraph{Proof of formula (\ref{B18}).}
Since the number of creation operators (w.r.t. the vacuum
$|0\rangle $) $\psi _i(-\bar y)$ is equal to the number of
annihilation operators $\psi ^\ast _+(w^{-1}_j)$,
after the normal re-ordering
$$
\psi ^\ast _+(w^{-1}_n)...\psi ^\ast _+(w^{-1}_1)\psi _{n-1}(\bar t)\ldots
\psi _0(\bar t)|0\rangle  = const\cdot |0\rangle
$$
and
$$const =
\langle 0|\psi ^\ast _+(w^{-1}_n)\ldots \psi ^\ast _+
(w^{-1}_1)\psi _{n-1}(\bar t)
\ldots \psi _0(\bar t)|0\rangle  =
$$
$$
= \left.\det
\left[\langle 0|\psi ^\ast _+(w^{-1}_i)\psi _{j-1}(\bar t)|0\rangle \right]
\right| _{i,j = 1,\ldots ,n}
$$
Now, using formulas (\ref{A20}), (\ref{A16}), one can get
$$
\langle 0|\psi ^\ast _+(w^{-1}_i)\psi _{j-1}(\bar t)|0\rangle  = w^{j-1}_i
e^{\xi (\bar t,w_i)}
$$
Thus, finally
$$const = \det [w^{j-1}_i e^{\xi (\bar t,w_i)}] =
\Delta (w)\exp [\sum ^n_{j=1}\xi (\bar t,w_j)]
$$
Substituting now (\ref{B18}) into (\ref{B17}) and using the obvious fact
that $P_-|0\rangle  = |0\rangle $
and $\psi _-(z_i)|0\rangle  = 0$, one obtains
\beq\label{B19}
\new
\begin{array}{c}
G_0P_+e^{\bar H(\bar t)}|n\rangle
= {1\over n!} \int _\gamma  \prod ^n_{i=1}A(z_i,w_i)e^{\xi (\bar t,w_i)}
dz_idw_i \Delta (w)\psi (z_1)\ldots \psi (z_n)|0\rangle
\end{array}
\eeq
At last, using one of the main formulas of the paper \cite{DJKM}
(which can be easily proved within the bosonization framework)
\beq\label{B20}
\psi (z_1)\ldots \psi (z_n)|0\rangle  = \Delta (z)
\exp \left[\bar H\left(\sum ^n_{i=1}\epsilon (z_i)\right)\right]|0\rangle
\eeq
where $\epsilon (z_i)$ is the vector with components
$\displaystyle{\epsilon _k(z_i) =
{1\over k} z^k_i}$, we obtain the desired result:
\beq\label{B21}
\new
\begin{array}{c}
G_0P_+e^{\bar H(\bar t)}|n\rangle
= {1\over n!} \int _\gamma  \prod ^n_{i=1}A(z_i,w_i)e^{\xi (\bar t,w_i)}
dz_idw_i
\Delta (w)\Delta (z) \exp
\left[\bar H\left(\sum ^n_{i=1}\epsilon (z_i)\right)\right]|0\rangle
\end{array}
\eeq
Thus, we finally get:
\beq\label{B22}
\new
\begin{array}{c}
\tau _n(t,\bar t) = \langle n|e^{H(x)}G_0P_+e^{\bar H(\bar t)}|n\rangle
= {1\over n!} \int _\gamma
\Delta (w)\Delta (z)\prod ^n_{i=1}A(z_i,w_i)e^{\xi (t,z_i)+\xi
(\bar t,w_i)}dz_idw_i
\end{array}
\eeq

\paragraph{Determinant representation.}

Using formulas (\ref{vac}), (\ref{B7}) and the property $[P_+,g_0] = 0$,
one obtains:
$$
\tau _n(t,\bar t) =
\langle 0|\psi ^\ast _0\ldots
\psi ^\ast _{n-1}e^{H(t)}G_0P_+e^{\bar H(\bar t)}\psi _{n-1
}\ldots \psi _0|0\rangle  =
$$
$$
=
\langle 0|e^{H(t)}\psi ^\ast _0(-t)\ldots
\psi ^\ast _{n-1}(-t)P_+G_0P_+\psi _{n-1}
(\bar t)\ldots \psi _0(\bar t)e^{\bar H(\bar t)}|0\rangle
$$
Since $\psi ^\ast _i(-t)$ and $\psi _i(\bar t)$ contains
only positive modes (see (\ref{A18}) and (\ref{A19})),
one gets from (\ref{B6}) and (\ref{B15})
\beq\label{B25}
\new
\begin{array}{c}
\tau _n(t,\bar t) =
\langle 0|\psi ^\ast _0(-t)\ldots
\psi ^\ast _{n-1}(-t)G_0\psi _{n-1}(\bar t)\ldots
\psi _0(\bar t)|0\rangle  = \\ \left.
\det \left[\langle 0|\psi ^\ast _i(-t)G_0\psi _j(\bar t)|0\rangle \right]
\right| _{i,j=0,...,n-1}
\end{array}
\eeq
The same arguments we used to derive (\ref{B17}) from
(\ref{B16}), when applied to
(\ref{B25}), allows one to conclude that only the linear term in
$A(z,w)$ contributes. Therefore, using (\ref{B10}), one obtains
\beq\label{B26}
\new
\begin{array}{c}
\langle 0|\psi ^\ast _i(-t)G\psi _j(\bar t)|0\rangle  = \int _\gamma
A(z,w)dzdw
\langle 0|\psi ^\ast _i(-t)
\psi _+(z)P_-\psi ^\ast _+(w^{-1})\psi _j(\bar t)|0
\rangle  = \\
= \int _\gamma  z^iw^jA(z,w)e^{\xi (t,z)+\xi (\bar t,w)}dzdw  =
\partial ^i_t \partial _{\bar t} ^j\int _\gamma
A(z,w)e^{\xi (t,z)+\xi (\bar t,w)}dzdw
\end{array}
\eeq
Thus, the final expression for the $\tau$-function in determinant
form is
\beq\label{B27}
\tau _n(t,\bar t) =
\left.\det \left[\partial ^i_{x_1}\partial _{\bar t_1}^j\int _\gamma
A(z,w)e^{\xi (t,z)+\xi (\bar t,w)}dzdw\right] \right|
_{i,j=0,...,n-1}
\eeq

\app{Fundamental representations of $SL_q(p)$}
In this Appendix we generalize the structure of fundamental representations
described in s.5.1 to the quantum case, i.e. to the group
$SL_q(p)$ with $q\neq 1$ (see \cite{GKLMM} and references therein).

The main new notion that we need is
$q$-antisymmetrization that is defined to be the sum over all
the permutations
\be
\left( [1,\ldots,k]_q\right) = \sum_P (-q)^{{\hbox{deg}}\ P}
\left(P(1),\ldots,P(k)\right)
\ee
where
\be
{\hbox{deg}}\ P = {\hbox{number\ of\ inversions\ in}}\ P
\nn
\ee
For instance, $q$-antisymmetrization is used to define the
$q$-determinant:
\be
{\det}_q A \sim  A^{[1}_{[1}\ldots A^{p]_q}_{p]_q}
= \sum_{P,P'} (-q)^{{\hbox{deg}}\ P + {\hbox{deg}}\ P'}
\prod_{a} A^{P(a)}_{P'(a)}
\label{fqdet}
\ee
Note that this quantity does not obligatory coincide with
$A^1_{[1}\ldots A^p_{p]_q}$. For instance, for $p=2$ (\ref{fqdet})
gives $\displaystyle{\frac{1}{[2]}(A^1_1A^2_2 - qA^1_2A^2_1 - qA^2_1A^1_2 +
q^2A^2_2A^1_1)}$, while $q$-antisymmetrizing over only the lower indices
would give merely $A^1_1A^2_2 -
qA^1_2A^2_1$. Moreover, $A^1_{[1}A^1_{2]} = A^1_1 A^2_2 - q A^1_2
A^1_1$ is not obliged to vanish and even
$A^1_{[1}A^1_{1]_q} = (1-q)(A^1_1)^2 \neq 0$.

The ``normal" properties of the $q$-antisymmetrization restore, if
one considers as $A^i_j$ elements of the coordinate ring $A(GL(p))$
(algebra of functions on the
quantum group $GL_q(p)$). These matrix elements satisfy the following
commutation relations (in the particular case of $p=2$, they reduce to
(\ref{core})) \cite{FRT}:
\be
\forall i, \forall j_1<j_2 \ \ \ \
A^i_{[j_1}A^i_{j_2]_q} = 0, \ \ \ \ (ab=qba,\ cd=qdc) \nn \\
\forall i_1<i_2, \forall j \ \ \ \
A^{[i_1}_jA^{i_2]_q}_j = 0, \ \ \ \ (ac = qca,\ bd = qdb) \nn \\
\forall i_1\neq i_2, j_1\neq j_2 \ \ \ \
A^{i_1}_{j_2}A^{i_2}_{j_1} = A^{i_2}_{j_1}A^{i_1}_{j_2},
 \ \ \ \ (bc = cb) \\
\forall i_1<i_2,\ j_1<j_2 \ \ \ \\
A^{i_1}_{j_1} A^{i_2}_{j_2} - A^{i_2}_{j_2}A^{i_1}_{j_1} =
(q - q^{-1}) A^{i_1}_{j_2}A^{i_2}_{j_1} \ \ \ \
(ad-da = (q-q^{-1})bc)
\ee
For $A \in GL_q(p)$
\be
{\det}_q A = A^1_{[1}\ldots A^p_{p]_q} = A^{[1}_1\ldots A^{p]_q}_p
\ee
and $A^{i_1}_{[1}\ldots A^{i_k}_{k]_q} = 0$ if
two arbitrary upper indices coincide (but all the lower indices are
different: even for $A \in GL_q(p)$ it is still correct that
$A^{1}_{[1}A^{1}_{1]_q} = (1-q)(A^1_1)^2 \neq 0$).

The notion of $q$-antisymmetrization is of importance when
constructing fundamental representations, since the $k$-th
fundamental representation
$F^{(k)}$ of $SL_q(p)$ is the $q$-skew degree of
$F = F^{(1)}$:
at $q\neq 1$ we obtain instead of (\ref{frep}):
\be
F^{(k)} = \left\{\Psi^{(k)}_{i_1\ldots i_k} \sim
\psi_{[i_1}\ldots \psi_{i_k]_q}, \ \ \ i_1<\ldots <i_k\right\}
\label{fqrep}
\ee
Note that now one should {\it manifestly} require for all
$i_a$ to differ.

All the formulas like (\ref{inttwfrep}) for the \iws remains
unchanged provided the antisymmetrization in them is replaced by
the $q$-antisymmetrization (with $q$-$\epsilon$-symbol defined in the
obvious way). Instead of (\ref{gkdet}), one now gets:
\be\label{gkdetq}
g^{(k)}\left({i_1\ldots i_k}\atop{j_1\ldots j_k}\right) \sim
{\det}_q g^{i_a}_{j_b}, \ \ \ \ i_1 < \ldots i_k,\ \hbox{or}\
j_1<\ldots j_k
\ee
while (\ref{gkgk}) turns into
\be
g^{(k)}\left({i_1\ldots i_k}\atop{[j_1\ldots j_k}\right)
g^{(k')}\left({i'_1\ldots i'_{k'}}\atop{j_{k+1}]_q j'_1\ldots j'_{k'-1}
}\right) = \\=
g^{(k+1)}\left({i_1\ldots i_k[i'_k}\atop{j_1\ldots j_{k+1}}\right)
g^{(k'-1)}\left({i'_1\ldots i'_{k'-1}]_q}\atop{j'_1\ldots j'_{k'-1}}
\right)
\label{gqgq}
\ee
\be
i_1 <\dots< i_k, \ \ \hbox{or}\ \ j_1<\ldots j_k,\ \ \ \ \
i_1'<\ldots<i_{k'}'\ \ \hbox{or}\ \ j_{k+1}<j_1'<\ldots<j_{k'-1}',\\
i_1<\ldots<i_k<i'_{k'},\ \hbox{or}\ j_1<\ldots<j_{k+1},\ \ \
i_1'<\ldots<i'_{k'-1},\ \hbox{or}\ j'_1<\ldots<j'_{k'-1}
\ee
Exactly similar to (\ref{gkgk}), these are just an identity for the
matrices from $GL_q(p)$.\footnote{It is important that $g$ belongs
to some $GL_q(p)$, i,e, its elements satisfy the proper
commutation relations. This is certainly {\it implied} in
derivation of (\ref{gqgq}), since is assumed that $g$ is the group
element. Be it not the case, one would need to understand
(\ref{gkdetq}) in the sense of (\ref{fqdet})
(in particular, one would need to write {\it and} instead of {\it or}
in (\ref{gkdetq})), and then one would run in a
contradiction with (\ref{gqgq}). To make it more transparent,
let us take $k = k' = 1$ so that (\ref{gqgq}) becomes
\be
g^{i}_{[j_1}g^{i'}_{j_2]_q} = g^{(2)}\left({i\ i'}\atop{j_1j_2}\right)
\nn
\ee
and the l.h.s. of this formula is equal to $g^i_{j_1}g^{i'}_{j_2} - q
g^{i}_{j_2}g^{i'}_{j_1}$, while the r.h.s. would be interpreted as
$g^{i}_{j_1}g^{i'}_{j_2} - q g^{i}_{j_2}g^{i'}_{j_1} -
q g^{i'}_{j_1}g^i_{j_2} + q^2g^{i'}_{j_2}g^i_{j_1}$.}
The feature that is essentially new as compared with
(\ref{gkgk}) is explicit restricting the indices
$i,j,i',j'$, which makes the translation to the language of generating
functions more sophisticated.

Let us note that, similar to the classical case, one can construct from the
quantum minors
(\ref{gkdetq}) local coordinates on the quantum flag space
\cite{Noumi}. Similar to the classical case, these coordinates satisfy a set
of (quantum) bilinear (Plucker) relations and present some quantum
generalization of the Plucker coordinates
\cite{Plucker}. Unfortunately, so far no consistent ways to parametrize the
Plucker coordinates is known.

Now we would like to discuss formula
(\ref{fqrep}) in detail. We start with the classical case, since the
structure of representations remains practically unchanged in the course of
quantization.

The Lie algebra $SL(p)$ is generated by the generators
$T_{\pm{\bf\alpha}}$ and Cartan operators $H_{\bf\beta}$ so that
$[H_\beta, T_{\pm{\bf\alpha}}] = \pm\2({\balpha\bbeta})
T_{\pm{\bf\alpha}}$. All the elements of representation are the
eigenfunctions of the Cartan generators
$H_{\bf\beta}$,
$H_{\bf\beta}|\blambda\rangle = \2({\bbeta\blambda})|\blambda\rangle$.
The highest weight vectors of the representations
$F^{(k)}$ -- ${\bmu}_k$ -- form a system dual to the {\it simple} roots
${\balpha}_i$, $i = 1,\ldots,r$:
$({\bmu}_i{\balpha}_j) = \delta_{ij}$, and
${\brho} \equiv \frac{1}{2}\sum_{{\bf\alpha}>0}{\balpha} =
\sum_i {\bmu}_i$.

The representation $F^{(1)}$ contains the states of the form
\be
\psi_i = T_{-(i-1)}\ldots T_{-2}T_{-1}\psi_1, \ \ \ i = 1,\ldots,p
\label{psii}
\ee
Moreover,
\be
T_{-i}\psi_j = \delta_{ij}\psi_{i+1}
\ee
(thus, for $T_- = \sum_{i=1}^r T_{-i}\ $  $T_-^i\psi_j = \psi_{j+i}$
and one obtains (\ref{simplfrep})), and
\be
\blambda(\psi_i) = {\bmu}_1 - {\balpha_1} - \ldots -
{\balpha}_{i_1}
\ee
Here $T_{\pm i} \equiv T_{\pm{\bf\alpha}_i}$ are the generators associated
with the simple roots. Denote the corresponding basis in the Cartan
algebra $H_i = H_{{\bf\alpha}_i}$,
${H_i}|{\blambda}\rangle = \2({{\balpha}_i{\blambda}}) |{\blambda}\rangle =
\lambda_i|{\blambda}\rangle$. Then,
\be\label{lambda}
\lambda_i^{(j)} \equiv \lambda_i(\psi_j) = {\f 2}(\delta_{ij} -
\delta_{i,j-1})
\ee
This formula, along with (\ref{psii}) and commutation relations of the
algebra, gives $||\psi_i||^2 = 1$. Since the co-multiplication in the
classical case is just
$\Delta(T) = T\otimes I + I\otimes T$, $\psi_{[1}\ldots \psi_{k]}$
exhaust the highest weight vectors (i.e. they are cancelled by all
$\Delta_k(T_{+i})$ and, thus, by all $\Delta_k(T_{+{\bf\alpha}})$).

Quantum UEA $U_q(SL(p))$ is given by the Chevalle generators
$T_{\pm i}$ corresponding to the simple roots and
$q^{\pm H_i}$, with the defining commutation relations
($A_{ij}$ is the Cartan matrix for $SL(p)$)
\cite{Drin,Ros}
\be\label{notations}
q^{H_i} T_{\pm j} q^{-H_i} = q^{\pm A_{ij}}T_{\pm j}\\
\phantom. [T_{+i},T_{-j}] = \delta_{ij}\frac{q^{2H_i}-q^{-2H_i}}{q - q^{-1}}
\ee
and the co-multiplication
\be\label{horoshij}
\Delta(T_{\pm i}) = q^{H_i}\otimes T_{\pm i} + T_{\pm i}\otimes
q^{-H_i}\\
\Delta(q^{\pm H_i}) = q^{\pm H_i}\otimes q^{\pm H_i}
\ee
Besides, the Chevalle generators satisfy additional (Serre) relations
\cite{Drin}. One can use, instead of the Chevalle basis, the basis
of generators corresponding to all roots. In this case, the Serre
identities are no longer needed (they become relations in the algebra).
However, in this case the co-multiplication for $T_{\pm\alpha}$,
which can be easily read off from (\ref{horoshij}) turns out to be
quite complicated. For instance, one obtains for the generator
corresponding to the root ${\balpha}$ of ``weight" 2,
i.e. $T_{\pm{\bf\alpha}}
= \pm [ T_{\pm {\bf\alpha}_i}, T_{\pm{\bf\alpha}_{i+1}}]$,
\be
\Delta(T_{-{\bf\alpha}}) = - [\Delta(T_{\alpha_i}),\Delta(T_{\alpha_{i+1}})]
= q^{H_\alpha}\otimes T_\alpha + T_\alpha \otimes q^{-H_\alpha} +\\+
(q^{1/2}-q^{-1/2})\left[(T_{-\alpha_i}\otimes T_{-\alpha_{i+1}})
(q^{H_{i+1}}\otimes q^{-H_i}) -
(T_{-\alpha_{i+1}}\otimes T_{-\alpha_i})(q^{H_i}\otimes q^{-H_{i+1}})\right]
\ee
With the co-multiplication given, one can easily check
formula (\ref{fqrep}). For instance, for $F^{(2)}$:
\be
\Delta(T_{+i}) (\psi_1\psi_2 - q\psi_2\psi_1) =
\delta_{i,1}(q^{\lambda_1^{(1)}}\psi_1\psi_1 - q^{1-\lambda_1^{(1)}}
\psi_1\psi_1) = 0
\ee
since $\lambda_1^{(1)} = \frac{1}{2}$. Thus, $\Psi^{(2)}_{12}
\equiv \psi_{[1}\psi_{2]_q}$ is actually the highest weight vector.
Similarly, ($i<j$):
\be
\Delta(T_{-l})\Psi^{(2)}_{ij} = \Delta(T_{-l})(\psi_i\psi_j - q\psi_j
\psi_i) = \nn \\
= \delta_{li} q^{-\lambda_i^{(j)}}(\psi_{i+1}\psi_j -
q^{1+2\lambda_i^{(j)}}\psi_j\psi_{i+1}) +
  \delta_{lj} q^{\lambda_j^{(i)}}(\psi_i\psi_{j+1} -
q^{1-2\lambda_j^{(i)}}\psi_{j+1}\psi_i)
\ee
In accordance with (\ref{lambda}), for $i<j$
always $\lambda_j^{(i)} = 0$, while
$\lambda_i^{(j)}= 0$ unless $j = i+1$, and
$\lambda_i^{(i+1)} = -{\f 2}$. Thus,
\be
\Delta(T_{-l})\Psi_{ij}^{(2)} = \delta_{lj}
\Psi^{(2)}_{i,j+1} + \delta_{li}\Psi^{(2)}_{i+1,j}(1 - \delta_{i+1,j})
\ee
This allows one to determine the action of any $\Delta(T_{-\alpha})$.
One can easily describe the action of all
$\Delta(T_{+i})$ and repeat the procedure for other representations
$F^{(k)}$.

\app{Manifest example: $\tau_2$ for $SL_q(p)$}
In this Appendix, we consider the simplest $\tau$-function for the group
$SL_q(p)$ in the simplest non-trivial (second fundamental)
representation \cite{KMM2}.
Denote through $\{u\}$ and $\{v\}$ the subsets $\{s\}$ such that
$i(s)=1$ and $i(s)=2$ respectively. Then,
\be\label{tau2}
\tau_2  = \tau_1^{00}(\{q\xi_u\},\{q^{-1}\xi_v\},\xi_s;\
\{\bar \xi_u\},\{\bar \xi_v\},\bar \xi_s) \tau_1^{11}(\{\xi_u\},
\{\xi_v\},\xi_s;\
  \{q\bar \xi_u\},\{\bar \xi_v\},\bar \xi_s) - \\
- q\tau_1^{01}(\{q\xi_u\},\{q^{-1}\xi_v\},\xi_s;\
  \{\bar \xi_u\},\{\bar \xi_v\},\bar \xi_s)
\tau_1^{10}(\{\xi_u\},\{\xi_v\},\xi_s;\
  \{q^{-1}\bar \xi_u\},\{q\bar \xi_v\},\bar \xi_s) - \\
- q\tau_1^{10}(\{q^{-1}\xi_u\},\{\xi_v\},\xi_s;\
  \{\bar \xi_u\},\{\bar \xi_v\},\bar \xi_s)
\tau_1^{01}(\{\xi_u\},\{\xi_v\},\xi_s;\
  \{q\bar \xi_u\},\{\bar \xi_v\},\bar \xi_s) + \\
+ q^2\tau_1^{11}(\{q^{-1}\xi_u\},\{\xi_v\},\xi_s;\
  \{\bar \xi_u\},\{\bar \xi_v\},\bar \xi_s)
\tau_1^{00}(\{\xi_u\},\{\xi_v\},\xi_s;\
  \{q^{-1}\bar \xi_u\},\{q\bar \xi_v\},\bar \xi_s)
\ee
Here $\xi_s$ denotes all times with $i(s)>2$. Consider now
parametrization A. In formula (\ref{tau2}), each set of $\{u\}$
and $\{v\}$ consists of a single element:
$\{u\}=\{s=1\}$, $\{v\}=\{s=2\}$. Then,
\be
  \tau_2  = \tau_1^{00}(q\xi_1,q^{-1}\xi_2,\xi_i;\ \bar
\xi_1,\bar \xi_2,\bar \xi_i) \tau_1^{11}(\xi_1,\xi_2,\xi_i;\ q\bar \xi_1,\bar
  \xi_2,\bar \xi_i) - \\ - q\tau_1^{01}(q\xi_1,q^{-1}\xi_2,\xi_i;\ \bar
  \xi_1,\bar \xi_2,\bar \xi_i) \tau_1^{10}(\xi_1,\xi_2,\xi_i;\ q^{-1}\bar
  \xi_1,q\bar \xi_2,\bar \xi_i) - \\ - q\tau_1^{10}(q^{-1}\xi_1,\xi_2,\xi_i;\
  \bar \xi_1,\bar \xi_2,\bar \xi_i)
\tau_1^{01}(\xi_1,\xi_2,\xi_i;\
  q\bar \xi_1,\bar \xi_2,\bar \xi_i) + \\
+ q^2\tau_1^{11}(q^{-1}\xi_1,\xi_2,\xi_i;\
  \bar \xi_1,\bar \xi_2,\bar \xi_i)
\tau_1^{00}(\xi_1,\xi_2,\xi_i;\
  q^{-1}\bar \xi_1,q\bar \xi_2,\bar \xi_i)
\ee
This expression can be rewritten in a more compact form using
the operators
\be
{\cal D}_1^L \equiv M_1^-D_1\otimes I, \ \ \ {\cal D}_1^R
\equiv M_1^+M_2^-\otimes D_1, \ \ \ \bar{\cal D}_1^L \equiv
{\bar D_1} \otimes \bar M_1^-\bar M_2^+, \ \ \ \bar{\cal
D}_1^R \equiv I\otimes \bar M_1^+{\bar D_1}\\
\tau_2 = \left({\cal D}_1^R\bar{\cal D}_1^R -
q{\cal D}_1^L\bar{\cal D}_1^R - q{\cal D}_1^R\bar{\cal D}_1^L
+ q^2{\cal D}_1^L\bar{\cal D}_1^L\right)\tau_1\otimes
\tau_1
=\left( {\cal D}_1^R-q{\cal D}_1^L\right)\cdot\left(
\bar{\cal D}_1^R-q\bar{\cal D}_1^L\right)
\tau_1\otimes \tau_1
\ee
These operators satisfy the commutation relations that form the algebra
similar to that generated by $\theta$ and $\chi$ in (\ref{comrel}):
$
{\cal D}_1^{L}    {\cal D}_1^{R}
= q{\cal D}_1^{R}{\cal D}_1^{L}, \ \
\bar {\cal D}_1^{L} \bar   {\cal D}_1^{R}
= q\bar{\cal D}_1^{R}\bar{\cal D}_1^{L}$.
Formula for $\tau_2$ can be rewritten in a more ``invariant" form in terms
of the operators
\be
\hbox{{\hr D}}^L_i\equiv D_i\otimes I,\ \ \
\hbox{{\hr D}}^R_i\equiv \prod_j M_j^{-\vec\alpha_i\vec\alpha_j}\otimes D_i,
\ \ \
\bar {\hbox{{\hr D}}}^L_i\equiv \bar D_i\otimes\prod_j \bar
M_j^{-\vec\alpha_i\vec\alpha_j},\ \ \
\bar {\hbox{{\hr D}}}^R_i\equiv I\otimes \bar D_i
\ee
that commute as
\be
\hbox{{\hr D}}_i^{L}    \hbox{{\hr D}}_j^{R}
= q^{\vec\alpha_i\vec\alpha_j}
\hbox{{\hr D}}_j^{R}\hbox{{\hr D}}_i^{L}, \ \ \ \ \ \
\bar {\hbox{{\hr D}}}_i^{L} \bar {\hbox{{\hr D}}}_j^{R}
= q^{\vec\alpha_i\vec\alpha_j}\bar{\hbox{{\hr D}}}_j^{R}\bar{\hbox{{\hr
D}}}_i^{L}
\ee
Indeed, one can write
\be
\tau_2 = M^-_1\otimes\bar M_1^+\left( \hbox{{\hr D}}_1^R-q
\hbox{{\hr D}}_1^L\right)\cdot\left(
\bar{\hbox{{\hr D}}}_1^R-q\bar{\hbox{{\hr D}}}_1^L\right)
\tau_1\otimes \tau_1
\ee

\app{Construction of the Whittaker function for $SL(p)$}
In this Appendix, we construct the Toda WF corresponding to the group
$SL(p)$ \cite{GKMMMO}. We start with notations and definitions. The
standard textbooks for us are \cite{Wa}.

\sapp{Notations}
Lie algebra of the Lie group $SL(p)$ is completely given by the
following commutation relations of the generators corresponding to the
simple roots (Chevalle generators):
\be\label{commrel1}
\phantom{fhg}[T_{\pm i},T_{0,j}]=\mp A_{ij}T_{\pm,i},\ \ [T_{+i},T_{-j}]=
\delta_{ij}T_{0,j},\ \ i,j=1,\ldots, p-1
\ee
and additionally the Serre relations
\be\label{Serre}
\phantom{fhg}\hbox{ad}_{T_{\pm i}}^{1-A_{ij}}\left(T_{\pm j}\right)=0
\ee
where $\hbox{ad}_x^k(y)\equiv \underbrace{[x,[x,...,[x,y]..]]}_{k\ \hbox{
times}}$. All other commutation relations can be obtained from
(\ref{commrel1})-(\ref{Serre}), the generators which correspond to positive
(negative) nonsimple roots being constructed from the positive
(negative) simple root (Chevalle) generators by the manifest formula
$[T_{\balpha},T_{\bbeta}]=N_{\balpha,\bbeta}T_{\balpha+\bbeta}$.
Here the generator $T_{\balpha+\bbeta}$ corresponds to the non-simple root
$\balpha+\bbeta$ and $N_{\balpha,\bbeta}$ are some non-zero structure
constants. Making use of nonsimple root generators, the Serre identities
are replaced by appropriate Lie algebra relations.

We also use the following notations:
$\balpha_i$ are the simple roots in the corresponding Lie algebra,
$A_{ij}\equiv\balpha_i\cdot\balpha_j$ is the Cartan matrix,
$\bmu_i$ are the fundamental weights that, by definition, lie in the
dual lattice
$\bmu_i\cdot\alpha_j=\delta_{ij}$ (i.e. $\bmu_i= A_{ij}^{-1}\balpha_j$) and
$\brho\equiv \2\sum_{\balpha\in\De^+}\balpha=\sum_i \bmu_i$, where $\De^+$
is the set of all the positive roots. In these notations,
$\phi_i=-\bmu_i\cdot\bphi$. Roots $\balpha_i$ can be also considered as
vectors in the $p$-dimensional Cartan plane of the group
$GL(p)$:  $\balpha_i=\bfe_{i+1}-\bfe_i$, $\bfe_i\cdot\bfe_j=\delta_{ij}$.

In order to define the Casimir operators, let us fix a representation
$\rho$ of UEA $U(SL(p))$. Define
the $L$-operator by the formula \cite{FRT}
\be
L\equiv \sum_{\balpha\in\De}\rho(T_{\balpha})\otimes T_{-\balpha}+
\sum_iA^{-1}_{ij}\rho(T_i)\otimes T_j
\ee
Then the $k$-th Casimir operator can be defined via the formula
\be
C_k\equiv \Tr_\rho L^k
\ee
where the trace is taken over the representation $\rho$. Indeed, since
the result does not depend on the representation $\rho$, one can choose
it to be the simplest one -- the first fundamental representation. Then,
the Casimir operator can be easily calculated. In particular, the
quadratic Casimir acquires the form
\be
C_2=\sum_{\balpha\in\De}
T_{\balpha}T_{-\balpha}+\sum_{ij}^{p-1}A_{ij}^{-1}T_{0,i}T_{0,j}
\ee
where the first sum goes over all (positive and negative) roots.

\sapp{Representations}
For arbitrary group $SL(p)$, one can define the (right) regular
representation only in general terms of the group acting on the
algebra of functions:
\beq\label{rrrep}
\pi_{reg}(h) f(g)=f(gh)
\eeq
Hence, we mainly use the group (not algebraic) terms (see, however,
\cite{GKMMMO}). It turns out that, for the generic group $SL(p)$,
one can describe restricting the space of functions onto irreducible
representations. To this end, one needs to consider the representations
induced by one-dimensional representations of the Borel subgroup. This
restricts the space of all functions onto the functions satisfying the
covariantness condition
\beq\label{irrep}
f_{\lambda}(bg)=\chi _{\lambda}(b)f_{\lambda}(g)
\eeq
where $b$ is an element of the Borel subgroup of the lower-triangle
matrices and $\chi _{\lambda}$ is the character of the Borel subgroup,
which is equal to
\beq \label{rep}
\chi
_{\lambda}(b)=\prod _{i=1}^{p-1} \mid b_{ii} \mid ^{(\blambda-\brho)\bfe_i}
(\mbox{sign}b_{ii} )^{\epsilon _i}
\eeq
where $\epsilon_i$ are equal either to 0 or to 1. For the sake of
simplicity, we consider only the representations without sign factors,
although the generic case can be also easily treated. The described
representations belong to the principal sphere series.

Thus, the representations are given by restricting the space of functions
onto the functions defined on the quotient group
$B\backslash G$ that can be identified with the (nilpotent) subgroup
of the strictly upper-triangle matrices $N_+$.
Given $\lambda$, there is the natural Hermitian bilinear form on the
space of matrix elements $X$, which is just given by the flat measure
\be
<f_L|f_R>_{\lambda}
=\int_{X=B\backslash G} \overline {f_{L,\lambda}(x)}f_{R,\lambda}(x)
\prod_{ij}dx_{ij}
\ee
This form becomes the scalar product when $\lambda$ is pure imaginary, which
leads to unitary irreducible representations of the principal series.

\sapp{Hamiltonian and Schroedinger equation}
We consider the following reduction conditions analogous to
(\ref{statecondG}):
\be\label{rcglNGR}
T_{+,i}|\psi_R>=i\mu_i^R|\psi_R>
\ee
and
\be\label{rcglNGL}
<\psi_L|T_{-,i}=i\mu_i^L<\psi_L|
\ee
where $\mu^{R,L}_i$ are cosmological constants. These conditions are
completely given by the Chevalle generators, since all others are
easily obtained from the commutation relations. More precisely,
the action of nonsimple root generators just vanishes. Now one can
easily get
\be
(\blambda^2-\brho^2)F^{(\blambda)}(\bphi)\equiv
(\blambda^2-\brho^2)<\psi_L|e^{-\bmu_i\bphi T_{0,i}}|\psi_R>=
<\psi_L|e^{-\bmu_i\bphi T_{0,i}}C_2|\psi_R>=\\=
<\psi_L|e^{-\bmu_i\bphi T_{0,i}}\left(2\sum_{\balpha\in\De^+}
T_{\balpha}T_{-\balpha}+ 2\sum_{ij} A^{-1}_{ij}T_{0,j}
+\sum_{ij}A_{ij}^{-1}T_{0,i}T_{0,j}\right)|\psi_R>=\\=
\left( {\partial^2\over\partial\bphi^2}+ 2\sum_i{\d (\balpha_i\bphi)}
-2\sum_i\mu^L_i\mu_i^Re^{\balpha_i\bphi}
\right)<\psi_L|e^{-\bmu_i\bphi T_{0,i}}|\psi_R>
\ee
and $\Psi^{(\blambda)}(\bphi)=e^{-\brho\bphi}F^{(\blambda)}(\bphi)$
satisfies the Toda chain Schroedinger equation
\be\label{LeslNG}
\left( {\partial^2\over\partial\bphi^2}
-2\sum_i\mu^L_i\mu_i^Re^{\balpha_i\bphi}
\right)\Psi^{(\blambda)}(\bphi)=\blambda^2
\Psi^{(\blambda)}(\bphi)
\ee

\sapp{Solution to the Schroedinger equation}
In order to calculate the Whittaker function explicitly, one should
manifestly construct solutions to the conditions (\ref{rcglNGR})
and (\ref{rcglNGL}). However, we do not use here the algebraic approach,
i.e. do not solve the corresponding differential
equation\footnote{This is done in \cite{GKMMMO}.}, but, instead, use the
group theory approach. That is, we find the functions
$f_{L,\lambda}(g)$ and $f_{R,\lambda}(g)$ in the representation space,
which satisfy the conditions
\be\label{statecondRN}
\pi_{\lambda}(z) f_{R,\lambda}(g)=f_{R,\lambda}(gz)=e^{\sum
i\mu_iz_{i,i+1}}f_{R,\lambda}(g) =e^{iTr (\mu z)}f_{R,\lambda}(g),\\ (\mu
)_{ij}=\delta _{i-1,j}\mu _i, \ \ \ g,z\in N_+
\ee
and
\be\label{statecondLNG}
\pi_{\lambda}(z^t) f_{L,\lambda}(g)=e^{iTr (\mu_L z)}f_{L,\lambda}(g)
\ee
which are the group counterparts of the conditions
(\ref{rcglNGR}) and (\ref{rcglNGL}).

Solving the first condition is equivalent to constructing a one-dimensional
representation of the group of the upper-triangle matrices. At first,
we construct an additive character of the group using the fact that, in the
product of two upper-triangle matrices, the elements lying on the
next-to-main diagonal are summed. Then, one can exponentiate this
character and obtain the one-dimensional representation we are looking
for:
\beq\label{fR}
f_{R,\lambda}^{\mu}(x)=e^{i\Tr \mu x}
\eeq

In order to find the function $f_{L,\lambda}(x)$ that satisfies
\ref{statecondLNG}), we use the inner automorphism of the group
$SL(p)$, which maps strictly upper-triangle matrices into strictly
lower-triangle ones. This automorphism is explicitly described by
the matrix
\beq
S=\left ( \begin{array}{ccccccc} 0 & 0      & 0
& 0     & 0  & 0      & 1 \\ 0 & \ldots & 0      & 0     & 0  & 1      & 0 \\
0 & \ldots & 0      & 0     & 1 & 0      & 0 \\ 0 & \ldots & 0      & 1    &
0  & \ldots & 0 \\ & \ldots & 1 & 0     & 0  & \ldots &   \\ 0 & 1      & 0
&\ldots & 0  & \ldots & 0 \\ 1 & \ldots & 0      & 0     & 0  & \ldots & 0
\end{array} \right )\
\eeq
i.e.
\beq
S_{ij}=\delta _{i+j,p+1}
\eeq
and
\beq\label{tmat}
(S^{-1}zS)_{ij}=
z_{p+1-i,p+1-j}
\eeq
In fact, it is necessary for the matrix $S$ to be an element of the
group $SO(p)$ (generally it is constructed as an element of the Weyl group
\cite{GKMMMO}). Therefore, it should be normalized so that it has the
unit determinant. One can do this multiplying the matrix by -1, which
does not effect formula (\ref{tmat}). Note that the automorphism $S$ maps
the element $z_{i,i-1}$ into the element
$z_{p-i,p+1-i}$ instead of $z_{i-1,i}$ (because of the condition
(\ref{rcglNGL}), we are only interested in the elements $z$ lying on the
next-to-main diagonal). The situation can be corrected by the proper
reflection of the matrix $\mu$. Namely, one should introduce the new matrix
$\tilde\mu_L\equiv \mu_{p-i}^L\delta_{i-1,j}$ so that the condition
(\ref{statecondLNG}) can be rewritten as
\be\label{newcond}
\pi_{\lambda}(SzS^{-1})f_L(g)=f_L(gSzS^{-1})
e^{iTr (\tilde\mu_L z)}f_{L,\lambda}(gS^{-1}S)
\ee
i.e.
\be
f_L(gS^{-1}z)=
e^{iTr (\tilde\mu_L z)}f_{L,\lambda}(gS^{-1})
\ee
We already know that the solution to this equation in upper-triangle
matrices is
\be
\left.f_L(gS^{-1})\right|_{B_-=0}=f_L(n_+)=e^{iTr (\tilde\mu_L n_+)}
\ee
where $n_+$ the upper-triangle part of $gS^{-1}$. We need, however,
the solution of (\ref{newcond})
for $g=x\in N_+$, i.e. some re-calculation is necessary:
\be
f_L(xS^{-1})=\chi_{\lambda}(xS^{-1})f_L(n_+)=\chi_{\lambda}(xS^{-1})
e^{iTr (\tilde\mu_L n_+)}
\ee
Let us calculate this function more explicitly. First of all, elements of the
diagonal matrix $h$ in the Gauss decomposition of any matrix $g$ are given by
the formula\footnote{The simplest way to derive this formula is to note that
the diagonal element $g_{ii}$ of the matrix $g$ depends only on
matrix elements lying in the left upper corner of the $i\times i$
submatrix. Then, one may use the induction over the rank of matrix to obtain
formula (\ref{gaussexphr}).}
\be\label{gaussexphr}
h_i={\Delta_i(g)\over\Delta_{i-1}(g)},\  \ \ h_1=\Delta_1(g),\ \ \
\Delta_p(g)= {\Delta_1(g)\over\Delta_{p-1}(g)}
\ee
where $\Delta_i(g)$ denotes the upper main $i\times i$ minor of the
matrix $g$. Thus, the diagonal part of $xS^{-1}$
is given by
the ratio $h_i={\Delta_i(xS^{-1})\over \Delta_{i-1}(xS^{-1})}$ and
\be\label{2}
\chi_{\lambda}(xS^{-1})=\prod_i\Delta_i^{(\blambda\balpha_i-1)}(xS^{-1})
\ee
In order to get the formula for the elements
$(n_+)_{k-1,k}$, note that they depend only on the
$k\times k$ submatrix placed in the upper left corner.
Consider the element $(k-1,k)$ of the matrix $b^{-1}\equiv n_+g^{-1}$,
$g=xS^{-1}$. It is equal to zero by definition of $b$ and, therefore,
one obtains
\beq
(n_+)_{k-1,k}(g^{-1})_{k,k}+(n_+)_{k-1,k-1}(g^{-1})_{k-1,k}=0
\eeq
Using the manifest representation for elements of the inverse matrix,
one gets the result:
\beq
(n_+)_{k-1,k}=\frac {\Delta_{k-1,k}(xS^{-1})}{\Delta_{k-1}(xS^{-1})}
\eeq
where $\Delta_{k-1,k}(xS^{-1})$ is defined to be the determinant of
the $(k-1)\times (k-1)$ submatrix of the matrix $xS^{-1}$ with
the $k$-th and $k-1$-th columns exchanged.

Thus, we have the manifest expressions for the functions
$f_L$ and $f_R$. In order to calculate the Whittaker function, one
only needs now to determine the action of the Cartan part of the
group element $g$ on $f_R(x)$. Note that, although the element
$xh$ with $x\in X$ and $h\in H$ does not belong to $X$, the element
$h^{-1}xh$ does. Hence, using (\ref{rrrep}), one can obtain
\be\label{Cartan}
\pi_{\lambda}(h)f_{R,\lambda}(x)=
f_{R,\lambda}(hh^{-1}xh)=h^{\blambda-\brho}f_{R,\lambda}(h^{-1}xh)
\ee
Then, using this formula and the manifest expression (\ref{fR}), we
get
\be
\pi_{\lambda}\left(e^{-\bmu_i\bphi T_{0,i}}\right)f_{R,\lambda}(x)=
e^{\left\{-\bmu_i\cdot\bphi\right\}\left\{(\blambda-\brho)\cdot
(\bfe_{i+1}-\bfe_i)\right\}}f_{R,\lambda}
(e^{\bmu_i\bphi T_{0,i}} xe^{-\bmu_i\bphi T_{0,i}})=\\=
e^{(\brho-\blambda)\bphi}f_{R,\lambda}(e^{\bmu_i\bphi T_{0,i}}
xe^{-\bmu_i\bphi T_{0,i}})=
e^{(\brho-\blambda)\bphi}e^{i\Tr x\mu
e^{\balpha\bphi}}
\ee
where the combination $\mu e^{\balpha\bphi}$ denotes the matrix
with the matrix elements
$\delta_{i-1,j}\mu_ie^{\balpha_i\bphi}$. Collecting all this together,
we finally obtain for the Whittaker function (\ref{res}):
\be
\Psi (\phi)=e^{-\blambda\bphi}
\int _{X=B\backslash G} \prod_{i<j}dx_{ij}\prod_{i=1}^{p-1}
\Delta_i^{-(\blambda\balpha_i+1)}(xS^{-1})\times
e^{ i\mu_{i}^Rx_{i,i+1}e^{\balpha_i\bphi} -i\mu_{p-i}^L
\frac{\Delta _{i,i+1}(xS^{-1})}{\Delta _{i}(xS^{-1})}}
\ee

\end{document}